\documentclass[a4paper,11pt]{amsart}
\usepackage{amscd}
\usepackage{amssymb}
\let\:=\colon

\def\CC{{\mathbf C}}
\def\DD{{\mathbf D}}
\def\FF{{\mathbf F}}

\def\NN{{\mathbf N}}
\def\PP{{\mathbf P}}
\def\QQ{{\mathbf Q}}

\def\ZZ{{\mathbf Z}}
\def\A{{\mathcal A}}
\def\B{{\mathcal B}}
\def\C{{\mathcal C}}
\def\D{{\mathcal D}}
\def\E{{\mathcal E}}
\def\F{{\mathcal F}}
\def\G{{\mathcal G}}

\def\I{{\mathcal I}}

\def\L{{\mathcal L}}

\newcommand{\oh}{{\mathcal{O}}}
\def\P{{\mathcal P}}
\def\Q{{\mathcal Q}}

\def\U{{\mathcal U}}
\def\V{{\mathcal V}}

\newcommand{\bbD}{\mathbb{D}}

\def\x{\times}
\def\*{\otimes}

\def\sub{\subseteq}
\def\Aut{\operatorname{Aut}}

\def\Proj{\operatorname{Proj}}

\def\Hom{\operatorname{Hom}}
\def\HomC{{\mathcal Hom}}

\def\CEnd{{\mathcal End}}
\def\Hilb{\operatorname{Hilb}}
\def\Pic{\operatorname{Pic}}
\def\Grass{\operatorname{Grass}}

\def\rank{\operatorname{rank}}

\def\top{\operatorname{top}}
\newcommand{\kker}{\ker\,\,}
\def\Isom{\operatorname{Isom}}

\def\coker{\operatorname{coker}}

\def\codim{\operatorname{codim}}
\def\ed{\operatorname{ed}}
\def\rank{\operatorname{rank}}
\def\Sym{\operatorname{Sym}}
\newcommand{\barr}{\overline}
\newcommand{\rarr}{\rightarrow}
\def\ttheref{
\mathop{\hbox{\boldmath$\cdot$}\raise1.0ex\hbox{\boldmath$\cdot$}
\hbox{\boldmath$\cdot$}
}\nolimits}

\def\tttheref
{\mathop{\hbox{\boldmath$\cdot$}
\raise2.1ex\hbox{\kern0.4em\hbox{\boldmath$\cdot$}}
\hbox{\kern0.4em\hbox{\boldmath$\cdot$}}}
\nolimits}

\def\doublerightarrow{
\mathop{\hbox{$\rightarrow$}
\kern-0.8em\hbox{$\rightarrow$}}\nolimits}

\def\doubleleftarrow{
\mathop{\hbox{$\leftarrow$}
\kern-0.8em\hbox{$\leftarrow$}}\nolimits}

\def\righttheref{
\mathop{\hbox{\boldmath$\cdot$}\raise1.0ex\hbox{\boldmath$\cdot$}
\kern-0.6em\hbox{$\rightarrow$}
}\nolimits}

\def\lefttheref{
\mathop{\hbox{$\leftarrow$}
\kern-0.6em\raise1.0ex\hbox{\boldmath$\cdot$}
\hbox{\boldmath$\cdot$}
}\nolimits}

\def\rrighttheref{
\mathop{\hbox{\boldmath$\cdot$}
\raise2.1ex\hbox{\kern0.5em\hbox{\boldmath$\cdot$}}
\kern-1.0em\hbox{\boldmath$\longrightarrow$}
}\nolimits}

\def\llefttheref{
\mathop{\hbox{$\boldmath\longleftarrow$}
\kern-1.1em\raise2.1ex\hbox{\boldmath$\cdot$}
\hbox{\kern0.5em\hbox{\boldmath$\cdot$}}}\nolimits}

\def\uptheref{
\mathop{\hbox{\boldmath$\cdot$}
\kern-0.2em\raise1.0ex\hbox{\boldmath$\nearrow$}
\kern0.2em\hbox{\boldmath$\cdot$}
}\nolimits}

\def\otimesover#1{\underset{#1}{\otimes}}
\def\coprodover#1{\underset{#1}{\coprod}}

\def\Sim{\,\,\raisebox{0.3ex}{$\frown$}\!
                \raisebox{-0.2ex}{$\smile$}\,\,}
\def\smallhbar{
\mathop{\hbox{$\scriptstyle h$}
\kern-0.25em\lower0.1ex\hbox{$\scriptstyle\Bar{}\,\,$}
}\nolimits}

\newcommand{\map}{\barr{\mathcal{M}}}
\newcommand{\smap}{\barr{M}}
\newcommand{\M}{\barr{M}}

\newcommand{\sto}{_{0,n}(X,\beta)}

\newcommand{\nk}{\barr{M}_{0,n}}

\newcommand{\moo}{\barr{M}_{0,0}}

\newcommand{\gotha}{\mathfrak{a}}

\newcommand{\smallcup}{\mathbin{\text{\scriptsize$\cup$}}}

\def\therefore
{\mathop{\hbox{$\cdot$}\raise1.0ex\hbox{$\cdot$}\hbox{$\cdot$}}
\nolimits}

\def\ttherefore{
\mathop{\hbox{\boldmath$\cdot$}\raise1.0ex\hbox{\boldmath$\cdot$}
\hbox{\boldmath$\cdot$}
}\nolimits}

\def\Therefore{\mathop{
\hbox{$\bullet$}\raise1.5ex\hbox{$\bullet$}\hbox{$\bullet$}}
\nolimits}
\newcommand{\largefrac}{\!\mbox{\large/}}
\newcommand{\mmid}{\!\mid\!}
\newcommand{\lav}{<\!}
\newcommand{\rav}{\!>}
\newcommand{\laav}{<\!\!<\!}
\newcommand{\raav}{\!>\!\!>}

\newtheorem{thm}{Theorem}[section]

\newtheorem{lem}[thm]{Lemma}
\newtheorem{prop}[thm]{Proposition}

\newtheorem{ex}{Example}

\theoremstyle{definition}

\theoremstyle{remark}

\numberwithin{equation}{section}

\begin{document}
\thispagestyle{empty}
\title[Rational curves on $N(3;2,3)$]
{Rational curves on the space of \\
determinantal nets of conics}
\author{Erik N. Tj\o tta}
\address{Institute of Mathematics\\ University of Bergen\\ 5007 Bergen,
 Norway}
\email{erik.tjotta@mi.uib.no}
\date{September 29, 1997.}
\begin{abstract}
We describe the Hilbert scheme components para\-me\-tri\-zing
lines and conics on the space of determinantal nets of
conics, $\NN$. As an application, we use the quantum Lefschetz
hyperplane principle to compute the instanton
numbers of rational curves on a complete intersection Calabi-Yau
threefold in $\NN$.  
We also compute the number of lines and
conics on some Calabi-Yau sections of non-decomposable
vector bundles on $\NN$. The paper contains a brief summary of
the A-model theory leading up to Givental-Kim's quantum Lefschetz 
hyperplane
principle. This work makes up my doctoral dissertation, defended 
Oktober 31, 1997.
\end{abstract}
\maketitle
\clearpage
\pagenumbering{roman}
\tableofcontents
\listoftables
\clearpage
\addtocontents{toc}{\vspace{0.2cm}}
\section*{Introduction}
\setcounter{page}{-1}
\pagenumbering{arabic}
\label{O-intro}
\bigskip

In recent years, observations in theoretical physics have
directed enormous attention to problems related to the
enumerative geometry of algebraic curves on manifolds. Ideas
 from string theory and topological field theory have led to
the development of quantum cohomology, and to the mirror
symmetry conjecture.

To a manifold $M$ that satisfies certain positivity
assumptions, one may associate a topological field theory
$A(M)$, where the correlation functions carry enumerative
information (in the form of Gromov-Witten (GW)-invariants)
about morphisms of algebraic curves to the manifold. In
\cite{W,KM} it was shown that this puts severe
restrictions on the GW-invariants, sometimes so strong that
all GW-invariants (at least in the genus zero case) can be
determined from a few simple ones. Subsequently, a big
effort was made to build a mathematical framework for this
formalism. This has now been completed within the realms of
symplective geometry \cite{RT}, and algebraic geometry
\cite{KM,BM,BF,B,LiT}.
This theory has provided enumerative geometers
with a whole new tool. An example of its power is the
celebrated recursion formula obtained by Kontsevich for the
number of rational curves of degree $d$ in $\PP^2$ meeting
$3d-1$ given points (See Example~\ref{exPP2} in
Section~\ref{potential}).

In string theory, the assumption is that space-time is a
ten-dimensional manifold which is locally a product
$M_{4}\times M$, where $M_{4}$ is the usual Minkowski
space-time, and $M$ is a Calabi-Yau complex threefold.
Initially, any Calabi-Yau is eligible. By a formal dimension
count, algebraic curves on $M$ are expected to be rigid, and
the $A$-model correlation functions should carry invariants
which reflect the number of curves of given genus in a
given cohomology class. Accordingly, the study of algebraic
curves on Calabi-Yau manifolds is a very natural and
interesting mathematical problem of ``physical'' interest. A
further feature which enhances the appeal of these
manifolds, is the following intriguing physical observation: To a
Calabi-Yau manifold $M$ one can associate two different
topological field theories $A(M)$ and $B(M)$, in such a
way that there should exist another Calabi-Yau manifold $W$
with $A(M)=B(W)$ and $A(W)=B(M)$.
Mathematically, it indicates an unexpected relationship
between the symplectic geometry of $M$ and the complex
geometry of $W$, and vice versa. This is the mirror
symmetry conjecture.

Interest into this field exploded with the advent of a
sensational computation by a group of physicists \cite{CdOGP}.
In their paper, the mirror symmetry conjecture was used to compute a
generating function for the numbers $n_{d}$ of rational curves of
degree $d$ on the quintic threefold in $\PP^4$. The
essential outcome of the computation was that two
differential operators
\begin{equation*}
D^4-5q(5D+1)(5D+2)(5D+3)(5D+4)\qquad{\rm and}\qquad
D^2\frac{1}{K(q)}D^2\,,
\end{equation*}
with $K(q)=5+\sum_{d>0}n_{d}d^3\frac{q^d}{1-q^d}$ and
$D=q\frac{\partial}{\partial q}$, were equivalent up to a
transformation (See Example~\ref{exU} in
Section~\ref{GQH} for a more precise statement). This enabled
the (conjectural) determination of $n_{d}$
for arbitrary large $d$. A result of
this order was far beyond reach of conventional
mathematical methods. Up to then, mathematicians had only
been able to compute the numbers $n_{1}$ and $n_{2}$, and soon
afterwards $n_{3}$ and
$n_{4}$. These computations confirmed the conjectured numbers.
A crucial point in the computation of \cite{CdOGP} was that a
mirror candidate of the quintic had been identified. Simply
stated: A difficult computation related to the $A$-model of
the quintic, became simple once it could be translated to
the $B$-model of its mirror. Soon afterwards, the
computations pioneered in \cite{CdOGP} were further developed and
successfully applied to a host of new cases where mirror
candidates had been identified, e.g., \cite{M1,LT,GMP,BvS,Me}.
One limitation to this procedure is that only for a relatively small
class\footnotemark
\footnotetext{Most notably, complete intersections in toric
varieties \cite{BB}, but see also \cite{Vo1,Bo,R,BCKvS}
for other mirror candidates and recent development.}
of Calabi-Yau manifolds does one know of good mirror candidates.
However, by studying some combinatorial aspects of mirror
computations, it was essentially conjectured in \cite{BvS},
that in some cases the $A$-model contains enough structure
to explain the computations. In a series of remarkable
papers \cite{G}, Givental proved the above mirror
computations by finding a way to perform them rigorously on
the $A$-model as indicated by \cite{BvS}. The core of the
proof is a quantum Lefschetz hyperplane principle, which
relates the quantum cohomology of a manifold to the quantum
cohomology of its complete intersections. For a deeper
background into the subject of mirror symmetry, the reader
is recommended \cite{M2,Vo2,CK}.

\bigskip

In this paper, we study rational curves on the
six-dimensional variety $\NN$ which para\-me\-tri\-zes
determinantal nets of conics. Alternatively, it may be
viewed as the variety para\-me\-tri\-zing trisecant planes of the
Veronese surface in $\PP^5$. Our main reason for this study
is that there are three vector bundles on $\NN$ whose
general global sections are Calabi-Yau threefolds. In fact,
one of the bundles is decomposable, hence the associated
Calabi-Yau section is a complete intersection in $\NN$. As
usual, the problem of computing (by conventional methods)
rational curves on the sections leads to the problem of
describing the curves in $\NN$. We shall describe the Hilbert
schemes of lines and conics on $\NN$, and use this to
compute the number of lines and conics on the Calabi-Yaus.
The highlight of the paper is an application of the quantum
Lefschetz hyperplane principle to the computation of the
instanton numbers of rational curves on
the complete intersection. A generating function for these
numbers will be obtained from a simple knowledge of the
cohomology of $\NN$ and the lines on $\NN$. Again, this
illustrates the power of quantum cohomology. To the best of
our knowledge, this is the first ``mirror computation''
performed without knowing the mirror.

The paper is organized as follows:
Sections~\ref{modulispaces} and \ref{GQH} contain a basic
review of the theory of quantum cohomology following
\cite{KM,BM,G,Kim}. No originality is claimed
for this, except for perhaps a few odd parts in
Section~\ref{gqs}. The cohomology of $\NN$ has essentially
been computed in \cite{ES2}. We give a review of this
computation in Section~\ref{trisecant}. In
Section~\ref{lines}, we describe the Hilbert scheme of lines
on $\NN$ as a projective bundle over
$\PP^2\times\PP^{2\vee}$. This enables us to describe the
Chow ring of the Hilbert space of lines. As an application we
compute the number of lines on the Calabi-Yau sections. The instantons for
the complete
intersection Calabi-Yau are computed in Section~\ref{qdmN}.
In Section~\ref{conics}, we give a description of the
Hilbert scheme of conics which suffices for computational
purposes. This is used to compute the (virtual) number of
conics on the Calabi-Yau sections. We have collected some
technical aspects of the computations in
Section~\ref{calculations} and in the appendices.

\bigskip

\noindent{\bf Preliminaries.}
An important reason for the success of \cite{G}
is the use of a localization principle in equivariant
cohomology \cite{AB}. The idea to use such a method in the
context of enumerative geometry first appeared in
\cite{ES3}. In our work, this principle will also play a
role. Below, we state a simple version that is suitable for
our purpose.

Let $M$ be a smooth orbifold with a $\CC^{\ast}$-action
whose fixpoints are isolated. Let $M^{\CC^{\ast}}$ denote
the set of fixpoints. Let $\B_{1},\ldots,\B_{s}$ be
$\CC^{\ast}$-bundles on $M$ of rank $r_{1},\ldots,r_{s}$.
If $f\in M$ is a fixpoint, there is an induced orbifold
$\CC^{\ast}$-action action on the fibers $\B_{if}$. Let
$\omega_{1}(f),\ldots,\omega_{m}(f)$ and
$\tau_{1}(\B_{i},f),\ldots,\tau_{r_{i}}(\B_{i},f)$ be the
weights of the $\CC^{\ast}$-action on $T_{f}M$ and $\B_{if}$
respectively. Since $M$ is an orbifold we must allow the
weights to be in $\QQ$.

Suppose $P(\ldots,c_{j}(\B_{i}),\ldots)$ is a polynomial in
the Chern classes $c_{j}(\B_{i})$. Since the Chern classes
$c_{j}(\B_{i})$ are symmetric polynomials in the Chern
roots $b_{i,1},\ldots,b_{i,r_{i}}$ of $\B_{i}$, we can find
a polynomial $\Tilde{P}$ such that
\\
\[P(\ldots,c_{j}(\B_{i}),\ldots)=
\Tilde{P}(\ldots,b_{i,1},\ldots,b_{i,r_{i}},\ldots)\,.\]

\medskip

\noindent{\bf Theorem.}
(Bott's formula)
\begin{equation*}
\int_{M}P(\ldots,c_{j}(\B_{i}),\ldots)=
\frac{1}{|\Aut(f)|}
\sum_{f\in M^{{\CC}^{\ast}}}
\frac{\Tilde{P}(\ldots,
\tau_{1}(\B_{i},f),\ldots,\tau_{r_{i}}(\B_{i},f),\ldots)}
{\omega_{1}(f)\ldots\omega_{m}(f)}\,.
\end{equation*}

\medskip

Another result of the same flavor is the classical
Bialynicki-Birula theorem \cite{Bi}. Let $M$ be a manifold
with a $\CC^{\ast}$-action as above. Assume
$M^{\CC^{\ast}}=\{f_{1},\ldots,f_{r}\}$,
and for each $i$, $1\le i\le p$, let $M_{i}=\{m\in
M|\lim_{t\mapsto 0}tm=m_{i}\}$. Consider the induced action
on $T_{f_{i}}M$, and let $T_{f_{i}}M^{+}$ be the subspace
where $\CC^{\ast}$ acts with positive weights.

\medskip

\noindent{\bf Theorem.} (Bialynicki-Birula)
{\em The $M_{i}$ form a cellular decomposition of $M$. Further,
each $M_{i}$ is an affine space with $\dim M_{i}=\dim
T_{f_{i}}M^{+}$.
}

\medskip

 From \cite{F} (Example~1.9.1, 19.1.11) it follows that
homological and rational equivalence on $M$ coincide, and that
the
Chow ring is the free abelian group generated by the classes
of $\smap_{i}$.

\bigskip

\noindent {\bf Notation.}
In this paper, we only work with even cohomology.
For a complex variety $X$ we use the notation
$A^k(X)=H^{2k}(X,\ZZ)$. If $R$ is a $\ZZ$-algebra, we let
$A_{R}^k(X)=A^k(X)\otimes R$. In many instances our
varieties will be smooth with a cellular decomposition,
hence $A^{\ast}(X)$ corresponds to the Chow ring. If $V$ is
a closed subvariety of $X$, then $V$ determines a cohomology
class in $A^{\ast}(X)$, which we also denote by $V$. The
bracket notation $[V]$ will be reserved for virtual
fundamental classes. If $X$ is smooth of dimension $n$, we
may identify cohomology classes with homology classes via
the Poincar\'e duality isomorphism
$A_{\QQ}^k(X)\to A_{n-k}^{\QQ}(X)$. This will be done
without warning throughout the text. Hence, if
$\alpha,\beta\in A^{\ast}(X)$, we use both
$\lav\alpha,\beta\rav$ and $\int_{\beta}\alpha$ to denote
the Poincar\'e pairing.

The notational conventions used are mostly as in \cite{H}.
If $\F$ is a sheaf on a scheme $X$, we write $H^i(\F)$ for
the $i$-th cohomology group unless there is an ambiguity. If
$f\:Y\to X$ is a morphism of schemes, we often use $\F_{Y}$
to denote the pullback $f^{\ast}\F$. For a
coherent sheaf $\F$, we use Grothendieck's $\PP(\F)=\Proj
(\Sym\F)$, and $\Grass_{r}(\F)$ denotes the Grassmanian of
rank $r$ subbundles. If $\L$ is a line with first Chern
class $\tau$, we shall often use the convenient notation
$\L=\oh(\tau)$.

\bigskip

\noindent {\bf Acknowledgments.} It is a pleasure to thank
I.~Ciocan-Fontanine, G.~Ellingsrud, G.~Fl\o y\-stad,
T.~Graber, P.~Meurer, K.~Ranestad, E.~R\o dland, and
D.~van Straten for valuable conversations on the subject
matter. A special thanks is due to B.~Kim for explaining me
the work of Givental and himself. More than anyone, I am
grateful to my supervisor S.~A.~Str\o mme, for introducing
me to a stimulating problem, and for always being available
for discussions.

A major part of this work came together during a
tremendously active year at the Mittag-Leffler Institute. I
would like to thank the Institute, and in particular
D.~Laksov, for support and for making the stay possible.

I also benefited much from a stay at Institut Henri
Poincar\'e. I am thankful to C.~Peskine for this invitation.

Finally, I would like to thank the Department of Mathematics
at the University of Bergen for a friendly working
environment.

This work has been supported by the Norwegian Research
Council for Science and Humanities, and by Nordisk
Forskerutdanningsakademi (NorFA).
\addtocontents{toc}{\vspace{0.2cm}}
\section{Moduli spaces of stable curves and maps}
\label{modulispaces}

\bigskip

The theory of quantum cohomology requires a good
compactification of the space of maps of curves into a
variety. One such compactification, at least when the
target variety is nice enough, is the moduli space of stable
maps introduced by Kontsevich \cite{K}. The moduli space of
stable maps specializes to the well known moduli space of curves
$\smap_{g,n}$ if the target variety is simply a point. In
this section we define, and state some relevant results
about, these moduli spaces. We only consider the genus zero
situation, but most results are true for higher genera.
For a nice introduction to this topic, and further details,
the reader is recommended \cite{FP}.\\

In Sections~\ref{modulispaces} and \ref{GQH}, $X$ will
always denote a smooth complex projective variety.

\bigskip

\subsection{Definitions}\label{definitions}

An $n$-pointed {\em stable curve}  of genus 0 $(C, s_{1}, \ldots, s_{n})$ is a
projective, connected, (at worst) nodal algebraic curve $C$ with $n$ marked
points
$\  s_{1}, \ldots, s_{n}$ such that:

\begin{enumerate}
\item[i)]$\dim H^1(C,\oh_{C})$=0
\item[ii)] $  s_{1}, \ldots, s_n$ are non-singular points on $C$
\item[iii)] each component of $C$ contains at least three special points,
where {\em
special} means singular or marked.
\end{enumerate}

A {\em stable map}
$(C,f, s_{1},\ldots,s_{n})$ from an $n$-marked genus 0 curve
$C$ to $X$ is a morphism $f\colon C \rarr X$ such that all components of
$C$ which are
mapped to a point are stable. Let $S$ be an algebraic scheme over
$\CC$. A {\em family of stable maps}
$(\C\to S,f,s_{1},\ldots,s_{n})$ of $n$-marked
genus $0$ curves to $X$ consists of a flat projective
morphism $\C\rarr S$ with sections $s_{i}$ and a morphism
$f\colon \C \to X$ such that each geometric fiber
$(\C_{s},f_{s},s_{1}(s),\ldots,s_{n}(s))$ is a stable map.

Two families of maps over $S$,
\[(\C\to S,f,s_{1},\ldots,s_{n})\qquad{\rm and}
\qquad(\C'\to S,f',s_{1}^{\prime},\ldots,s_{n}^{\prime})\]
are isomorphic if there exists an $S$-isomorphism $\varphi\colon
\C\rarr\C'$ such that
for each $i$ the following diagram commutes:
\[
 \begin{array}{ccc}
          \C
         &\xrightarrow[\!\!\displaystyle{
         \lower3.0ex\hbox{$
         {\scriptstyle s_{i}}\hbox{
\kern-0.7em \raise1.0ex\hbox{$\nwarrow$}}
         \lower2.0ex\hbox{$\quad
         {\!\!S\!\!}\hbox{\kern-0.3em\raise6.0ex\hbox{$\scriptstyle \varphi$}}
         \quad$}
         \hbox{\kern-0.1em \raise0.7ex\hbox{$
         \nearrow
\hbox{\kern-0.7em \lower1.0ex\hbox{$\scriptstyle s_{i}^{\prime}$}}
                                            $}}
         $}}\!\!]
         {\displaystyle{
         \hbox{$
         \nearrow
         \raise2.5ex\hbox{$\quad
         {\!\!X\!\!\!\!\!\!\!}\hbox{\kern-0.6em\lower5.0ex\hbox{$\Sim$}}
         \,\,$}
         \kern-0.5em\nwarrow
         $}}}
         &\C'\\
          &&
         \end{array}
\]
Let $\beta\in A_{1}(X)$. Define a contravariant functor:
\[\map\sto\colon {\underline{\rm{Sch}}}^0 \rarr
{\underline{\rm{Sets}}}\]
by letting $\map\sto(S)$ be the set of isomorphism classes of families
of stable maps $(\C\to S,f,s_{1},\ldots,s_{n})$ such that
$f_{\ast}[\C_{s}]=\beta$ for all geometric fiber $\C_s$.
Whenever $X$ is a point, we shall simply
write $\map_{0,n}$.

\medskip

\noindent{\em Remark~1.} The stability conditions i)--iii)
ensure that the automorphism groups of the maps
$(C,f ,s_{1},\ldots,s_{n})$ are finite.

\medskip

\noindent{\em Remark~2.} We shall allow labelling of marked
points by other indexing sets than $\{1,\ldots,n\}$, and we
will write $\map_{0,A}(X,\beta)$ if the indexing set is $A$.

\bigskip

\subsection{Natural maps}\label{natural}

The functor of stable maps comes equipped with the
following natural structures:
\\\\
{\em Evaluation.}
For $1\le i\le n$ define morphisms
\[e_{i}\colon \map\sto\rarr \HomC(\,-\,, X)\]
by
\[e_{i}(C,f,s_{1},\ldots,s_n)=f(s_{i})\, .\]
\\
{\em Forgetting point (Contraction).}
For $1\le i\le n+1$ define morphisms
\[\pi_{i}\colon \map_{0,n+1}(X,\beta)\to\map\sto\]
by forgetting the $i$-th marking and contracting eventual nonstable
components of the resulting map.
\\\\
{\em Forgetting the map.}
Define a morphism
\[\eta\colon \map\sto\to\map_{0,n}\]
by forgetting the map and contracting eventual
nonstable components of the resulting marked curve.
\\\\
{\em Gluing.}
Define a morphism
\[\varphi_{AB}^{\beta'\beta''}\colon \map_{0,A\cup\{\bullet\}}
(X,\beta')\x_X\map_{0,B\cup\{\bullet\}}(X,\beta'')\to
\map_{0,A\cup B}(X,\beta'+\beta'')\]
by gluing the maps in $\bullet$.

\bigskip

\subsection{Representability}\label{representability}

The moduli problem of stable curves and maps has been
studied in detail in \cite{DM,Kn,K}. Unfortunately, the
functor $\map\sto$ is in general not representable.

\begin{thm}
There exists a projective, coarse moduli space $\smap\sto$ for
the functor $\map\sto$.
\end{thm}

This will suffice for most of what we are going to consider.
However, we shall occasionally need a universal object for
the moduli problem, something that forces us to work in the category
of stacks \cite{DM}. A satisfactory intersection theory for
Deligne-Mumford stacks has been developed in \cite{V}.

\begin{thm}
The moduli stack of stable maps
$\smap\sto$ is a proper, separated Deligne-Mumford
stack of finite type with universal family
$(\smap_{0,n+1}(X,\beta),e_{n+1},s_{1},\ldots,s_{n})$

\begin{equation}
        \begin{CD}
                {\smap_{0,n+1}(X,\beta)} @>{e_{n+1}}>>X \\
                @V{\pi}VV \\
                   {\smap\sto}
        \end{CD}
\end{equation}
where $\pi$ forgets the $(n+1)$th point, and the
sections $s_{i}$, $i=1,\ldots,n$, are defined by the map
\[s_i\colon\smap\sto\simeq\smap\sto\x\smap_{0,3}\to\smap_{0,n+1}
(X,\beta)\,,\]
which glues along the $i$-th point of first factor.
\end{thm}

Note that we abuse notation and write $\smap\sto$ regardless of
whether we are working with it as a stack, or as a scheme.

\medskip

\noindent{\em Remark~1.} By deformation arguments (See \cite{FP}),
if $\smap\sto$ is non-empty, then
\[\dim\smap\sto\ge \dim X+\lav c_{1}(TX),\beta\rav +n-3\,.\]
The right-hand side of the equation is called the expected
dimension, which we will denote $\ed_\beta$.
If $\dim H^1(C,f^{\ast}TX)=0$ for all points
$[C,f,s_1,\ldots,s_n]$ in $\smap\sto$, then $\smap\sto$ is a
smooth Deligne-Mumford stack (i.e., an orbifold) of expected
dimension. Or, as a scheme, $\smap\sto$ has only finite
quotient singularities.

\medskip

\noindent{\em Remark~2.} If $X$ is a point and $n\ge 3$, we use
$\smap_{0,n}$ to
denote the moduli space instead of $\smap_{0,n}(pt,0)$. If
$\beta=0$, then $\smap_{0,n}(X,0)=\smap_{0,n}\times X$. When
$n\le 2$, it will be convenient for us to let $\smap_{0,n}$ be a
point, even if the expected dimension is negative.
For $n\ge 3$, $\smap_{0,n}$ is smooth of dimension $n-3$.

\bigskip

\subsection{Boundary divisors}
\label{bound.divisors}

In this section, we describe some divisors on the moduli
space of stable curves. Relations among these divisors are
the core of quantum cohomology.

Let $A\cup B$ be a partition of $\{1,\ldots,n\}$. The image
of the gluing maps $\smap_{0,A\cup\{\bullet\}}\times
\smap_{0,B\cup\{\bullet\}}\to\smap_{0,A\cup B}$ define
divisors which we denote $D(A|B)$. If
$\smap_{0,n}\to\smap_{0,4}$ is the contraction map, then the
inverse image of $D(12|34)$ is the divisor $\sum D(A|B)$,
where the sum is over all partitions $A\cup B$ of
$\{1,\ldots,n\}$ such that $1,2\in A$ and $3,4\in B$.
Since $\smap_{0,4}\simeq\PP^1$, the divisors $D(12|34)$
and $D(13|24)$ are linearly equivalent, hence we have the
fundamental relations
\begin{equation}
\label{1-1}
\sum_{\substack{1,2\in A\\ 3,4\in B}}D(A|B)=
\sum_{\substack{1,3\in A\\ 2,4\in B}}D(A|B)
\end{equation}
in the Chow ring of $\smap_{0,n}$. In fact, this 
describes the Chow ring of $\smap_{0,n}$:

\begin{thm}\cite{Ke}
The Chow ring $A^{\ast}(\smap_{0,n})$ is generated by the
divisors $D(A|B)$. The ideal of relations among the divisors
is generated by the relations \eqref{1-1} and the relations
$D(A|B)\cdot D(C|D)=0$, for which there are no inclusions among
the sets $A$, $B$, $C$, $D$.
\end{thm}
\addtocontents{toc}{\vspace{0.2cm}}
\section{Gravitational quantum cohomology}
\label{GQH}

\bigskip

In this section we give a rapid review of some basics of the
quantum theory found in \cite{W,KM,D,G}.
 Our main focus is to explain the
quantum Lefschetz hyperplane principle. In
Sections~\ref{gravdesc}--\ref{Dubrovin} we explain some
aspects of the WDVV-theory for gravitational descendents,
which are correlators generalizing the Gromov-Witten
classes. In Sections~\ref{flat}--\ref{QLefschetz}, we
have extracted the relevant parts of \cite{G}
and \cite{Kim}, necessary to state the quantum Lefschetz
hyperplane principle.

Throughout
the remainder of this paper, we will use the
following notation. Let $T_{0},\ldots,T_{m}$
and $T^0,\ldots,T^m$ be homogeneous bases for
$A_{\QQ}^{\ast}(X)$ such that $T_{0}=1$, $\{T_{1},\ldots,T_{r}\}$
is a basis for $A_{\QQ}^1(X)$, and $\lav T_{i},T^j\rav=\delta_{ij}$.
We shall also use $p_{i}=T_{i}$ to denote the divisor classes.
The coordinates of $T_{i}$ on $A_{\CC}^{\ast}(X)$ are
denoted by $t_{i}$. Also let $q_{i}$ be formal parameters 
for $i=1,...,r$.

\bigskip

\subsection{Gravitational descendents}\label{gravdesc}

For $n\ge 0$ consider the natural maps
\begin{equation}
        \begin{CD}
           \smap\sto @>(e,\eta)>>  X^n\times\nk @>{\eta}>> \smap_{0,n}  \\
         @. @V{e}VV \\
          {} @.X^n
        \end{CD}
\end{equation}
where $e=(e_{1},\ldots,e_{n}) $ is the evaluation map, and
$\eta$ forgets the map. By a slight abuse of notation, we
shall also use $e$ and $\eta$ to denote the projection maps.
Let $\gamma_{i}$ denote
arbitrary classes in $A^{\ast}(X)$.
In \cite{KM,B,BF,BM} a virtual fundamental class
$[\smap\sto]\in A_{\ed_{\beta}}(\smap\sto)$ has been constructed
such that the (tree level) {\em Gromov-Witten (GW) classes}
\begin{equation}
(\gamma_1,\ldots,\gamma_n)_{0,\beta}=\eta_{\ast}\left(
\prod_{i=1}^{n}e_{i}^{\ast}(\gamma_{i})\cdot[\smap\sto]\right)
\end{equation}
satisfy the axioms introduced by Kontsevich and Manin.
Here $[\smap\sto]$ has been identified with the
Poincar\'e dual of its image in $X^n\times \smap_{0,n}$.

More generally, we shall consider the classes of (tree level)
{\em gravitational descendents}
\[(\tau_{d_1}\gamma_1,\ldots,\tau_{d_n}\gamma_n)_{0,\beta}=
\eta_{\ast}\left(\prod_{i=1}^{n}c_1(\L_{i})^{d_{i}}
e_{i}^{\ast}(\gamma_{i})\cdot[\smap\sto]\right)\]
and the {\em invariants}
\[\lav\tau_{d_1}\gamma_1,\ldots,\tau_{d_n}\gamma_n\rav_{0,\beta}=
\int_{[\smap\sto]}\prod_{i=1}^nc_{1}(\L_{i})^{d_{i}}e_{i}^{\ast}
(\gamma_i)\,,\]
where $\L_i$ is the line bundle $s_{i}^{\ast}\omega_{\pi}$,
and $d_{i}$ are non-negative integers. We extend this definition 
to include negative $d_{i}$'s by defining 
$(\tau_{d_1}\gamma_1,\ldots,\tau_{d_n}\gamma_n)_{0,\beta}=0$ 
if some $d_{i}<0$. Here, $\omega_{\pi}$ is
the relative canonical line bundle of the universal map, hence the
fibers of $\L_{i}$ over a point $[C,f,s_{1},\ldots,s_{n}]$ are
$T_{s_{i}}C^{\vee}$. Occasionally, we abbreviate the
notation and write $(\otimes_{i=1}^{n}\tau_{d_{i}}\gamma_{i})_{0,\beta}$
instead of
$(\tau_{d_1}\gamma_1,\ldots,\tau_{d_n}\gamma_n)_{0,\beta}$.
Gravitational descendents satisfy the following analogues of
the KM-axioms:\footnotemark
\footnotetext{Note that we only work with even
dimensional cohomology, as opposed to \cite{KM} which is
formulated with odd cohomology as well.}
\\\\
{\em Fundamental class.}
If $\beta\ne 0$, then
\[
\lav\tau_{d_1}\gamma_1,\ldots,\tau_{d_n}\gamma_n,1\rav_{0,\beta}\,=
\sum_{i=1}^{n}\lav\tau_{d_1}\gamma_1,\ldots,\tau_{d_i-1}\gamma_i,\ldots,
\tau_{d_n}\gamma_n\rav_{0,\beta}.
\]
\\\\
 {\em Divisor.}
If $\beta\ne 0$ and $D\in A^{1}(X)$, then
\begin{multline*}
\pi_{\ast}(\tau_{d_1}\gamma_{1},\ldots,\tau_{d_n}
\gamma_{n},D)_{0,\beta}=\\
\biggl(\int_{\beta}
{D}\biggr)(\tau_{d_1}\gamma_{1},\ldots,\tau_{d_n}\gamma_{n})_{0,\beta}
+
\sum_{i=1}^n(\tau_{d_1}\gamma_{1},\ldots,\tau_{d_{i}-1}D\cdot\gamma_{i},\ldots,
\tau_{d_n}\gamma_{n})_{0,\beta}.
\end{multline*}
\\\\
{\em Splitting.}
Consider the commutative diagram, where the upper left and
lower right squares are fibred:
\begin{equation}
\label{diagsplit}
        \begin{CD}
          \nk @<\varphi_{AB}<<
          \M_{0,A\smallcup\{\bullet\}}\times
          \M_{0,B\smallcup\{\bullet\}}
           @=  \M_{0,A\smallcup\{\bullet\}}\times
           \M_{0,B\smallcup\{\bullet\}} \\
         @A\eta AA    @AA{\eta_{A}\times\eta_{B}}A
          @AA{\eta_{A}\times\eta_{B}}A         \\
          \nk(X,\beta) @<{G }<<
          \coprod_{\beta}
          @>\psi_{1}\times\psi_{2}>>{\coprod_{\beta}^{\prime}}   \\
         @. @VV{(e_{A},e_{B},e_{\bullet})}V
          @VV{(e_{A},e_{B},e_{\bullet},e_{\bullet})}V \\
         {} @.X^n\times X @>\Delta>> X^n\times X^2
                \end{CD}
\end{equation}
Here, $\Delta$ is the diagonal map on the last factor, $G$
is the gluing map,
$\psi_{1}$ and $\psi_{2}$ are projection maps,
\[\coprod_{\beta}=\coprod_{\beta_{1}+\beta_{2}=\beta}\M_{0,A\smallcup
\{\bullet\}}(X,\beta_{1})
\times_{X}\M_{0,B\smallcup\{\bullet\}}(X,\beta_{2})\,,\]
and
\[{\coprod_{\beta}}'=\coprod_{\beta_{1}+\beta_{2}=\beta}
\M_{0,A\smallcup\{\bullet\}}(X,\beta_{1})
\times\M_{0,B\smallcup\{\bullet\}}(X,\beta_{2})\,.\]
In \cite{BM} and \cite{B} it is shown that the virtual
fundamental class defines an orientation for $\smap$.
An important consequence of this is:
\[
\varphi_{AB}^![\smap_{0,n}(X,\beta)]
=\sum_{\beta_1+\beta_2=\beta}\Delta^!\bigl(
[\smap_{0,A\cup\{\cdot\}}(X,\beta_{1})]
\times[\smap_{0,B\cup\{\cdot\}}(X,\beta_{2})]\bigr).
\]
Multiply each side of this equation with
$\prod_{i=1}^{n}c_1(\L_{i})^{d_{i}}
e_{i}^{\ast}(\gamma_{i})$,
push down by $\eta_{A}\times\eta_{B}$, and use the projection formula 
and commutativity of the diagram 
to find
\[
\varphi_{AB}^{\ast}
(\tau_{d_1}\gamma_1,\ldots,\tau_{d_n}\gamma_n)_{0,\beta}=
\sum_{\beta_1+\beta_2=\beta}\sum_f(T_{f},\otimesover{a\in
A}\tau_{d_a}\gamma_a)_{0,\beta_1}
(T^f,\otimesover{b\in B}\tau_{d_b}\gamma_b)_{0,\beta_2}\,.
\]
\\\\
{\em Mapping to a point.}
If $\beta=0$ then $[\smap\sto]=\M_{0,n}\times X$ and
\begin{equation*}
\lav\tau_{d_1}\gamma_1,\ldots,\tau_{d_n}\gamma_n\rav_{0,n}=
\int_{X}{\prod_{i=1}^n\gamma_i}\cdot
\int_{\M_{0,n}}{\prod_{i=1}^n{c_{1}(\L_{i})^{d_{i}}} }\,,
\end{equation*}
where
\begin{equation*}
\prod_{i=1}^n{c_{1}(\L_{i})^{d_{i}}}=
\left\{
     \begin{array}{lcl}
     \frac{(n-3)!}{d_{1}!\cdots d_{n}!}&&{\rm if}\quad\sum_{i=1}^n d_{i}=n-3\\
     &&\\
     0       &&{\rm else}\,.
     \end{array}
\right.
\end{equation*}
See \cite{W} and \cite{K} for the last claim.

\medskip

\noindent{\em Remark~1.}
 If $\smap\sto$ is of expected dimension, then
$[\smap\sto]$ is just the fundamental class of $\smap\sto$.

\medskip

\noindent{\em Remark~2.} With $d_i=0$ for all $i$, we recover the
axioms in \cite{KM}.

\medskip

\noindent{\em Remark~3.} Note that if all $\gamma_{i}$ are
homogeneous, then
$\lav\tau_{d_1}\gamma_1,\ldots,\tau_{d_n}\gamma_n\rav_{0,\beta}=0$
unless $\sum_{i=1}^n\codim(\gamma_{i})+d_{i}$ is equal to
the expected dimension of $\smap\sto$. The same claim is
also true unless $\beta=0$ or $\beta$ is effective, since $\smap\sto$ is
then empty.

\medskip

\noindent{\em Remark~4.} The fundamental class and divisor
properties follow from the relations
\begin{equation}
c_{1}(\L_{i})
-\pi^{\ast}c_{1}(\L_{i})=s_{i}(\smap\sto)
\quad{\rm for}\quad i=1, \ldots, n\,.
\end{equation}
Rigorous proofs of these claims can be found in \cite{P}
when $X$ is convex. These proofs hold in a more general
setting, since the virtual fundamental
classes satisfy
$\pi^{\ast}[\smap\sto]=[\M_{0,n+1}(X,\beta)]$.

\begin{ex}\label{exCY3fold}{\em
If $\ed_{\beta}=0$ and $n=0$, then $\lav\,\,\rav_{0,\beta}$
is simply the degree of $[\moo(X,\beta)]$. An interesting example is
the class of Calabi-Yau threefolds. A variety is Calabi-Yau
if its canonical bundle is trivial. Hence, if $X$ is
Calabi-Yau of dimension $3$, the expected dimension of
$\moo(X,\beta)$ is $0$ for all effective $\beta$. Let
$N_{\beta}^X$ denote the degree of $[\moo(X,\beta)]$. The
numbers $N_{\beta}^X$ are related to the virtual numbers\footnotemark
\footnotetext{The virtual number $n_{\beta}^{X}$ is defined in
such a way that if the rational curves of class $\beta$ in
$X$ are rigid, then $n_{\beta}$ is the actual number of such curves.
See \cite{Katz} for a discussion.}
$n_{\beta}^{X}$ of rational curves of class $\beta$ in $X$.
The Aspinwall-Morrison \cite{AM} formula describes this
relationship. For simplicity, we shall assume that $\Pic
(X)\simeq\ZZ$, and identify $\beta$ with the degree
$d=\lav p,\beta\rav$, where $p$ is the ample generator of
$\Pic (X)$. Then the Aspinwall-Morrison formula reads
\begin{equation}
N_{d}^{X}=\sum_{k|d}k^{-3}n_{d/k}^{X}\,.
\end{equation}
In \cite{K}, this formula was explained as follows: If
$d>1$, then $\moo(X,d)$ is rarely of expected dimension
since the locus of multiple coverings of lower degree
curves will be of positive dimension. Hence we need to
compute the contributions to $[\moo(X,d)]$ from this locus.
Suppose $C\simeq\PP^1$ is a smooth
rational curve of degree $d/k$ in $X$, and $C$ is rigid,
i.e., $N_{C/X}\simeq\oh_{\PP^1}(-1)\oplus\oh_{\PP^1}
(-1)$. Then the locus of $k:1$ coverings of $C$ in
$\moo(X,d)$ is a $2k-2$ dimensional component which is
isomorphic to $\moo(\PP^1,k)$. By cohomology and base change
theorems in \cite{H}, the obstruction sheaf
$R^1\pi_{\ast}e_{1}^{\ast}N_{C/X}$ is a vector bundle of
rank $2k-2$ on $\moo(\PP^1,k)$. The contribution to
$N_{d}^{X}$ from the component $\moo(\PP^1,k)$ is the
integral
\[\int_{\moo(\PP^1,k)}c_{\top}(R^1\pi_{\ast}e_{1}^{\ast}N_{C/X})\,.\]
An evaluation of this integral using Bott's formula yields
$1/k^3$. See also \cite{GrP} for more general multiple
cover computations.
}\end{ex}

\bigskip

\subsection{The potential}\label{potential}
For non-negative integers $i$, and $j=0,\ldots,m$,
let $t_{i}^j$ be formal variables
such that $t_{0}^j=t_j$, and let 
$T=\sum_{i,j}t_{i}^j\tau_{i}T_{j}$. The 
{\em gravitational potential} is defined as 
\[
\Phi=
\sum_{n\ge 0}\frac{1}{n!}\sum_{\beta}{q^{\beta}\lav T^{\otimes
n}\rav_{0,\beta}}\,.\]
Using the linearity of the gravitational descendents we obtain a
 formal power series in $\QQ[[q_{i},q_{i}^{-1},t_{i}^{j}]]$,
\[
\Phi=
\sum_{n\ge 0}\sum
\sum_{\beta}
\lav\otimes_{i,j}(\tau_{i}T_{j})^{\otimes n_{i}^j}\rav_{0,\beta}
q^{\beta}\prod_{i,j}
\frac{(t_{i}^j)^{n_{i}^j}}
{n_{i}^j!}\,,
\]
where the second summation is over all partitions
$\sum_{i,j}n_{i}^j=n$. 
We use Witten's notation to denote the partial derivatives of $\Phi$:
\[\laav\tau_{d_1}T_{j_{1}},\ldots,\tau_{d_k}T_{j_{k}}\raav_{0}=
\partial_{d_1}^{j_{1}}\ldots\partial_{d_k}^{j_{k}}\Phi\,,\]
where $\partial_{k}^l=\frac{\partial}{\partial t_{k}^l}\,.$ This
is extended linearly to define
 $\laav\tau_{d_1}\gamma_{j_{1}},\ldots,\tau_{d_k}\gamma_{j_{k}}\raav_{0}$ 
for arbitrary classes $\gamma_{i}=\sum a_{ij}T_{j}$ 
in $A_{\QQ}^{\ast}(X)$. 
It is easy to verify that 
\begin{equation}\label{to-tre}
\laav\tau_{d_1}\gamma_1,\ldots,\tau_{d_k}\gamma_k\raav_{0}=
\sum_{n\ge 0}\frac{1}{n!}\sum_{\beta}{q^{\beta}
\lav\tau_{d_1}\gamma_1,\ldots,\tau_{d_k}\gamma_k,
T^{\otimes n}\rav_{0,\beta}}\,.
\end{equation}
Also, let
\[\lav\tau_{d_1}\gamma_1,\ldots,\tau_{d_k}\gamma_k\rav_{0}=
\laav\tau_{d_1}\gamma_1,\ldots,\tau_{d_k}\gamma_k\raav_{0}
\mid_{\{t_{i}^j=0|i>0\}}\,.\]

\medskip

\noindent{\em Remark.} Several variations of the potential 
appear in the literature on GW-invariants. If we restrict 
$\Phi$ to $\{q_{i}=1\}\cup\{t_{i}^j=0|i>0\}$ we recover the potential
in \cite{FP}. On the other hand, restricting $\Phi$ to 
$\{t_{i}=0|1\le i\le r\}\cup\{t_{i}^j=0|i>0\}$ yields the 
potential defined in \cite{GP}. Using the fundamental 
and divisor properties one can deduce that these are 
equivalent if we identify $q_{i}=e^{t_{i}}$.

The entire structure of (gravitational) quantum cohomology
is entailed in the fact that the potential is a solution to
the following equations:

\begin{thm}{\em (WDVV equations)}\label{th2-1}
\begin{equation*}
\sum_{f}{
\laav\tau_{a}\gamma_{a},\tau_{b}\gamma_{b},T_{f}\raav_{0}\,
\laav T^{f},\tau_{c}\gamma_{c},\tau_{d}\gamma_{d}\raav_{0}}=
\sum_{f}{
\laav\tau_{a}\gamma_{a},\tau_{c}\gamma_{c},T_{f}\raav_{0}\,
\laav
T^{f},\tau_{b}\gamma_{b},\tau_{d}\gamma_{d}\raav_{0}}\,.
\end{equation*}
\end{thm}

\begin{proof}
Recall the fundamental relations \eqref{1-1} on
$\smap_{0,n+4}$,
\begin{equation}
\label{e:eq}
\sum D(A\cup\{ab\}|B\cup\{cd\})\simeq
\sum D(A\cup\{ac\}|B\cup\{bd\})\,,
\end{equation}
where the summations run over all partitions $A\cup
B=\{1,\ldots,n\}$.
Multiply each side of \eqref{e:eq} by
$(T^{\otimes n},\tau_{a}\gamma_{a},\tau_{b}\gamma_{b},
\tau_{c}\gamma_{c},\tau_{d}\gamma_{d})_{0,\beta}$ and use
splitting to establish
\begin{multline*}
\sum\sum_{\beta_{1}+\beta_{2}=\beta}\sum_{f}
\lav
T^{\otimes|A|},\tau_{a}\gamma_{a},\tau_{b}\gamma_{b},T^{f}\rav_{0,\beta_{1}}
\lav
T^{\otimes|B|},\tau_{c}\gamma_{c},\tau_{d}\gamma_{d},T^{f}\rav_{0,\beta_{2}}
=\\
\sum\sum_{\beta_{1}+\beta_{2}=\beta}\sum_{f}
\lav
T^{\otimes|A|},\tau_{a}\gamma_{a},\tau_{c}\gamma_{c},T_{f}\rav_{0,\beta_{1}}
\lav
T^{\otimes|B|},\tau_{b}\gamma_{b},\tau_{d}\gamma_{d},T_{f}\rav_{0,\beta_{2}}\,,
\end{multline*}
where the first summations are as in \eqref{e:eq}. Since there
are  $\binom{n_{\,}}{n_{1}}$
partitions of $\{1,\ldots,n\}$ such that
$|A|=n_{1}$, we have
\begin{eqnarray*}
\sum_{n_{1}+n_{2}=n}\sum_{\beta_{1}+\beta_{2}=\beta}\!\!\!\sum_{f}
\frac{n!}{n_{1}!n_{2}!}\!\!
\lav\! T^{\otimes
n_{1}},\tau_{a}\gamma_{a},\tau_{b}\gamma_{b},T^{f}\!\rav_{0,\beta_{1}}
\lav\! T^{\otimes
n_{2}},\tau_{c}\gamma_{c},\tau_{d}\gamma_{d},T^{f}\!\rav_{0,\beta_{2}}
=\\
\sum_{n_{1}+n_{2}=n}\sum_{\beta_{1}+\beta_{2}=\beta}\!\!\!\sum_{f}
\frac{n!}{n_{1}!n_{2}!}\!\! \lav\! T^{\otimes
n_{1}},\tau_{a}\gamma_{a},\tau_{c}\gamma_{c},T_{f}\!\rav_{0,\beta_{1}}
\lav \!T^{\otimes
n_{2}},\tau_{b}\gamma_{b},\tau_{d}\gamma_{d},T_{f}\!\rav_{0,\beta_{2}}\,.
\end{eqnarray*}
Dividing by $n!$, multiplying by $q^{\beta}$, and summing over all $n$ and
$\beta$ yields the desired result.
\end{proof}

\begin{thm}{\em (Witten)}
\label{thm}
\begin{equation}
\laav\tau_{a+1}\gamma_{a},\tau_{b}\gamma_{b},\tau_{c}\gamma_{c}\raav_{0}\,=
\sum_{f}
\laav\tau_{a}\gamma_{a},T_{f}\raav_{0}
\laav T^{f},\tau_{b}\gamma_{b},\tau_{c}\gamma_{c}\raav_{0}\,.
\end{equation}
\end{thm}

\begin{proof}
This result is a consequence of the following well-known
relation:
\[c_{1}(\L_{1})=\sum_{K\cup L=\{ 4,\ldots,n\}}
\eta^{\ast}D(K\cup\{1\}|L\cup\{2,3\} )\,.
\]
 Applying it to
 $\laav\tau_{a+1}\gamma_{a},\tau_{b}\gamma_{b},
 \tau_{c}\gamma_{c}\raav_{0}$
we find
\begin{equation*}
\sum_{n\ge 0}\sum_{\beta}\frac{1}{n!}
\sum_{K\cup L=\{ 4,\ldots,n+3\}}q^{\beta}\int_{\smap_{0,n+3}}
\biggl(D(K\cup\{1\}|L\cup\{2,3\} )
(\tau_{a}\gamma_{a},\tau_{b}\gamma_{b},\tau_{c}\gamma_{c},T^{\otimes n}
)_{0,\beta}\biggr)\,.
\end{equation*}
The result follows by splitting and a combinatorial argument
as in the proof of Theorem~\ref{th2-1}.
\end{proof}

We now give an example due to Kontsevich of the remarkable
enumerative implications of the WDVV-equations.

\medskip

\begin{ex}
\label{exPP2} {\em (Quantum cohomology of ${\PP^2}$)
Let $T_{i}=p^i$, $i=0,1,2$, where $p$ is the ample
generator of $\Pic(\PP^2)$. For  $d\ge 1$, let
\[K_{d}=\lav\underbrace{ T_{2},\ldots,T_{2}}_{3d-1}\rav_{0,d}\,.\]
Then $K_{d}$ is the number of rational curves of degree $d$
through $3d-1$ points, and the GW-potential of $\PP^2$ is
$\Phi|_{\{t_{i}^j=0|i>0\}}=\frac{1}{2}(t_{0}t_{1}^2+
t_{0}^2t_{1})+\Gamma$, where
\[\Gamma=\sum_{d=1}^{\infty}K_{d}q^d
\frac{t_{2}^{3d-1}}{(3d-1)!}\,.\]
A direct application of Theorem~\ref{th2-1} yields
\[
\partial_{2}\partial_{2}\partial_{2}\Gamma=
(\partial_{1}\partial_{1}\partial_{2}\Gamma)^2-
(\partial_{1}\partial_{1}\partial_{1}\Gamma)\cdot
(\partial_{1}\partial_{2}\partial_{2}\Gamma)\,,\]
and equating the coefficients on both sides, we find
Kontsevich's recursion formula:
\[K_{d}=\sum_{d_{1}+d_{2}=d}K_{d_{1}}K_{d_{2}}d_{1}^2d_{2}
\bigg[d_{2}\binom{3d-4}{3d_{1}-2}-
d_{1}\binom{3d-4}{3d_{1}-1}\bigg]\,.\]
}\end{ex}

\bigskip

\subsection{Quantum cohomology rings}\label{QCR}

Let us consider the free $\QQ [[q,q^{-1},t]]$-module
\(A^{\ast}(X)\otimes\QQ[[q,q^{-1},t]]\)
generated by $T_{0},\ldots,T_{m}$.
The {\em quantum product} is defined as
\[
T_{i}\ast T_{j}=\sum_{f}\lav T_{i},T_{j},T_{f}\rav_{0}T^{f}\,.
\]
Extending this product $\QQ[[q,q^{-1},t]]$-linearly puts a
$\QQ[[q,q^{-1},t]]$-algebra structure on
$A^{\ast}(X)\otimes\QQ[[q,q^{-1},t]]$. This ring is called the (big) quantum
cohomology ring $QH_{b}^{\ast}(X)$.

\begin{prop}
$QH_{b}^{\ast}(X)$ is an associative, commutative ring with
identity $T_{0}=1$.
\end{prop}

\begin{proof}
Commutativity is obvious. That $T_{i}\ast T_{0}=T_{i}$ for
all $i$ follows from the fundamental class and mapping to a
point property.
Finally, associativity is equivalent to the WDVV equations
for the correlators $\lav\,\,\,\rav_{0}$.
\end{proof}

We shall mainly work with the (small) quantum cohomology
ring $QH^{\ast}(X)$ which is the specialization of
$QH_{b}^{\ast}(X)$ obtained by restricting 
$t=(t_{1},\ldots,t_{m})=(0,\ldots,0)$. Using
\eqref{to-tre}, one can see that the product on
$QH^{\ast}(X)$ is
\[
T_{i}\ast
T_{j}=\sum_{\beta}q^{\beta}\sum_{f}\lav
T_{i},T_{j},T_{f}\rav_{0,\beta}T^{f}\,. \]
We make the ring $QH^{\ast}(X)$ graded by putting
$\deg q_{i}=\lav T^i,c_{1}(TX)\rav$, and letting $T_{i}$ keep
its grading from $A^{\ast}(X)$.

\medskip

\noindent{\em Remark.}
 If the gravitational potential $\Phi$ is a formal
power series we can also define
\[T_{i}\ast T_{j}=\sum_{f}\laav T_{i},T_{j},T_{f}\raav_{0}T^{f}.\]
This puts a $\QQ[[t_{i}^j]]$-algebra
structure on
 $A^{\ast}(X)\otimes\QQ[[t_{i}^j]]$.
This ring is still commutative and associative, but due to
the more complicated nature of the fundamental class
property, $T_{0}$ is no longer the identity.

\medskip

\begin{ex}
\label{exPPn} {\em
The quantum cohomology ring of $\PP^n$ is
$QH^{\ast}(\PP^n)=\QQ[p,q]/(p^{n+1}-q)$, where $p$ is the
ample generator of $\Pic(\PP^n)\simeq\ZZ$.
}\end{ex}

\bigskip

\subsection{Dubrovin's formalism}\label{Dubrovin}

In this section we describe yet another way one may view
the WDVV-equations, as developed in \cite{D}.
First we review some general theory about integrability
conditions and flat connections.

Let $M$ be a vector bundle on ${\CC}^m$ with coordinates
$x_{i}$ on ${\CC}^m$ and $m_{j}$ on the fibers, and let
$\hat{m}_{j}=(x_{1},\ldots,x_{m},0,\cdots,0,m_{j},0,\ldots,0)$
be sections of $M$.

Suppose we have a connection $\nabla\:M\to
M\otimes\Omega_{\CC^m}$ with covariant derivatives
\[\nabla_{i}(g_{j}\hat{m}_{j})=
\frac{\partial g_{j}}{\partial x_{i}}\hat{m}_{j}+
\sum_{k}\Gamma_{ij}^k g_{j}\hat{m}_{k}\,,\]
where $\Gamma_{ij}^k$ are holomorphic functions on $\CC^m$.
The connection defines a distribution \[\bbD\:T{\CC}^m\to
TM\] on $TM$ by letting
\[\frac{\partial}{\partial x_{p}}\mapsto\frac{\partial}{\partial x_{p}}-
\sum_{j,k}\Gamma_{pj}^k m_{j}\frac{\partial}{\partial
m_{k}}\,.\]
By Frobenius' theorem, a distribution is integrable if it is
closed under the Lie bracket operation $[\;,\;]$. It is
straightforward to verify the
following facts:
\\\\
$\bbD$ is closed under  $[\;,\;]$ if and only if $\nabla$
is flat.
\\
A section $F\:{\CC}^m\to M$ is an integral
of $\bbD$ if and only if it is flat, i.e., $\nabla(F)=0$.

\bigskip

We now apply this theory to the affine space $V=A_{\CC}^{\ast}(X)$ and
the trivial vector bundle $M=V\times V$. Consider the ``connection''
defined by the covariant derivatives
\begin{equation}
\label{e:eqq}
\nabla_{i}(F)=\hbar\frac{\partial F}{\partial
t_{i}}-T_{i}\ast F\,.
\end{equation}
Here, $\hbar$ is just a formal parameter. If $\hbar=1$, then
\eqref{e:eqq} is a formal connection. It is straightforward to verify:
\begin{prop}
$\nabla$ is flat if and only if the quantum product $\ast$ is associative.
\end{prop}

\bigskip

\subsection{Flat sections}\label{flat}

 From the discussion in the previous section, we can conclude that
there should exist formal solutions to the system of first
order differential equations
\begin{equation}
\label{e:nabla}
\hbar\frac{\partial}{\partial t_{i}}=T_{i}\ast
\end{equation}
where $i=0,\ldots,m$. Via gravitational descendents we can actually
obtain closed formulas for these solutions.

For $a=0,\ldots,m$, let
\begin{equation}
\label{e:r1}
\vec{S}_{a}=T_{a}+\sum_{n}\sum_{b}{\hbar}^{-n-1}\lav\tau_{n}T_{a},T_{b}\rav_
{0}T^b.
\end{equation}

\noindent{\bf Claim.} $\vec{S}_{a}$ is a solution to the
system of first order differential equations \eqref{e:nabla}.

\begin{proof}
This is a consequence of the equations in Theorem~\ref{thm},
specialized to the partial derivatives $\lav\,\,\,\rav_{0}$. On one
side of
\eqref{e:nabla} we find
\[{\hbar}\frac{\partial\vec{S}_{a}}{\partial t_{i}}=
\sum_{n}\sum_{b}{\hbar}^{-n}\lav\tau_{n}T_{a},T_{b},T_{i}\rav_{0}T^b\,,\]
and on the other side
\[
T_{i}\ast \vec{S_{a}}=\sum_{b}\lav T_{i},T_{a},T_{b}\rav_{0}T^b+
\sum_{n}\sum_{b,c}{\hbar}^{-n-1}\lav\tau_{n}T_{a},T_{b}\rav_{0}\lav
T^b,T_{i},T_{c}\rav_{0}T^c\,,\]
which by Theorem~\ref{thm} equals
\[\sum_{b}\lav T_{i},T_{a},T_{b}\rav_{0}T^b+
\sum_{n}\sum_{c}{\hbar}^{-n-1}\lav\tau_{n}T_{a},T_{i},T_{c}\rav_{0} T^c \,.\]
\end{proof}

We shall restrict to only consider the system \eqref{e:nabla} with
$i=0,\ldots,r$ and the small quantum product.
A formal calculation
using the fundamental class, divisor, and mapping to a
point properties yields:
\\
The system of differential equations
\begin{equation}
{\hbar}\frac{\partial}{\partial t_{0}}=T_{0}\ast\,\,\,,\qquad\qquad
{\hbar}\frac{\partial}{\partial t_{i}}=p_{i}\ast\qquad
{\rm for}\quad i=1,\ldots,r\,,
\end{equation}
has solutions
\begin{equation}
\label{e:r2}
\vec{s}_{a}=e^{t_{0}/\hbar}\biggl(e^{pt/{\hbar}}T_{a}+\sum_{\beta\ne 0}q^{\beta}
\sum_{b}T^b\int_{[\smap_{0,2}(X,\beta)]}
\frac{e_{1}^{\ast}(e^{pt/{\hbar}}T_{a})}{{\hbar}-c}
e_{2}^{\ast}(T_{b})\biggr)\,,
\end{equation}
where $a=0,\ldots,m$, $pt=\sum p_{i}t_{i}$, $c=c_{1}(\L_{1})$, and
$q_{i}=e^{t_{i}}$.

\bigskip

\subsection{Quantum $\D$-modules}\label{QDmodule}

The quantum $\D$-module $Q\D(X)$ of $X$ is the $\D$-module
generated by all $\lav\vec{s},1\rav$ for flat sections
$\vec{s}$. In other words, if $\D$ is the ring of linear
differential operators $\CC[{\hbar}\frac{\partial}{\partial
t_{i}},q_{i},{\hbar}\mid i=0,\ldots,r]$, with
\[\left[{\hbar}\frac{\partial}{\partial t_{i}},q_{i}\right]=
{\hbar}\frac{\partial}{\partial t_{i}}\circ q_{i}
-q_{i}\circ {\hbar}\frac{\partial}{\partial t_{i}}={\hbar}q_{i}\,,\]
then the quantum $\D$-module is the set of $P\in \D$ such
that $P(\lav\vec{s},1\rav)=0$ for all flat sections $\vec{s}$.
Using the formula \eqref{e:r2}, we can write up a set of
generators for the quantum $\D$-module in vector form:
\begin{equation}
s^X=\sum_{a=0}^m\lav\vec{s}_{a},1\rav T^a=
e^{{(t_{0}+pt)}/{\hbar}}\left(
1+\sum_{\beta\ne 0}q^{\beta}e_{1\ast}\left(\frac{1}{{\hbar}-c}\right)
\right)
\end{equation}
where $e_{1}\:\M_{0,2}(X,\beta)\to X$.

\medskip

\noindent{\em Remark.}
The $\D$-module can be useful in finding relations in
$QH^{\ast}(X)$. If
$P({\hbar}\frac{\partial}{\partial t_{i}},q_{i},{\hbar})$
is in $Q\D(X)$, then $P(p_{i},q_{i},0)=0$ in $QH^{\ast}(X)$. A similar
result holds in $QH_{b}^{\ast}(X)$.

\medskip

We close this section with a description of the quantum
$\D$-module in two cases.

\medskip

\begin{ex}\label{exPP5} $\PP^{4}$.
{\em Let $p$ be the ample generator of $\Pic(\PP^4)$. Then
the multiplication matrix $p\ast-$ with respect to the
basis $1,p,\ldots,p^4$ is
\begin{equation*}
 p\ast 1=p\,,\qquad p\ast p=p^2\,,\qquad
 p\ast p^2=p^3\,,\qquad p\ast p^3=p^4\,,\qquad
 p\ast p^4=q\,.
\end{equation*}
Hence a section $\vec{s}=s_{0}+s_{1}p+\cdots+s_{4}p^4$ is
flat if and only if
\begin{equation*}
\hbar\frac{d}{dt}\begin{bmatrix}s_{0}\\ s_{1}\\ s_{2}\\ s_{3}\\
s_{4}\end{bmatrix}
=\begin{bmatrix}
0&0&0&0&q\\
1&0&0&0&0\\
0&1&0&0&0\\
0&0&1&0&0\\
0&0&0&1&0
\end{bmatrix}
\begin{bmatrix}s_{0}\\ s_{1}\\ s_{2}\\ s_{3}\\s_{4}\end{bmatrix}=
\begin{bmatrix}qs_{4}\\ s_{0}\\ s_{1}\\ s_{2}\\ s_{3}\end{bmatrix}\,.
\end{equation*}
Solving for $s_{4}$, we find
\begin{equation*}
\hbar\frac{d}{dt}s_{4}=s_{3}\,,\quad
(\hbar\frac{d}{dt})^{2}s_{4}=s_{2}\,,\quad
(\hbar\frac{d}{dt})^{3}s_{4}=s_{1}\,,\quad
(\hbar\frac{d}{dt})^{4}s_{4}=s_{0}\,,\quad
(\hbar\frac{d}{dt})^{5}s_{4}=qs_{4}\,,
\end{equation*}
hence $(\hbar\frac{d}{dt})^5-q$ generates the quantum
$\D$-module of $\PP^4$. Note that $p^5-q=0$ in
$QH^{\ast}(\PP^4)$.
}\end{ex}

\medskip

\begin{ex}\label{exCY}
{\em
 Suppose $Y$ is a Calabi-Yau threefold with
 $\Pic(Y)\simeq\ZZ$. Let $p$ be the ample generator of
 $\Pic(Y)$. Let $N_{d}^Y$ and $n_{d}^Y$ be as in
 Example~\ref{exCY3fold} of Section~\ref{gravdesc}. From the
 definition of the quantum product in $QH^{\ast}(Y)$, we find
 \begin{equation}
 p\ast 1=p\,,\qquad p\ast p=K(q)p^2\,,\qquad
 p\ast p^2=p^3\,,\qquad p\ast p^3=0\,,
 \end{equation}
 where
 \[K(q)=1+\frac{1}{\deg(Y)}\sum_{d>0}\lav
 p,p,p\rav_{0,d}q^d=1+\frac{1}{\deg(Y)}\sum_{d>0}
 d^3N_{d}^Yq^d\,.\]
 A section $\vec{s}=s_{0}+s_{1}p+s_{2}p^2+s_{3}p^3$ is flat
 if and only if
\begin{equation*}
\hbar\frac{d}{dt}\begin{bmatrix}s_{0}\\ s_{1}\\ s_{2}\\ s_{3}\end{bmatrix}
=\begin{bmatrix}
0&0&0&0\\
1&0&0&0\\
0&K(q)&0&0\\
0&0&1&0
\end{bmatrix}
\begin{bmatrix}s_{0}\\ s_{1}\\ s_{2}\\ s_{3}\end{bmatrix}
\,.
\end{equation*}
Solving for $s_{3}$, we find that
$(\hbar\frac{d}{dt})^2\frac{1}{K(q)}(\hbar\frac{d}{dt})^2$ generates the
quantum $\D$-module of $Y$. Note that by the
Aspinwall-Morrison formula
\begin{equation}
K(q)=1+\frac{1}{\deg(Y)}\sum_{d>0}d^3n_{d}^Y
\frac{q^d}{1-q^d}\,.
\end{equation}
 }\end{ex}

\bigskip

\subsection{Generalized quantum structures}\label{gqs}

Let $\B$ be a vector bundle of rank $b$ on $X$ that is generated by its
global sections and suppose $Y\sub X$ is the zero scheme of a
general global section of $\B$. The idea is to find a way to relate
the quantum theory of $Y$ to the quantum theory of $X$.

First, we look at the virtual fundamental class
$[\smap_{0,n}(Y,\beta)]$. Note that by the normal bundle
sequence and the adjunction formula
$c_{1}(TY)=c_{1}(TX)-c_{1}(\B)$, hence the expected dimension
of $\smap_{0,n}(Y,\beta)$ is
\[\dim X-b+\lav\beta,c_{1}(TX)-c_{1}(\B)\rav+n-3\,.\]

The inclusion map $i\:Y\to X$ induces a map
\[j\:\nk(Y,\beta)\to\nk(X,i_{\ast}\beta)\]
 such that
\begin{equation}
        \begin{CD}
    \coprodover{i_{\ast}\beta=d}\nk(Y,\beta) @>{j}>>\nk(X,d)\\
        @VV{e} V             @VV{e}V\\
        Y^n @>i>>    X^n
        \end{CD}
\end{equation}
commutes.
Let
\begin{equation}
        \begin{CD}
                 {} @. \B\\
                 @. @VV{}V \\
                \smap_{0,n+1}(X,d) @>e_{n+1}>>X \\
                @VV\pi V \\
                   \nk(X,d)
        \end{CD}
\end{equation}
be the universal map. Tensoring
${\oh}_{\nk(X,d)\times X}\to {\oh}_{\smap_{0,n+1}(X,d)}$
by $e_{n+1}^{\ast}\B$ and pushing down yields a map
\begin{equation}
H^{0}(X,\B)_{\nk(X,d)}\xrightarrow{\nu}
\pi_{\ast}e_{n+1}^{\ast}\B\,.
\end{equation}
For any point $[C,f,s_{i}]$ in $\nk(X,d)$, we have
$H^1(C,f^{\ast}\B)=0$, hence base change theorems in \cite{H}
imply that $\B_{d}=\pi_{\ast}e_{1}^{\ast}\B$ is a vector bundle
with fibers $H^0(C,f^{\ast}\B)$. Note that $\B_{d}$ is independent
of the marked points, in the sense that
$\pi_{i}^{\ast}\B_{d}=\B_{d}$, where $\pi_{i}$ forgets the
$i$-th point. The map $\nu$ sends a section
$s\in H^0(X,\B)$ to the restriction
$s\circ f\in H^0(C,f^{\ast}\B)$, so the stable map
$f\:C\to X$ factors through $i$ if and only if
$[C,f,s_{i}]$ lies in the zero locus of $\nu(s)$. In other
words, if $o$ is the zero section, then the diagram
\begin{equation}
        \begin{CD}
      \coprodover{i_{\ast}\beta=d}\nk(Y,\beta) @>{}>>\nk(X,d)\\
        @V{} VV             @VV{\nu(s)}V\\
        \nk(X,d) @>o>>    \B_{d}
        \end{CD}
\end{equation}
is fibred.
The construction of the classes $[\smap_{0,n}(Y,\beta)]$ is such
that
\[[\nk(Y,\beta)]=\left(o![\nk(X,d)]\right)_{\nk(Y,\beta)}\,,\]
where $(\;)_{\nk(Y,\beta)}$ indicates the proper component
of $\coprod_{i_{\ast}\beta=d}\nk(Y,\beta)$. From the excess
intersection formula \cite{F},
\begin{equation}
\label{e:e1}
j_{\ast}\sum_{i\ast\beta=d}[\nk(Y,\beta)]=B_{d}\cdot[\nk(X,d)]\,,
\end{equation}
with $B_{d}=c_{{\top}}(\B_{d})$. Note that if
$C\simeq\PP^1$, then by Riemann-Roch
\[\dim H^0(\PP^1,f^{\ast}\B)=b+\int_{d}c_{1}(\B)\,,\]
hence the rank of $\B_{d}$ is $b+\int_{d}c_{1}(\B)$, and
\eqref{e:e1} is a class of expected dimension.

We shall use this construction to define a new associative
quantum product on $A^{\ast}(X)\otimes\QQ[[q,q^{-1}]]$. The
product will account for GW-invariants on $Y$ involving
classes from $i^{\ast}A(X)$.

Define new correlators on $A^{\ast}(X)$
\[(\tau_{d_1}\gamma_1,\ldots,\tau_{d_n}\gamma_n)_{0,d}^{\B}=
\eta_{\ast}\!\left(\prod_{i=1}^n c_{1}(\L_{i})^{d_{i}}
e_{i}^{\ast}(\gamma_{i})\cdot B_{d}\cdot[\nk(X,d)]\right)\,.\]
As before let
$\lav\,\,\,\rav_{0,d}^{\B}=\int_{[\smap_{0,n}(X,d)]}(\,\,)_{0,d}^{\B}\,.$
Note that, by \eqref{e:e1},
\begin{equation}
(\tau_{d_1}\gamma_1,\ldots,\tau_{d_n}\gamma_n)_{0,d}^{\B}=
\sum_{i_{\ast}\beta=d}\bigl(\tau_{d_1}i^{\ast}(\gamma_1),\ldots,
\tau_{d_n}i^{\ast}(\gamma_n)\bigr)_{0,\beta}.
\end{equation}
Hence if for some $j$, $\, i^{\ast}(\gamma_j)=0$, then
$(\tau_{d_1}\gamma_1,\ldots,\tau_{d_n}\gamma_n)_{0,d}^{\B}=0$.

It is easy to see that the new correlators satisfy the fundamental class and
divisor properties of \ref{gravdesc}. We shall also see
that there is a splitting principle.
Let $\rho_{j}\colon\B_{d}\to e_{j}^{\ast}\B$ be the map
which on the fibers $H^{0}(C,f^{\ast}\B)$ of $\B_{d}$
evaluates the sections at the $j$th marked point. Associated
to the gluing map
\[\varphi\colon\M_{0,A\smallcup\{\bullet\}}(X,\beta_{1})
\times_{X}\M_{0,B\smallcup\{\bullet\}}(X,\beta_{2})\to
\M_{0,A\smallcup B}(X,\beta_{1}+\beta_{2})\,,\]
there is an exact sequence
\[0\to\varphi^{\ast}\B_{\beta_{1}+\beta_{2}}\to
\psi_{1}^{\ast}\B_{\beta_{1}}\oplus\psi_{2}^{\ast}\B_{\beta_{2}}
\stackrel{\rho_{\bullet}^1-\rho_{\bullet}^2}
{\rarr}e_{0}^{\bullet}\B\to
0\,,\]
where $\psi_{i}$ are projection maps. Hence, if
$B=c_{{\top}}(\B)$ and
 $B_{j,d}'=c_{{\top}}({\kker}\rho_{j})$, we have
\begin{equation}
c_{{\top}}(\varphi^{\ast}B_{\beta_{1}+\beta_{2}})\cdot
e_{0}^{\ast}B=B_{\beta_{1}}\cdot B_{\beta_{2}}\qquad{\rm
and }\qquad
B_{d}=B_{j,d}^{\prime}\cdot e_{j}^{\ast}B\,.
\end{equation}
Let $A\cup B$ be a partition of $\{1,\cdots,n\}$. By the
commutativity of the diagram \eqref{diagsplit},
\begin{equation}
\label{e:e5}
\begin{split}
\varphi_{AB}^{\ast}\!
(\tau_{d_1}\gamma_1,\ldots,\tau_{d_n}\gamma_n)_{0,d}^{\B}
&=(\eta_{A}\!\times\!\eta_{B})_{\ast}\varphi_{AB}^!
\left(\prod_{i=1}^n c_{1}(\L_{i})^{d_{i}}
e_{i}^{\ast}(\gamma_{i})\cdot[\nk(X,d)]\cdot B_{d}\right)\\
&=(\eta_{A}\!\times\!\eta_{B})_{\ast}(\psi_{1}\!
\times\!\psi_{2})_{\ast}
(G^{\ast}(B_{d}))
(e_{\bullet}\!\times\!
e_{\bullet})^{\ast}(\Delta)\!
\sum
\end{split}
\end{equation}
where
\[\sum=\!\!\sum_{\beta_{1}\!+\beta_{2}\!=d}\biggl(
\prod_{a\in A}\! c_{1}(\L_{a})^{d_{a}}
e_{a}^{\ast}(\gamma_{a})\!\cdot\![\M_{0,A
\smallcup\{\bullet\}}(X,\beta_{1})]
\times\!\!
\prod_{b\in B}\! c_{1}(\L_{b})^{d_{b}}
e_{b}^{\ast}(\gamma_{b})\!\cdot\![\M_{0,B\smallcup\{\bullet\}}
(X,\beta_{2})]\biggr)\,.\]
Since
\begin{equation}
\label{e:e4}
\begin{split}
(\psi_{1}\times\psi_{2})_{\ast}(G^{\ast}(B_{d}))
(e_{\bullet}\times e_{\bullet})^{\ast}(\Delta)&=
\sum_{f}B_{\beta_{1}}
e_{\bullet}^{\ast}(T_{f})\times B_{\bullet,
\beta_{2}}^{\prime}
 e_{\bullet}^{\ast}(T^f)\\
& =
\sum_{f}B_{\bullet, \beta_{1}}^{\prime}
 e_{\bullet}^{\ast}(B\cdot T_{f})\times B_{\bullet,
 \beta_{2}}^{\prime}
 e_{\bullet}^{\ast}(T^f)\,,
\end{split}
\end{equation}
we only have to sum over a basis for $A^{\ast}(X)/{\kker}i^{\ast}$.

Define a new pairing $<\,,\,>_{\B}$ on
$A^{\ast}(X)$ by
$\lav\gamma,\gamma'\rav_{\B}=\lav\gamma,B\cdot\gamma'\rav$.
This pairing is non-degenerate on
$A^{\ast}(X)/{\kker}i^{\ast}$, so we can choose bases $\{T_{f}\}$ and $\{^e
T\}$ for  $A^{\ast}(X)/{\kker}i^{\ast}$ such that
$\lav T_{f},^e T\rav_{\B}=\delta_{ef}$. Then \eqref{e:e4} becomes
\[\sum_{f}B_{\beta_{1}}e_{\bullet}^{\ast}(T_{f})\times
B_{\beta_{2}}e_{\bullet}^{\ast}({}^fT)\,,\]
and with \eqref{e:e5} we have the following splitting principle
for the new correlators:
\[
\varphi_{AB}^{\ast}(\tau_{d_1}\gamma_1,\ldots,\tau_{d_n}
\gamma_n)_{0,d}^{\B}
=\sum_{\beta_{1}+\beta_{2}=d}\sum_{f}
(\otimesover{a\in A}
\tau_{d_a}\gamma_a,T_{f})_{0,\beta_{1}}^{\B} ({}^f T,
\otimesover{b\in B}
\tau_{d_{b}}\gamma_{b})_{0,\beta_{2}}^{\B}\,. \]
This means that the results of Sections~\ref{gravdesc}--\ref{QDmodule}
are true if we replace  $(\;)_{0,d}$
and $<\,,\,>$ by $(\;)_{0,d}^{\B}$ and $<\,,\,>_{\B}$. In particular, 
if we
define the partial derivatives $\ll\,\gg_{0}^{\B}$ and
$<\,>_{0}^{\B}$
as in Section~\ref {potential} we have the following  structures:
\\\\
{\em WDVV equations.} Theorem~\ref{th2-1} and
Theorem~\ref{thm} hold for  $\ll\,\gg_{0}^{\B}$ and
$<\,>_{0}^{\B}$.
\\\\
$Q H^{\ast}(\B)$. Define products ${\ast}_{\B}$
on $A^{\ast}(X)\otimes \QQ[[q,q^{-1},t]]$ and
$A^{\ast}(X)\otimes \QQ[[q,q^{-1}]]$ by

\[T_{i}\ast_{\B}T_{j}=\sum_{k}\lav
T_{i},T_{j},T_{k}\rav_{0}^{\B}{^k T}\]
 and
\[T_{i}\ast_{\B}T_{j}=\sum_{d,k}\lav
T_{i},T_{j},T_{k}\rav_{0,d}^{\B}q^d\,{^k T}\,,\]
respectively, where $q^d$ still denotes
$\prod q_{i}^{\lav p_{i},d\rav}$ with
the old pairing. This makes
\[QH_{b}(\B)=A^{\ast}(X)\otimes\QQ[[q,q^{-1},t]]\]
 and
 \[QH^{\ast}(\B)=
A^{\ast}(X)\otimes\QQ[[q,q^{-1}]]\]
into commutative and
associative rings. The ring $Q H^{\ast}(\B)$ is graded if
we let
\[\deg q_{i}=\lav T^i,c_{1}(X)-c_{1}(\B)\rav_{X}\,.\]

\medskip

\noindent{\em Remark.} {\em $QH^{\ast}(Y)^{{}^{\perp}}$.}
We can split
\[A^{\ast}(Y)=i^{\ast}A^{\ast}(X)\oplus
i^{\ast}A^{\ast}(X)^{{}^{\perp}}\,,\]
where $\perp$ is with respect to the Poincar{\'e} pairing
$<\,,\,>_{Y}$ on $A^{\ast}(Y)$. Let 
$p_{i}^{\prime}$ be generators of $A_{\QQ}^1(Y)$
 with exponential coordinates $q_{i}^{\prime}$. Let
$\perp\:A^{\ast}(Y)\to
i^{\ast}A^{\ast}(X)$
denote the projection map. 
The product ${\ast}_{Y}^{{}^\perp}=\perp\circ{\ast}_{Y}$
defines a ring structure on
$QH^{\ast}(Y)^{{}^{\perp}}=i^{\ast}A^{\ast}(X)\otimes
\QQ[[q',q^{\prime-1}]]$. The arguments that produced 
associativity of $QH^{\ast}(\B)$ also show that 
$QH^{\ast}(Y)^{{}^{\perp}}$ is associative.
If $i^{\ast}\:A_{\QQ}^1(X)\to
A_{\QQ}^1(Y)$ is an isomorphism, then  $i^{\ast}$ induces an isomorphism
between $QH^{\ast}(\B)/{\ker i^{\ast}}$ and $QH^{\ast}(Y)^{{}^{\perp}}$.

A word of caution: There is no reason to believe that
$QH^{\ast}(Y)^{{}^{\perp}}$ is associative in general when
$Y\subseteq X$ without $Y$ being a zero scheme of a bundle.
\\\\
{\em Quantum $\D$-modules}. As in
\ref{Dubrovin}--\ref{QDmodule} there is a flat connection
$\nabla^{\B}$ defined by $\ast_{\B}$, with flat sections
\[\vec{s}_{a}^{\,\B}=e^{t_{0}/\hbar}
\!\!\left(\!e^{pt/{\hbar}}T_{a}\!+\!\sum_{d\ne 0}q^{d}
\sum_{b}{^b}T\int_{[\smap_{0,2}(X,\beta)]}\!
\frac{e_{1}^{\ast}(e^{pt/{\hbar}}T_{a})}{{\hbar}-c}
e_{2}^{\ast}(T_{b})
B_{d}\!\right)\,,\]
and a quantum $\D$-module $Q\D(\B)$ generated by
\[s^{\B}=\sum_{a}\lav\vec{s}_{a}^{\B},1\rav_{\B}{}^a T=
e^{(t_{0}+pt)/\hbar}\left(B+\sum_{d\ne 0}q^d
e_{1\ast}(\frac{B_{d}}{\hbar-c})\right)\,,\]
where $e_{1}\:\smap_{0,2}(X,d)\to X$ and $c=c_{1}(\L_{1})$.

\medskip

\noindent{\em Remark.}
Similar, corresponding to the product ${\ast}_{Y}^{{}^\perp}$ on
$QH^{\ast}(Y)^{{}^{\perp}}$ there is the quantum $\D$-module
$Q\D^{{}^{\perp}}(Y)$ generated by
\[s_{Y}^{{}^{\perp}}=
e^{(t_{0}^{\prime}+p't')/\hbar}\left(
1+\sum_{d\ne 0}q^{\prime d}
(\perp\circ e_{1})_{\ast}(\frac{1}{\hbar-c})
\right)\]
where $e_{1}\:\smap_{0,2}(Y,d)\to Y$ and $c=c_{1}(\L_{1})$.
Note that $Q\D^{{}^{\perp}}(Y)\supseteq
Q\D(Y)\,.$

\bigskip

\subsection{Quantum Lefschetz hyperplane principle}\label{QLefschetz}

Suppose
$\B=\oplus\oh(L_{i})$ is a $\rank k$
decomposable vector bundle generated by global sections, and
with Chern  roots $L_{i}\in A^1(X)$. The
classical Lefschetz hyperplane theorem states that
$i^{\ast}:A_{\QQ}^{l}(X)\to A_{\QQ}^{l}(Y)$
 is an isomorphism for $l\le\dim X-(k+1)$ and injective
for $l=\dim X-k$.

The quantum
Lefschetz hyperplane  principle relates $Q\D(X)$ and
$Q\D^{{}^{\perp}}(Y)$.
Let
\begin{equation}
\label{Lef}
I=e^{(t_{0}+pt)/\hbar}\left(1+\sum_{\beta\ne 0}
q^{\beta}\prod_{j=1}^k\prod_{m=1}^{\lav L_{j},\beta\rav}
(L_{j}+m\hbar)e_{1\ast}(\frac{1}{\hbar-c_{1}})\right)\,,
\end{equation}
and let
$q_{i}$ be of degree $\deg q_{i}=\lav
T^i,c_{1}(TX)-c_{1}(\B)\rav$ in this expression.
The quantum Lefschetz hyperplane principle states that the
$\D$-module generated by $i^{\ast}I$ and
$s_{Y}^{{}^{\perp}}$ are equivalent. More precisely: $i^{\ast}I$ and
$s_{Y}^{{}^{\perp}}$ coincide up to a unique weighted homogeneous
change of variables
\begin{equation*}
\begin{split}
&t_0^{\prime}\mapsto t_0+f_0\hbar+f_{-1}\\
&t_i^{\prime}\mapsto t_i+f_i
\end{split}
\end{equation*}
where $f_{-1},\ldots,f_r$ are weighted homogeneous power series
of
$q_{1},\ldots,q_r$, ${\deg} f_i=0$ for all $i\ge 0$ and ${\deg}
f_{-1}=1$.
The change of variables is determined by the
coefficients of $1=\left(\frac{1}{\hbar}\right)^0$ and
$\frac{1}{\hbar}$ in the expansion of $i^{\ast}I$ and
$s_{Y}^{{}^{\perp}}$ as power series of $\frac{1}{\hbar}$.

A version of this principle was first conjectured in
\cite{BvS}. This principle has been proven by Givental \cite{G} for convex
toric manifolds and by Kim \cite{Kim} for homogeneous manifolds. In
fact the proofs hold for all manifolds acted upon by a torus
with a finite number of fixpoints and one-dimensional orbits, 
but then only up to the existence of a
virtual equivariant fundamental class as defined in \cite{G}.

We conclude this part by illustrating how this formalizes the
mirror calculation of \cite{CdOGP} for the quintic in $\PP^4$.

\medskip

\begin{ex} \label{exU}
{\em (Mirror calculation)

Suppose $Y\subseteq\PP^4$ is a general quintic
hypersurface. The solutions to the differential equations
$\hbar\frac{d}{dt_{0}}-1=0$ and
$(\hbar\frac{d}{dt})^5-q=0$ are explicitly known
\cite{G} as
\begin{equation*}
s_{\PP^4}=e^{(t_{0}+pt)/\hbar}\biggl(
1+\sum_{d=1}^{\infty}q^d\frac{1}{\prod_{m=1}^d(p+m\hbar)^5}\biggr)\,.
\end{equation*}
For the quintic, \eqref{Lef} becomes
\begin{equation}
\label{mi1}
\begin{split}
i^{\ast}I=&e^{(t_{0}+p't)/\hbar}\biggl(
1+\sum_{d=1}^{\infty}q^d\frac{\prod_{m=1}^{5d}(5p'+m\hbar)}
{\prod_{m=1}^d(p'+m\hbar)^5}\biggr)\\
=&e^{t_{0}/\hbar}\bigl(I_{0}+\frac{I_{1}}{\hbar}p'+
\frac{I_{2}}{\hbar^2}p^{\prime
2}+\frac{I_{3}}{\hbar^3}p^{\prime 3}\bigr)\,.
\end{split}
\end{equation}
where $i^{\ast}p=p'$ and $I_{i}\in\CC[[q]][t]$. On the other
hand,
$(\hbar\frac{d}{dt'})^2\frac{1}{K(q')}(\hbar\frac{d}{dt'})^2=0$
has explicit solutions
\begin{equation}
\label{mi2}
\begin{split}
s_{Y}^{\perp}=&e^{t_{0}^{\prime}/\hbar}\biggl(
e^{p't'/\hbar}+\frac{p^{\prime 2}}{\hbar^{2}5}\sum_{d=1}^{\infty}n_{d}d^3
\sum_{d=1}^{\infty}\frac{e^{p't'/\hbar}q^{\prime d}}{(p'/\hbar+kd)^2}\biggr)\\
=&e^{t_{0}^{\prime}/\hbar}\biggl(1+\frac{t'}{\hbar}p'+\cdots\biggr)\,,
\end{split}
\end{equation}
where $K(q')=1+\frac{1}{5}\sum_{d=1}^{\infty}n_{d}d^3
\frac{q^{\prime d}}{1-q^{\prime d}}$. The quantum Lefschetz
hyperplane principle identifies \eqref{mi2} with
\[e^{t_{0}^{\prime}/\hbar}\biggl(\frac{I_{0}}{I_{0}}+
\frac{I_{1}}{I_{0}}\frac{p'}{\hbar}+
\frac{I_{2}}{I_{0}}\frac{p^{\prime 2}}{\hbar^2}+
\frac{I_{3}}{I_{0}}\frac{p^{\prime 3}}{\hbar^3}\biggr)\]
if we put $t'=\frac{I_{1}}{I_{0}}$. But \eqref{mi1}, with $\hbar=1$ is
a complete set of solutions to the Picard-Fuchs equation
\begin{equation*}
\left(\frac{d}{dt}\right)^4I=5q\bigl(5\frac{d}{dt}+1\bigr)
\bigl(5\frac{d}{dt}+2\bigr)
\bigl(5\frac{d}{dt}+3\bigr)\bigl(5\frac{d}{dt}+4\bigr)\,.
\end{equation*}
This proves the mirror computation of \cite{CdOGP}.
}\end{ex}

\addtocontents{toc}{\vspace{0.2cm}}
\section{Trisecant planes of the Veronese surface in $\PP^5$}
\label{trisecant}

\bigskip

In this section we introduce our main object of study,
$\NN$, the space of trisecant planes to the Veronese surface
in $\PP^5$, which is the same as the space of determinantal
nets of conics. A related space has been studied in
\cite{EPS,ES1,ES2}. We give a summary of
results found there, modified to our case. In
Section~\ref{quotient}, $\NN$ is constructed as a
geometric quotient, and in
Sections~\ref{vecbundles}--\ref{multstructure} we look at
the geometry of $\NN$. In particular, Ellingsrud and Str\o
mme's derivation of the Chow ring of $\NN$ is covered.

\bigskip

\subsection{The quotient construction}\label{quotient}

Let $V$ be a $\CC$-vector space of dimension $3$. We shall
study the space of $3\times 2$ matrices with
entries in $V$, modulo row and column operations. Let $E$ and $F$
be $\CC$-vector spaces of dimension $3$ and $2$ respectively.
Let the
group $GL(E)\times GL(F)$ act on $\Hom(F,E\otimes V)$ by
\[(g,h)\alpha=(g\otimes id_V)\circ\alpha\circ h^{-1}\,,\]
i.e., so that,
\begin{equation*}
\begin{CD}
                F@>{\alpha}>>E\otimes {V} \\
                @VhVV @VVg\otimes id_{V}V \\
                F @>{(g,h)\alpha}>>E\otimes {V}
        \end{CD}
\end{equation*}
commutes for $(g,h)\in GL(E)\times GL(F)$ and $\alpha\in
\Hom(F,E\otimes V)$. The normal subgroup
\[\Gamma=\{k(id_E,id_F)\mmid k\in\CC^*\}\]
acts trivially on $\Hom(F,E\otimes V)$.
Let $G=GL(E)\times GL(F)\largefrac\Gamma$ and consider the
induced action of $G$ on $\Hom(F,E\otimes V)$. A point of
$\Hom(F,E\otimes V)$ will be called {\em stable}
(resp. {\em semistable}) if the corresponding point in
$\PP(\Hom(F,E\otimes V)^{\vee})$ is {\em stable} (resp. {\em semistable})
in the sense of \cite{MF} for the induced action of $G\cap
SL(\Hom(F,E\otimes V))$.

Suppose $A=(L_{ij})$ is a $3\times 2$ matrix representing
$\alpha$,
with $L_{ij}\in V$. Let $E_{\alpha}\subseteq S_{2}V$ be the span
of all $2\times 2$ minors of $A$. Then $E_{\alpha}$ does not depend on the
matrix representative $A$.

\begin{prop}
\label{propstable} \cite{EPS}.
For any $\alpha\in\Hom(F,E\otimes V)$ the following are equivalent:
\begin{enumerate}
\item[i)] $\alpha$ is stable
\item[ii)] $\alpha$ is semistable
\item[iii)] $\dim E_{\alpha}=3$
\item[(iv)] $\alpha$ cannot be represented by a matrix of type
$\bigl(\begin{smallmatrix} 0&\ast\\
                 0&\ast\\
                 \ast&\ast
   \end{smallmatrix}\bigr)$
or
$\bigl(\begin{smallmatrix} 0&0     \\
                 \ast&\ast\\
                 \ast&\ast
   \end{smallmatrix}\bigr)$.
\end{enumerate}
\end{prop}

By geometric invariant theory  there exists a projective
geometric quotient \[\NN={\Hom(F,E\otimes V)}^{ss}\largefrac
G\,.\]

\begin{prop}
\label{propfree}
\cite{EPS}. The group $G$ acts freely on $\Hom(F,E\otimes V)^{ss}$,
hence $\NN$ is smooth.
\end{prop}

Proposition~\ref{propstable} implies that there is a morphism
\begin{equation}
\label{e:p2-1}
\iota:\NN\to{\Grass}_{3}(S_{2}V)
\end{equation}
given by $\iota(\alpha)=E_{\alpha}$, where $(\alpha)$ is the class
of $\alpha\in\Hom(F,E\otimes V)^{ss}$.

Let ${\Hilb}_{\PP(V)}^3$ be the Hilbert scheme of subschemes
of length $3$ in
$\PP(V)$. There is a morphism
\[b:{\Hilb}_{\PP(V)}^3\to{\Grass}_{3}(S_{2}V)\]
given by $b([Z])=H^0(\PP(V),\,\I_Z(2))$ for $Z\subseteq\PP(V)$ a
zero-dimensional subscheme of length 3. Using the Hilbert-Burch
theorem \cite{E}, one may show that $\iota$ is an embedding
and that the image of $b$ can be identified with
$\NN\subseteq\Grass_{3}(S_{2}V)$. In fact
\begin{equation}
\label{e:p2-11}
b:{\Hilb}_{\PP(V)}^3\to \NN
\end{equation}
blows up the locus of nets of type $\lav
x_0^2,x_0x_1,x_0x_2\rav\subseteq S_{2}V$, where
$\{x_0,x_1,x_2\}$ is a basis for $V$.

\medskip

\noindent{\em Remark~1.}
The variety $\NN$ is often denoted $N(3;2,3)$. As the
notation indicates it belongs to a bigger class of
varieries $N(q;n,m)$, which parametrize $n\times m$ matrices
with entries from a $q$-dimensional vector space, modulo
row and column operations. Topological properties of $N(q;n,m)$
have been studied in \cite{Dr}. In \cite{EPS,ES1} the
variety $N(4;2,3)$ is considered. It is related to the
component of the Hilbert scheme containing the twisted cubics
in an analog manner to \eqref{e:p2-11}.

\medskip

\noindent{\em Remark~2.}
If $\Grass_{3}(S_{2}V)$ is seen as the variety of 2-planes in
$\PP(S_{2}V)$ then \eqref{e:p2-1} identifies
$\NN$ as the variety parametrizing trisecant planes to the
Veronese surface in $\PP(S_{2}V)$.

\bigskip

\subsection{Universal vector bundles}
\label{vecbundles}

Let $E'$ and $F'$ denote the trivial vector bundles
\begin{equation*}
E\otimes\frac{\det F}{\det E}\times\Hom(F,E\otimes V)^{ss}
\quad{\rm and}\quad F\otimes\frac{\det F}{\det E}\times\Hom(F,E\otimes V)^{ss}
\end{equation*}
on $\Hom(F,E\otimes V)^{ss}$. The bundles $E'$
and $F'$ along with the tautological map $F'\to E'\otimes V$ are
$G$-equivariant hence descend to bundles $\E$ and $\F$ on
$\NN$ and a map $\A:\F\to\E\otimes V$ of vector
bundles.

The embedding $\iota:\NN\to\Grass_{3}(S_{2}V)$ identifies the bundles
$\E=\iota^{\ast}\U$, where $\U$ is the tautological subbundle on
$\Grass_{3}(S_{2}V)$,
and the map $\A$ sits inside the complex
\begin{equation}
\label{diag31}
0\to\F\xrightarrow{\A}\E\otimes V\xrightarrow{m}S_{3}V\,,
\end{equation}
where $m$ is the multiplication map. In fact,
\eqref{diag31} is exact outside the locus of nets of type
$\lav x_{0}^2,x_{0}x_{1},x_{0}x_{2}\rav$, hence, by analytic
continuation, the complex uniquely determines
$\F$, in the sense that if $\G$ is a rank two bundle on
$\NN$, and if $\G\to\E\otimes V$ is a map of vector bundles
such that $0\to\G\to\E\otimes V\xrightarrow{m}S_{3}V$ is a
complex, then $\F$ and $\G$ are isomorphic.

The differential of the quotient map $\Hom(F,E\otimes
V)^{ss}\to\NN$ gives rise to an exact sequence
\begin{equation}
\label{e:e5-1-1}
0\to\oh_{\NN}\to\CEnd(\E)\oplus\CEnd(\F)\to\HomC(\F,\E\otimes V)\to
T\NN \to 0\,,
\end{equation}
where local sections $(g,h)$ of $\CEnd(\E)\oplus\CEnd(\F)$ are
sent to $(g\otimes id_V)\circ\A {\rm -}\A\circ h$.

\bigskip

\subsection{A torus action}
\label{torusaction}
 Let $T\subset GL(V)$ be a maximal torus and consider the action of
$T$ on $V$. It induces an action on $\NN$ and the bundles
$\E$ and $\F$. Let $\{x_0,x_1,x_2\}$ be a basis for $V$ such that the action of
$T$ on $V$ is diagonal, and let $\lambda_i$ the corresponding
characters so that $t\cdot x_i=\lambda_i(t)x_i$ for all
$t\in T$. The action on $\PP(V)$ has isolated fixpoints
$\{(x_0,x_1),(x_0,x_2),(x_1,x_2)\}$ and
isolated fixlines $\{(x_0),(x_1),(x_2)\}$. The blow-up map
\eqref{e:p2-11} is $T$-equivariant, hence if $\alpha\in \NN$ is a
fixpoint, the locus of base points
$B(\E_{\alpha})$ of the net $\E_{\alpha}$
must be supported on the triangle of fixpoints and fixlines in $\PP(V)$.
If $B(\E_{\alpha})$ is not a line, then it is projectively
equivalent to one of the following types: three distinct
points, the union of a point doubled on a line and another
point on the line, or the full first-order neighborhood of
a point.
\begin{table}[h!]\centering
\caption{Fixpoints\label{tablefixpoints}}
\begin{tabular}{||c|c|c|c|c||}
\hline\hline
&&&&\\
Matrix&$\E_{\alpha}$&$\F_{\alpha}$&Isotropy&$B(\E_{\alpha})$\\
       &&&&\\
       \hline
       &&&&\\
$\begin{pmatrix} x_1&0\\ x_2&x_2\\ 0&x_0 \end{pmatrix}$
       &$\begin{array}{c} x_1 x_2\\ x_0x_1\\
       x_0x_2\end{array}$
       &$\begin{array}{c}x_{1}x_{2}\otimes x_{0}-x_{0}x_{1}
       \otimes x_{2}\\
x_{0}x_{1}\otimes x_{2}-x_{0}x_{2}\otimes x_{1}
\end{array}$
       &$S_{3}$
       &$\tttheref$
       \\
       &&&&\\
       \hline
       &&&&\\
$\begin{pmatrix} x_1&0\\ x_0&x_1\\ 0&x_2\end{pmatrix}$
        &$\begin{array}{c} x_1^2\\ x_1x_2\\ x_0x_2\end{array}$
        &$\begin{array}{c}x_{1}^2\otimes x_{2}-x_{1}x_{2}\otimes
        x_{1}\\
x_{1}x_{2}\otimes x_{0}-x_{0}x_{2}\otimes x_{1}
\end{array}$
        &$\{id\}$
        &$\rrighttheref$
        \\
       &&&&\\
       \hline
       &&&&\\
$\begin{pmatrix} x_1&0\\ x_1&x_2\\ 0&x_2\end{pmatrix}$
        &$\begin{array}{c}x_1^2\\ x_1x_2\\ x_2^2\end{array}$
        &$\begin{array}{c}x_{1}^2\otimes x_{2}-x_{1}x_{2}\otimes
        x_{1}\\
x_{1}x_{2}\otimes x_{2}-x_{2}^2\otimes x_{1}
\end{array}$
        &$x_{1}\leftrightarrow x_{2}$
        &$\odot$
        \\
        &&&&\\
        \hline
        &&&&\\
$\begin{pmatrix} x_1&0\\ 0&x_1\\ x_0&x_2\end{pmatrix}$
        &$\begin{array}{c}x_1^2\\ x_1x_2\\ x_0x_1\end{array}$
        &$\begin{array}{c}x_{1}^2\otimes x_{2}-x_{1}x_{2}
        \otimes x_{1}\\
x_{1}^2\otimes x_{0}-x_{0}x_{1}\otimes x_{1}
\end{array}$
        &$x_{0}\leftrightarrow x_{2}$
        &\rule[0.5ex]{2.0em}{0.08ex}
        \\
        &&&&\\
       \hline\hline
\end{tabular}
\end{table}

In Table~\ref{tablefixpoints} a list of the fixpoints in
$\NN$ is presented (up to a permutation of $(x_0,x_1,x_2)$). For each
fixpoint $\alpha$ the table contains a matrix
representative, generators for the vector spaces
$\E_{\alpha}$ and $\F_{\alpha}$, a description of
$B(\E_{\alpha})$, and the isotropy group of the
$S_{3}$-action on the fixpoints induced from permutation of
$(x_0,x_1,x_2)$.

\medskip

\noindent {\em Remark.} Up to projective equivalence there
are five type of points in $\NN$: the four in
Table~\ref{tablefixpoints}, and the net $\lav
x_{1}^2,x_{1}x_{2},x_{2}^2+x_{0}x_{1}\rav$, whose locus of base
points is a
point with a second-order direction.

Let $\gamma:\CC^{\ast}\to T$ be a one-parameter subgroup
such that $\lambda_i(\gamma(t))=t^{\omega_i}$ for integer weights
$\{\omega_0,\omega_1,\omega_2\}$, and consider the induced
$\CC^{\ast}$-action on $\NN$. The weights can be chosen
such that the fixpoints of the $\CC^{\ast}$-action are the
same as the fixpoints of the $T$-action. Since the fixpoints are
isolated we can apply the Bialynicki-Birula theorem
to $\NN$. Counting the fixpoints in
Table~\ref{tablefixpoints} we find $13$, and the weights of
the $T$-action on the tangent spaces at fixpoints can be
found by using \eqref{e:e5-1-1}. In fact, in
Section~\ref{calculations} we shall look closer at how such
computations are performed.

\begin{thm}\label{betti}\quad
\begin{enumerate}
\item[i)] $\chi(\NN)=13$
\item[ii)] The Betti numbers $b_{p}=\dim A^p(\NN)$ are
\[(b_{0},b_{1},\ldots,b_{6})=(1,1,3,3,3,1,1).\]
\end{enumerate}
\end{thm}

\medskip

\noindent{\em Remark.} In \cite{Dr} a recurrence formula
can be found for the Betti numbers of $N(q;n,m)$ when $n$
and $m$ are relatively prime.

\bigskip

\subsection{Multiplicative structure of $A^{\ast}(\NN)$}
\label{multstructure}

In \cite{ES2} a detailed computation of $A^{\ast}(N(4;2,3))$
can be found. The same
method works for our case.

\begin{prop}
The Chow ring $A^{\ast}(\NN)$ is generated as a $\ZZ$-algebra by
the Chern classes of $\E$ and $\F$. Their monomial
values in $A^6(\NN)$ are (with $\gamma_i=c_i(\E)$ and
$\delta_i=c_i(\F)$) as in Table~\ref{tablemonomialvalues}.
\end{prop}
\begin{table}[h!]\centering
\caption{Monomial values\label{tablemonomialvalues}}
\begin{tabular}{||c|c|c|c|c|c|c||}
\hline\hline
&&&&&&\\
$\!\gamma_1^4\gamma_2\!=\!27\!$ &
$\!\gamma_1^4\delta_2\!=\!18\!$ &
$\gamma_1^3\gamma_3\!=\!5$ &
$\!\gamma_1^2\gamma_2^2\!=\!14\!$&
$\!\gamma_1^2\gamma_2\delta_2\!=\!9\!$&
$\gamma_1^2\delta_2^2\!=\!6$&
$\!\gamma_1\gamma_2\gamma_3\!=\!3\!$
\\&&&&&&\\\hline&&&&&&\\
$\!\gamma_1\gamma_3\delta_2\!=\!2\!$&
$\gamma_2^3\!=\!9$ &
$\gamma_2^2\delta_2\!=\!5$ &
$\gamma_2\delta_2^2\!=\!3$ &
$\!\gamma_3^2\!=\!1\!$ &
$\gamma_1^6\!=\!57$ &
$\delta_2^3\!=\!2$\\&&&&&&\\
\hline\hline
\end{tabular}
\end{table}

\medskip

\noindent{\em Remark.}
 Since $\wedge^3 E'\simeq\wedge^2 F'$ as $G$-bundles,
$c_1(\E)=c_1(\F)$.

\begin{proof}
First, we find some generators of $A^{\ast}(\NN)$, then we
look at relations among them, and finally we describe the class of a
point.\\
{\em Generators.}
Let
\[\eta\:P=\Isom(\E,E)\times_{\NN}\Isom(\F,F)\to\NN\]
be the principal $GL(E)\times GL(F)$-bundle associated to $\E$ and
$\F$. It is a result of Grothendieck \cite{Gr} that
$\eta^{\ast}:A^{\ast}(\NN)\to A^{\ast}(P)$ is
surjective with kernel generated by $c_i(\E)$ and $c_i(\F)$, $i>0$.
As $\eta$ factors through the quotient map $\Hom(F,E\times
V)^{ss}\to\NN$, and $A^{\ast}(\Hom(F,E\otimes V)^{ss})\simeq \ZZ$,
this shows that $c_i(\E)$, $c_i(\F)$ generate $A^{\ast}(\NN)$
over $\ZZ$.

\noindent{\em Relations.}
We find relations among the generators by displaying
$A^{\ast}(\NN)$ as the quotient of a polynomial ring.

First let $\FF=\FF(\E)\times_{\NN}\PP(\F)\xrightarrow{\rho}\NN$ be
the fiber product of the flag bundles with universal flags
\[\rho^{\ast}\E=\E^3\doublerightarrow\E^2\doublerightarrow\E^1
\doublerightarrow\E^0=0\]
and
\[\rho^{\ast}\F=\F^2\doublerightarrow\F^1\doublerightarrow\F^0=0\,.\]

Let $R=\ZZ[e_1,e_2,e_3,f_1,f_2]$ and
$R^S=\ZZ[\gamma_1,\gamma_2,\gamma_3,\delta_1,\delta_2]$ be
polynomial rings where $\gamma_i$ and $\delta_i$ are the elementary
symmetric polynomials in the $e_i$ and $f_i$ respectively. Note that
$R$ is a free $R^S$-module with basis
\[B=\{e_1^i e_2^j f_1^k\mid 0\le i\le 2,\, 0\le j,\, k\le 1\}.\]

Consider the following commutative diagram of ring homomorphisms:
\begin{equation}
\label{e:firk}
        \begin{CD}
          R^S @>\iota>>R\\
         @V\lambda VV @VV\mu V  \\
          A^{\ast}(\NN)@>\rho^{\ast}>> A^{\ast}(\FF)
        \end{CD}
\end{equation}
where:

\begin{enumerate}
\item[] $\rho^{\ast}$ is the injective
pullback map.
\item[]
$\iota$ is the inclusion map.
\item[]
$\lambda$ maps $\gamma_i$, $\delta_j$ to
the Chern classes $c_i(\E)$, $c_j(\F)$, respectively.
\item[] $\mu$ maps $e_i$, $f_j$ to
the Chern roots $c_1(\E^i)$--$c_1(\E^{i-1})$,
$c_1(\F^i)$--$c_1(\F^{i-1})$, respectively.
\end{enumerate}

For $r\in R$ let $c(r)\subseteq R^S$ be the ideal generated by the
coordinates of $r$ with respect to the basis in $B$. For $J\subset
R$ an ideal, let $c(J)$ be the ideal generated by $\{c(r)\mid r\in
J\}$. A short argument using the surjectivity of $\mu$ and
the commutativity of \eqref{e:firk} yields
$\ker(\lambda)=c(\ker(\mu))$.

\medskip

\noindent
{\bf Claim 1.}
$\ker(\mu)$ contains the following elements:
\begin{enumerate}
\item[i)] $r_0=e_1+e_2+e_3-f_1-f_2$
\item[ii)] $r_1=(e_1-f_1)^3(e_1-f_2)^3$
\item[iii)] $r_2=(e_1-f_2)^3(e_2-f_2)^3$,
\end{enumerate}

\begin{proof}
i) follows from $c_1(\E)=c_1(\F)$.

\noindent
For ii) and iii) consider the maps
  \begin{equation*}
\begin{CD}
          \rho^{\ast}\F   \\
         @VV\A V @. \\
          {\rho^{\ast}\E\otimes V}@>{}>>
         \kern-1.8em\hbox{$\rightarrow$}
         \kern1.0em\hbox{$\E^{1}\otimes V$}
        \end{CD}
\quad\, {\rm and}\quad\,
\begin{CD}
          \F^{(1)}@>>{}> \rho^{\ast}\F \\
         @. @VV\A V  @. \\
          @. {\rho^{\ast}\E\otimes V}@>{}>>
         \kern-1.8em\hbox{$\rightarrow$}
         \kern1.0em\hbox{$\E^{2}\otimes V$}
        \end{CD}
        \end{equation*}
where $\F^{(1)}=\ker(\rho^{\ast}\F\to\F^1)$. The stability condition
iv) of Proposition~\ref{propstable} prevents the maps
$\rho^{\ast}\F\to\E^1\otimes V$ and $\F^{(1)}\to
\rho^{\ast}\E^2\otimes V$ to have any
zeros, hence
\begin{eqnarray*}
&0=c_6(\rho^{\ast}\F^{\vee}\otimes\E^1\otimes V)=
\mu\left((e_1-f_1)^3(e_1-f_2)^3\right)\\
&0=c_6({\F^{(1)}}^{\vee}\otimes\rho^{\ast}\E^2\otimes V)=
\mu\left((e_1-f_2)^3(e_2-f_2)^3\right)\,.
\end{eqnarray*}
\end{proof}

\noindent
{\bf Claim 2.} Let ${\gotha}$ be the ideal generated by
$c(r_i)$, $i=0,1,2$.
Then
\[A_{\QQ}^{\ast}(\NN)=(R^S\largefrac{\gotha})\otimes_{\ZZ}\QQ\,.\]

\begin{proof}
See \cite{ES2}, Theorem~6.9.
\end{proof}

\noindent
{\em Class of a point.}
We need to express the class of a point in terms of the
generators $c_{i}(\E)$, $c_{i}(\F)$. By Gauss-Bonnet,
$13=\chi(\NN)=\deg c_{6}(T\NN)$, hence the class of a point is
$(1/13)c_{6}(T\NN)$.
A formal Chern class computation,
using the exact sequence \eqref{e:e5-1-1}, yields an expression
of desired kind for
$c_{6}(T\NN)$.

The numbers in Table~\ref{tablemonomialvalues} are obtained
from the above description using the {\em Schubert} package
in {\texttt{MAPLE}} \cite{KS}.
\end{proof}

\bigskip

\subsection{Calabi-Yau sections}
\label{cysections}
 From the resolution \eqref{e:e5-1-1} it follows that
$c_{1}(T\NN)=3p$
where $p=-c_1(\E)$ is ample, hence $\NN$ is Fano of index 3.

Consider
the following bundles:
\[3\oh(p),\quad \F^{\vee}\oplus\oh(2p),\quad
S_{2}\F^{\vee}\,.\]
All three are generated by their global sections and have $3p$
as their first Chern
class. Hence zero schemes of general global sections of
the above bundles are Calabi-Yau threefolds. We shall denote
these by $X$, $Y$, and $Z$, respectively.

Questions on the enumerative geometry of rational curves
on these Calabi-Yau´s, lead to questions about
$\moo(\NN,d)$. For the time being we simply
record that the expected dimension of this moduli space is
$3+3d$.

\addtocontents{toc}{\vspace{0.2cm}}
\section{Lines on $\NN$}
\label{lines}

 \bigskip

We show that the functor of stable maps
$\barr{\mathcal M}_{0,0}(\NN,1)$
 (\footnotemark)\
 is represented by a $\PP^2$-bundle over $\PP(V)\times\PP(V^{\vee})$.
 This enables us to find an explicit representation of the
 Chow ring $A^{\ast}(\moo(\NN,1))$ in terms of generators
 and relations. As an application we calculate the number
 of lines on the Calabi-Yau threefold sections.
\footnotetext{
 This is of course the same as the Hilbert functor of lines
 in $\NN$.}

 \bigskip

 \subsection{The moduli space}
 \label{modulspace}

 The idea is as follows:
 A line $L\subseteq\Grass_{3}(S_{2}V)$ is determined by a
 pencil $P\subseteq S_{2}V$ and a web $W\subseteq S_{2}V$
 such that
  $P\subseteq W$. The line $L$ parametrizes the
 family of nets $P\subseteq E_{t}\subseteq W$.
  From the relationship between $\NN$ and $\Hilb_{\PP(V)}^3$
 it follows that the line $L$ is
 contained in $\NN$ if and only if $\{B(E_{t})\mid t\in
 L\}$ describes a (non-constant) family of lines/three
 points in $\PP(V)$. As $B(W)\subseteq B(E_{t})\subseteq
 B(P)$ for all $t\in L$, this implies that $L\subseteq \NN$
 if and only if both i) and ii) hold:
 \begin{enumerate}
 \item[i)]$B(P)$ is one dimensional. As $B(P)$ is the
 intersection of two conics, $B(P)$ must be the union of a
 point $p$ and a line $l$, with the point possibly embedded
 onto the line.
 \item[ii)]$W\largefrac P$ is a pencil of conics through
 $p$, so $W\largefrac P$ corresponds to a point in
 \\
 $\PP(H^0({\I}_{p}(2))\largefrac H^0({\I}_{p\cup
 l}(2)))$.
  \end{enumerate}
 We show that this correspondence can be relativized.

 Let $PP^{\vee}=\PP(V)\times \PP(V^{\vee})$ and let
 $\PP=\PP(V)\times PP^{\vee}\xrightarrow{pr}PP^{\vee}$
 be the projection map. Let $\oh_{pr}(1)$ denote the
tautological line bundle. Consider the universal line and
point in $\PP(V)$
\begin{equation*}
\Tilde{l}=\{(x,p,l)\in\PP\mid x\in l\}\qquad {\rm and}\qquad
\Tilde{p}=\{(x,p,l)\in\PP\mid x= p\}\,,
\end{equation*}
and the exact sequence
\begin{equation}
0\to{\I}_{\Tilde{p}}\I_{\Tilde{l}}(2)\to
{\I}_{\Tilde{p}}(2)\to
{\I}_{\Tilde{p}}\otimes\oh_{\Tilde{l}}(2)\to 0\,.
\end{equation}
Since the sheaves in this sequence are flat, and
$H^1(\I_{p}\I_{l}(2))=H^1(\I_{p}(2))=H^1(\I_{p}\otimes\oh_{l}(2)) =0$
for arbitrary point-line pairs $(p,l)$ in $\PP(V)$, base
change theorems in \cite{H} imply that
\begin{equation}
\label{e:e6-1a}
0\to pr_{\ast}{\I}_{\Tilde{p}}\I_{\Tilde{l}}(2)\to
 pr_{\ast}{\I}_{\Tilde{p}}(2)
 \to  \V
 \to R^{1}pr_{\ast}{\I}_{\Tilde{p}}\I_{\Tilde{l}}(2)\,,
\end{equation}
with
$\V=pr_{\ast}({\I}_{\Tilde{p}}\otimes\oh_{\Tilde{l}}(2))$,
is an exact sequence of vector bundles. Further, the fibers
of $\V$ over point-line pairs $(p,l)$ are
$H^0(\I_{p}\otimes\oh_{l}(2))=H^0(\I_{p}(2))/H^0(\I_{p}\I_{l}(2))$.

\begin{thm}
\label{thmfunctor}
The functor of stable maps $\barr{\mathcal M}_{0,0}(\NN,1)$ is represented
by $\PP(\V)$. In particular, $\smap_{0,0}(\NN,1)$ is smooth and
of dimension $6$.
\end{thm}

We shall actually prove a more general version which will
be useful later. If $f\:\PP^1\to C$ is a parametrization of
a smooth curve $C$ of degree $d$ in $\NN$, then the pullback
of $\E$ decomposes into parts
\begin{equation}
\E_{\PP^1}\simeq\oh_{\PP^1}(-a_{1})\oplus\oh_{\PP^1}(-a_{2})
\oplus\oh_{\PP^1}(-a_{3})
\end{equation}
such that $d=\sum a_{i}$ and $a_{i}\ge 0$ for all $i$. Let
$H_{d}$ denote the components of $\Hilb_{\NN}^{dt+1}$
containing the set of rational curves of degree $d$, and
let $H^{0,d}$ denote the closure in $H_{d}$ of the locus of
smooth curves with splitting type $(0,0,d)$.
The scheme $H^{0,d}$ comes equipped with a universal family
of degree $d$ curves $\C_{d}\subseteq H^{0,d}\times\NN$.
In particular, $H^{0,1}=\moo(\NN,1)$ and $\C_{1}=\smap_{0,1}(\NN,1)$.

Consider the Hilbert scheme $\PP(S_{d}\V)$ parametrizing
flat families of degree $d$ curves in $\PP(\V^{\vee})$
over $PP^{\vee}$ and let
$\C_{d}^{\prime}\subseteq\PP(S_{d}\V)\times_{PP^{\vee}}
\PP(\V^{\vee})$ be the universal family. In the proof of the
theorem we show that $\PP(S_{d}\V)$ and $H^{0,d}$ are
naturally isomorphic when $d=1,2$.

Before we start, we introduce some more notations: Let
\[\rho\:\PP(S_{d}\V)\to PP^{\vee}\qquad{\rm and}\qquad
\rho^{\vee}\:\PP(\V^{\vee})\to PP^{\vee}\] denote the
projection maps, and let $\oh(\check{\sigma})$ denote the
tautological line bundle on $\PP(\V^{\vee})$.

\begin{proof}
The morphism
\begin{equation}
\label{e:e6-2}
\PP(\V^{\vee})\to\NN\,,
\end{equation}
which maps
$\V_T^{\vee}\to\oh_T(\check{\sigma})$ to nets
$pr_{\ast}{\I}_{\Tilde{p}}\I_{\Tilde{l}}(2)_T\oplus
\oh_T(-\check{\sigma})$ is surjective, embeds each
fiber $\PP(\V_{(p,l)}^{\vee})$ onto planes in $\NN$, and induces
morphisms
\begin{equation}
\label{e:e6-3}
\PP(S_{d}\V)\to H^{0,d}
\end{equation}
for $d=1,2$. From our initial analysis of the lines in $\NN$,
it follows that the morphism \eqref{e:e6-3} is surjective
for $d=1$. This
argument can easily be extended to the case $d=2$.

Consider the induced morphisms
$\PP(S_{d}\V){\times}_{PP^{\vee}}\PP(\V^{\vee})
\to H^{0,d}\times\NN$.
To prove that \eqref{e:e6-3} is an isomorphism we must show that the
universal family $\C_{d}$ lifts to the
universal family
$\C_{d}^{\prime}$.
We shall prove that the images $P_1\subseteq\NN$ and
$P_2\subseteq\NN$ of two different fibers
$\PP(\V_{(p_1,l_1)}^{\vee})$ and $\PP(\V_{(p_2,l_2)}^{\vee})$ are
such that either $P_1\cap P_2$ is empty or $P_1\cap P_2$ is a
point. It then follows that $\C_{d}^{\prime}$ is the unique component of
dimension $\dim\C_{d}$ in the preimage of $\C_{d}$.

If $[N]$ is a point in $P_i$ then either $B(N)$ is the line $l_i$,
or $B(N)$ is supported on $p_i$ and a length two subscheme of
$l_i$. Keeping this scenario in mind, it is easy to deduce
the following (assuming $(p_1,l_1)\ne(p_2,l_2)$):

If $l_1=l_2$,then $P_1\cap P_2$ is empty unless $p_1\in
l_1$ and $p_2\in l_2$. If so, $[N]\in P_1\cap P_2$ if and only if
$B(N)$ is the line $l_1=l_2$, hence $P_1\cap P_2$ is a point.

 If $l_1\ne l_2$, then $P_1\cap P_2$ is empty unless
$p_1\in
l_2$ and $p_2\in l_1$. If so, there are three cases to consider:
\begin{enumerate}
\item[i)] $p_1\not\in l_1$ and $p_2\not\in l_2$,
\item[ii)] $p_1\not\in l_1$ and $p_2\in l_2$,
\item[iii)] $p_1=p_2=l_1\cap l_2$.
\end{enumerate}
In all three cases, if $[N]\in P_{1}\cap P_{2}$, then $B(N)$
must contain the point $l_{1}\cap l_{2}$.
Without loss of generality we may assume that
$l_{1}=(x_{1})$, $l_{2}=(x_{2})$, $p_{1},p_{2}\in (x_{0})$ in
i), and $p_{1}\in (x_{0})$ in ii). It is straightforward to
verify that in each case $P_{1}\cap P_{2}$ consists of only
one point. These represent the nets:
\[
{\rm i)}\,\lav x_{0}x_{1},x_{1}x_{2},x_{0}x_{2}\rav\,,\qquad
{\rm ii)}\,\lav x_{0}x_{1},x_{1}x_{2},x_{2}^2\rav\,,\qquad
{\rm iii)}\,\lav x_{1}^2,x_{1}x_{2},x_{2}^2\rav\,.\]
\end{proof}

\noindent{\em Remark 1.} Note that if $P_{1}$ and $P_{2}$
is as in the proof, and their intersection is non-empty, then
the union $P_{1}\cup P_{2}$ is determined up to
projective equivalence by the intersection point,
which again must be a point of type 1)-4) in
Table~\ref{tablefixpoints}.

\medskip

\noindent{\em Remark 2.} On $\PP(\V^{\vee})$ there
is a natural $T$-equivariant complex
\begin{equation}
0\to pr_{\ast}\I_{\Tilde{l}}(1)\otimes\bigl(\wedge^2
pr_{\ast}\I_{\Tilde{p}}(1)\oplus\oh(-\check{\sigma})\bigr)\to
\bigl(pr_{\ast}\I_{\Tilde{p}}\I_{\Tilde{l}}(2)\oplus
\oh(-\check{\sigma})\bigr)\otimes V\xrightarrow{m} S_{3}V\,,
\end{equation}
so the morphism \eqref{e:e6-2} identifies the bundles
\begin{equation}
\begin{split}
\E_{\PP(\V^{\vee})}=&
pr_{\ast}\I_{\Tilde{p}}\I_{\Tilde{l}}(2)\oplus\oh(-\check{\sigma})
\\
\F_{\PP(\V^{\vee})}\simeq& pr_{\ast}\I_{\Tilde{l}}(1)\otimes
\bigl(\wedge^2
pr_{\ast}\I_{\Tilde{p}}(1)\oplus\oh(-\check{\sigma})\bigr)\,.
\end{split}
\end{equation}

\bigskip

\subsection{The Chow ring $A^{\ast}(\moo(\NN,1))$}
\label{chowen}

Let $\oh(\tau)$ and $\oh(\check{\tau})$ be the tautological
line bundles on $\PP(V)$ and $\PP(V^{\vee})$, and let
$\oh(\sigma)$ be the tautological bundle on $\PP(\V)$ . Then
\[A^{\ast}(\PP(V)\times\PP(V^{\vee}))=\ZZ[\tau,\check{\tau}]
\largefrac
(\tau^3,\check{\tau}^3)\,.\]
Since $\moo(\NN,1)$ is a projective bundle we have, by
Grothendieck's theorem \cite{F},
\begin{equation}
A^{\ast}(\moo(\NN,1))=\ZZ[\tau,\check{\tau},\sigma]
\largefrac
(\tau^3,\check{\tau}^3,r)\,,
\end{equation}
where $r=\sum_{i=0}^3(-1)^i\rho^{\ast}c_{i}(\V)\sigma^{3-i}\,.$

We shall determine the total Chern class $c(\V)$ in terms
of the generators $\tau$ and $\check{\tau}$.

Let $p_{1}\:\PP\to\PP(V)$ and $p_{2}\:\PP\to\PP(V)$ be the
projection maps on the first and second components. On $\PP$
there are exact sequences
\begin{equation}
\label{diagonal}
        \begin{CD}
        @. @. 0@.@. \\
        @. @. @V{}VV     \\
        {} @. @. p_{1}^{\ast}\Omega(1)@. \\
        @. @. @V{}VV     \\
          0 @>{}>> p_{2}^{\ast}\Omega(\tau)
           @>{\alpha}>> V_{\PP} @>{}>> \oh(\tau)@>{}>> 0 \\
         @. @. @V{\beta}VV     \\
          {} @. @. \oh_{pr}(1) \\
           @. @. @V{}VV     \\
           {} @. @. 0
        \end{CD}
\end{equation}
where $\Omega$ is the sheaf of K\"ahler differentials on
$\PP(V)$. The map
$s=\beta\circ\alpha\:p_{2}^{\ast}\Omega(\tau)\to\oh_{pr}(1)$
induces a section of
$\Omega(\tau)^{\vee}\otimes\oh_{pr}(1)$ whose zero scheme
is $\Tilde{p}$. Since $s$ is
regular \cite{F}, the twisted Koszul complex
\begin{equation}
0\to\wedge^2 \rho_{2}^{\ast}\Omega(\tau)\otimes
\oh_{pr}(-1)\to
\rho_{2}^{\ast}\Omega(\tau)\to \I_{\Tilde{p}}(1)\to 0
\end{equation}
is exact. Pushing down this sequence we find
$pr_{\ast}\I_{\Tilde{p}}(1)\simeq\Omega(\tau)$, and twisting
it by $\oh_{pr}(1)$ before pushing down yields
an exact sequence
\begin{equation}
\label{n4-2a}
0\to\wedge^2\Omega(\tau)\to\Omega(\tau)\otimes V\to
pr_{\ast}\I_{\Tilde{p}}(2)\to 0\,.
\end{equation}

Using a similar diagonal argument as in \eqref{diagonal},
one can show that
$\I_{\Tilde{l}}(1)\simeq\oh(-\check{\tau})$. Consequently
$\I_{\Tilde{p}}\I_{\Tilde{l}}(2)\simeq\I_{\Tilde{p}}(1)\otimes
\oh(-\check{\tau})$, so
\begin{equation}
\label{n4-2}
pr_{\ast}\I_{\Tilde{l}}(1)\simeq\oh(-\check{\tau})\qquad{\rm
and}\qquad
pr_{\ast}\I_{\Tilde{p}}\I_{\Tilde{l}}(2)\simeq pr_{\ast}\I_{\Tilde{p}}(1)\otimes
\oh(-\check{\tau})\,.
\end{equation}
The exact sequence \eqref{e:e6-1a}, along with \eqref{n4-2}
and \eqref{n4-2a} results in
\begin{equation}
\label{a2-5}
\begin{split}
c(\V)&=c(pr_{\ast}{\I}_{\Tilde{p}}(2))
c(pr_{\ast}\I_{\Tilde{p}}\I_{\Tilde{l}}(2))^{-1}\\
&=c(\Omega(\tau))^3c(\wedge^2\Omega(\tau))^{-1}c(\Omega(\tau)
\otimes\oh(-\check{\tau}))^{-1}\,,
\end{split}
\end{equation}
where $c(\Omega(\tau))=1-\tau+\tau^2$.

\bigskip

\subsection{Lines on the Calabi-Yau threefolds}
\label{linesCY}

Recall the construction of Section~\ref{gqs}. Let
\begin{equation*}
        \begin{CD}
                \smap_{0,1}(\NN,1) @>e_{1}>>\NN \\
                @V\pi VV \\
                \smap_{0,0}(\NN,1)
        \end{CD}
\end{equation*}
be the universal map, and let $\B$ denote any of the three
bundles $3\oh(p)$, $\F^{\vee}\oplus\oh(2p)$ and $S_{2}\F^{\vee}$.
Then $\B_{1}=\pi_{\ast}e_{1}^{\ast}\B$ is locally free with
fibers $H^0(\PP^1,f^{\ast}\B)$ over degree one maps
$f\:\PP^1\to\NN$, and the virtual numbers of lines on the Calabi-Yau
sections $X$, $Y$, $Z$ are given by the integrals
\begin{equation}
\label{n4-3a}
\int_{\moo(\NN,1)}c_{6}(3\oh(p)_{1})\,,\qquad
\int_{\moo(\NN,1)}c_{6}(\F_{1}^{\vee}\oplus\oh(2p)_{1})\,,\qquad
\int_{\moo(\NN,1)}c_{6}(S_{2}\F_{1}^{\vee})\,,
\end{equation}
respectively. With the description of the ring
$A^{\ast}(\moo(\NN,1))$ in the previous section, we are in
a position to evaluate the degree of any cycle. Hence, to
evaluate the above integrals we seek expressions for the
integrands in terms of the generators $\tau$, $\check{\tau}$,
and $\sigma$.

Consider the universal family
\begin{equation}
        \begin{CD}
        \C_1 @>e>> \PP(\V^{\vee})@> >>\NN\\
        @V\pi VV             @VV{}V\\
        H^{0,1} @>>{}>   PP^{\vee}\,.
        \end{CD}
\end{equation}
 From the remark in Section~\ref{modulspace},
\begin{equation}
\begin{split}
\E_{\C_{1}}=&pr_{\ast}\I_{\Tilde{p}}\I_{\Tilde{l}}(2)_{\C_{1}}
\oplus\oh_{\C_{1}}
(-\check{\sigma})\\
\F_{\C_{1}}\simeq &pr_{\ast}\I_{\Tilde{l}}(1)_{\C_{1}}\otimes
\bigl(\wedge^2pr_{\ast}\I_{\Tilde{p}}(1)_{\C_{1}}\oplus\oh_{\C_{1}}
(-\check{\sigma})\bigr)\,,
\end{split}
\end{equation}
hence
\begin{equation}
\label{a2-7}
\begin{split}
&\oh(p)_{1}=\pi_{\ast}\wedge^3\E_{\C_{1}}^{\vee}=
\wedge^2pr_{\ast}\I_{\Tilde{p}}\I_{\Tilde{l}}(2)^{\vee}\otimes
\pi_{\ast}\oh_{\C_{1}}(\check{\sigma})\\
&\F_{1}^{\vee}\simeq
pr_{\ast}\I_{\Tilde{l}}(1)^{\vee}\otimes
\bigl(\wedge^2pr_{\ast}\I_{\Tilde{p}}(1)^{\vee}\oplus
\pi_{\ast}\oh_{\C_{1}}(\check{\sigma})\bigr)\,.
\end{split}
\end{equation}
By a diagonal argument as in the previous section,
$\I_{\C}\simeq\oh(-\sigma)\otimes\oh(-\check{\sigma}))$, where
\[0\to\I_{\C}\to\oh_{\PP(\V)\times_{PP^{\vee}}\PP(\V^{\vee})}
\to\oh_{\C}\to 0\]
is the ideal sequence. This implies that
$\pi_{\ast}\oh_{\C}(\check{\sigma})=\coker\left(\oh(-\sigma)\to
\V_{H^{0,1}}^{\vee}\right)$. Also,
$\oh(2p)_{1}=S_{2}(\oh(p)_{1})$ and
$(S_{2}\F^{\vee})_{1}=S_{2}(\F_{1}^{\vee})$.

The desired expressions for the integrands can now be
obtained from formal Chern class computations using the
identifications \eqref{n4-2},\eqref{a2-5}, and \eqref{a2-7}.
These
computations and the evaluations of the integrals are effectively
done using the
{\em Schubert} package in  {\texttt{MAPLE}} \cite{KS}.

\begin{thm}\label{thm-int}
The number of lines on the Calabi-Yau threefolds $X$, $Y$,
and $Z$ which are given by the integrals \eqref{n4-3a} are
$147$, $216$, and $144$, respectively.
\end{thm}

There remains to show that the virtual numbers given by the
integrals \eqref{n4-3a} are in fact the actual numbers of
lines.

\begin{lem}
If $\B$ is one of the bundles $3\oh(p)$,
$\F^{\vee}\oplus\oh(2p)$ or $S_{2}\F^{\vee}$ then the map
\begin{equation}
\label{e:map2}
H^{0}(\NN,\B)_{\moo(\NN,1)}\xrightarrow{\nu}
\B_{1}
\end{equation}
is surjective.
\end{lem}

The proof of Theorem~\ref{thm-int} follows from this.
Indeed, if $s\in H^0(\NN,\B)$, then
$V(\nu(s))\subseteq\moo(\NN,1)$ describes the  locus of
lines on the zero scheme $V(s)\subseteq \NN$. From the lemma
we have that $\B_{1}$ is generated by $H^0(\NN,\B)$,
hence, by Kleiman's Bertini theorem \cite{Kl}, it follows that
for a general global section $s$, $V(\nu(s))$ is nonsingular of
codimension $\rank(\B_{1})$. Since
$\rank(\B_{1})=\dim\moo(\NN,1)$, the
number of points is finite and given by
\[\int_{\moo(\NN,1)}c_{6}(\B_{1})\,.\]
Further, the section $s$ can be chosen such that
$\nu(s)$ has zeros outside any closed subset, for
instance the subset $\PP(\V_{I})\subseteq\moo(\NN.1)$, where
$I\subseteq\PP(V)\times\PP(\V^{\vee})$ is the incidence
correspondence.
\begin{proof}
Let $f\:{\PP}^1\to\NN$ be a map of degree one. It is enough
to show surjectivity on the fibers of \eqref{e:map2}
\begin{equation}
H^{0}(\NN,\B)\to
H^{0}({\PP}^1,f^{\ast}\B)\,.
\end{equation}
Let $pl\:\NN\to{\PP}^{19}$ be the Pl{\"u}cker embedding.
For all $m$,
$pl^{\ast}\oh_{{\PP}^{19}}(m))=\oh(mp)$, and the map
\begin{equation}
\label{e:e10-1}
H^{0}({\PP}^{19},\oh_{{\PP}^{19}}(m))\to
H^{0}(\PP^1,f^{\ast}\oh(mp))
\end{equation}
factors through the map
\begin{equation}
\label{e:e10-2}
H^{0}(\NN,\oh(m))\to
H^{0}({\PP}^1,f^{\ast}\oh(mp))\,.
\end{equation}
The map \eqref{e:e10-1} is surjective for $m=1$ and $2$, hence so is
\eqref{e:e10-2}.

Next, since $\deg f^{\ast}\F^{\vee}=1$ and
$f^{\ast}\F^{\vee}$ is generated by its global sections
\begin{equation}
\label{e:e10-3}
S_{2}V^{\vee}\otimes V^{\vee}\to
H^{0}({\PP}^1,f^{\ast}\F^{\vee})
\end{equation}
is surjective, where \eqref{e:e10-3} is induced from
the inclusion $\F\to S_{2}V\otimes V_{\NN}$.
The map \eqref{e:e10-3} factors through
$H^{0}(\NN,\F^{\vee})\to
H^{0}({\PP}^1,f^{\ast}\F^{\vee})\,,$
so this map is also surjective.

Finally, $H^{0}({\PP}^1,f^{\ast}S_{2}\F^{\vee})=
S_{2}H^{0}({\PP}^1,f^{\ast}\F^{\vee})$\,,
so right exactness of $S_{2}$ implies that
\begin{equation}
\label{e:e10-4}
S_{2}(S_{2}V^{\vee}\otimes V^{\vee})\to
H^{0}({\PP}^1,f^{\ast}S_{2}\F^{\vee})
\end{equation}
is surjective. The map \eqref{e:e10-4} factors through
$H^{0}(\NN,S_{2}\F^{\vee})\to
H^{0}({\PP}^1,f^{\ast}S_{2}\F^{\vee})$\,,
and so the claim follows.
\end{proof}

\bigskip

\subsection{Fixlines}
 \label{fixlines}

The torus action $T$ on $\NN$ induces an action on $\moo(\NN,1)$.
The projection
\[\rho\:\moo(\NN,1)\to\PP(V)\times\PP(V^{\vee})\]
is
$T$-equivariant. Let $L\subseteq\NN$ be a line in $\NN$.
If $[L]\in\moo(\NN,1)$ is a fixpoint, then
$\rho([L])=(p,l)$ must correspond to a fixpoint and a fixline
in $\PP(V)$.
There are two projective equivalence classes for the
point-line pair, depending on whether the point is incident
to the line or not. Up to the $S_{3}$-action we may assume
$p=(x_{1},x_{2})$ and $l=(x_{0})$, or $l=(x_{1})$. In
addition, the line $L\subseteq\PP(\V_{(p,l)}^{\vee})$ must
be fixed. For each choice of $(p,l)$ there are three
fixlines in $\PP(\V_{(p,l)}^{\vee})$. Since the lines
$\lav x_{0}x_{1},x_{0}x_{2},x_{1}^2+tx_{0}x_{2}\rav$ and
$\lav x_{0}x_{1},x_{0}x_{2},x_{2}^2+tx_{0}x_{1}\rav$ are
swapped by the permutation $x_{1}\leftrightarrow x_{2}$,
this leaves us with five $S_{3}$-classes of fixlines.
Table~\ref{tablefixlines} contains a complete list of
the fixlines in $\moo(\NN,1)$, up to the $S_3$-action on $V$, along
with the isotropy subgroups of this action.
\begin{table}[h!]\centering
\caption{Fixlines\label{tablefixlines}}
\begin{tabular}{||c|c|c||}
\hline\hline
&&\\
Type&$\E_t$&Isotropy\\&&\\\hline&&\\
1&$\lav x_{0}x_{1},x_{0}x_{2},x_{1}^2+tx_{2}^2\rav$
 &$x_{1}\leftrightarrow x_{2}$
\\&&\\\hline&&\\
2&$\lav x_{0}x_{1},x_{0}x_{2},x_{1}x_{2}+tx_{1}^2\rav$
 &$\{id\}$
\\&&\\\hline&&\\
3&$\lav x_{1}^2,x_{1}x_{2},x_{0}x_{1}+tx_{2}^2\rav$
 &$\{id\}$
\\&&\\\hline&&\\
4&$\lav x_{1}^2,x_{1}x_{2},x_{2}^2+tx_{0}x_{2}\rav$
 &$\{id\}$
\\&&\\\hline&&\\
5&$\lav x_{1}^2,x_{1}x_{2},x_{0}x_{2}+tx_{0}x_{1}\rav$
 &$\{id\}$
\\&&
\\
\hline\hline
\end{tabular}
\end{table}

\medskip

\noindent
{\em Remark.}
 Up to projective equivalence there are six types of lines on
$\NN$: The five described above and the line
\[{\rm Type~6)}\qquad\lav
x_{1}^2,x_{1}x_{2},x_{2}^2+x_{0}x_{1}+tx_{0}x_{2}\rav\,.\]

\bigskip

\subsection{Positivity}
\label{positivity}

Let $L$ be a line in $\NN$. Since $L\simeq \PP^1$, the
vector bundle $T\NN_{L}$ decomposes into a sum of line
bundles. In Appendix~\ref{A-tangent} we show that
\begin{equation*}
T\NN_{L}=\left\{\begin{array}{lc}
\oh(2)\oplus\oh(1)\oplus 4\oh
&\mbox{if $L$ is of type
1)--2)}\\
\oh(2)\oplus2\oh(1)\oplus 4\oh\oplus\oh(-1)
&\mbox{ if $L$ is of type
3)--6)}
\end{array}\right.\,.
\end{equation*}
In particular, we see that $\NN$ is not convex, since
$H^1(f^{\ast}T\NN)\ne 0$ if $f\:C\to\NN$ is any multiple
cover of a line of type 3)--6).

\addtocontents{toc}{\vspace{0.2cm}}
\section{Quantum $\D$-module of $\NN$}
\label{qdmN}\addtocontents{toc}{\vspace{0.2cm}}

\bigskip

We now put the machinery of Section~\ref{GQH} to use. Our
goal is to compute the instanton numbers of rational curves
for the Calabi-Yau section $X$. This enables us
(in principle) to
determine all virtual numbers of rational curves on $X$.
Remarkably, the only geometrical input this
procedure requires is the Chow ring
$A^{\ast}(\NN)$, the resolution \eqref{e:e5-1-1}, and the
description of the fixlines, Table~\ref{tablefixlines}.

In what follows $T_{0},\ldots,T_{12}$ will refer to the
basis of $A_{\QQ}^{\ast}(\NN)$ where $T_{0}=1$, and the
rest is as in
Table~\ref{tablebasis} , with $\gamma_{i}=(-1)^ic_{i}(\E)$
and $\delta_{2}= c_{2}(\F)$. The sign $(-1)^i$ makes the
classes effective. Note
that $T_{11}$ is the class of a line, and $T_{12}$ is the
class of a point. As before, $p=\gamma_{1}$, $T^j$ denotes
the dual generators, and $t_{0}$ (resp. $t$) are
coordinates on $A_{\CC}^0(\NN)$ (resp. $A_{\CC}^1(\NN)$).
\begin{table}[h!]\centering
\caption{Basis\label{tablebasis}}
\begin{tabular}{||c|c|c|c|c|c||}
\hline\hline
&&&&&\\
$T_1=\gamma_{1} $&$T_2=\gamma_{1}^2 $ &$T_3=\gamma_{2}$&
$T_4=\delta_{2} $&
$T_5=\gamma_{1}^3$&$T_6=\gamma_{1}\gamma_{2}$\\&&&&&\\\hline&&&&&\\
$T_7=\gamma_{1}\delta_{2} $&$T_8=\gamma_{1}^4$ &
                $T_9=\gamma_{1}^2\gamma_{2} $&
$T_{10}=\gamma_{1}^2\delta_{2}$&$T_{11}=\frac{1}{57}\gamma_{1}^5$&
$T_{12}=\frac{1}{57}\gamma_{1}^6$\\&&&&&\\
\hline\hline
\end{tabular}
\end{table}

\bigskip

\subsection{The $A$-connection}
\label{aconnection}

We begin by determining the system of first order differential
equations
\begin{equation}
\label{e:e7-1-0}
\hbar \frac{\partial}{\partial t_{0}}F=p\ast F\,,\qquad
\hbar\frac{\partial}{\partial t}F=p\ast F
\end{equation}
for vector valued functions $F\:A_{\CC}^{0}(\NN)
\oplus A_{\CC}^{1}(\NN)\to
A_{\CC}^{\ast}(\NN)$. The operator $p\ast-$ will be represented
as a matrix with respect to the basis
$T_{0},T_{1},\ldots,T_{12}$.

Recall
\begin{equation*}
\begin{split}
p\ast T_a&=\sum_d q^d\sum_b\lav p,T_a,T_b\rav_{0,d}T^b\\
&=\sum_d q^d\sum_b\lav p,T_a,T_b\rav_{0,d}\sum_c
g^{bc}T_c\,,
\end{split}
\end{equation*}
where the matrix $(g^{bc})$ is the inverse of the matrix $(g_{bc})$
with $g_{bc}=\lav T_b, T_c\rav$.
By the mapping to a point and the divisor properties,
\begin{equation}
\label{e:e7-1-1}
p\ast T_a=p\cdot T_a+\sum_{q>0} q^d\sum_{b,c}d\lav
T_a,T_b\rav_{0,d}
g^{bc}T_c\,.
\end{equation}
The intersection product
\,$p\cdot T_a=\sum_{b,c}(\int_{\NN}pT_{a}T_{b})g^{bc}T_c$\,
is determined from the list of monomial values in
Table~\ref{tablemonomialvalues}. Hence, to determine \eqref{e:e7-1-1}
there remains to compute all two-point numbers
$\lav T_a,T_b\rav_{0,d}$ for all $d$. However, $\lav
T_a,T_b\rav_{0,d}=0$ unless
\begin{equation}
\label{e:e7-1-2}
\codim T_a+\codim T_b={\rm exp}.\dim\smap_{0,2}(\NN,d)=3d+5\,.
\end{equation}

\begin{thm}
\label{thm7-1}
There are twelve two-point line numbers satisfying \eqref{e:e7-1-2},
with values as in Table~\ref{tabletwopointlinenumbers}, and
one two-point conic number satisfying \eqref{e:e7-1-2}, with
value
\[\lav T_{11},T_{12}\rav_{0,2}=1\,.\]
\end{thm}
\begin{table}[h!]\centering
\caption{Two-point line numbers\label{tabletwopointlinenumbers}}
\begin{tabular}{||c|c|c|c||}
\hline\hline&&&\\
$\lav T_2,T_{12}\rav_{0,1}=3$&$\lav T_5,T_{11}\rav_{0,1}=8$&
$\lav T_8,T_{8}\rav_{0,1}=603$&$\lav T_9,T_9\rav_{0,1}=121$\\
&&&\\\hline&&&\\
$\lav T_{3},T_{12}\rav_{0,1}=0$&$\lav T_6,T_{11}\rav_{0,1}=3$&
$\lav T_8,T_9\rav_{0,1}=270$&$\lav T_9,T_{10}\rav_{0,1}=81$\\
&&&\\\hline&&&\\
$\lav T_4,T_{12}\rav_{0,1}=0$&$\lav T_7,T_{11}\rav_{0,1}=2$&
$\lav T_8,T_{10}\rav_{0,1}=180$&$\lav T_{10},T_{10}\rav_{0,1}=54$\\
&&&\\
\hline\hline
\end{tabular}
\end{table}

The line numbers are computed using Bott's formula on the
moduli space $\smap_{0,2}(\NN,1)$. In Section~\ref{representationsmotwo}
we explain this
computation. The number $\lav T_{11},T_{12}\rav_{0,2}=1$
comes out as a solution of the WDVV-relations (see the remark
below)
once the two-point line numbers are known. Another
way to compute this number is to use Bott's formula over
the space $\smap_{0,2}(\NN,2)$. In Section~\ref{conics} we
give a description of the space $\smap_{0,0}(\NN,2)$ which
contains enough information for such a computation to be
carried out.

Putting this data into \eqref{e:e7-1-1} we find that $p\ast-$ is represented
by the matrix
\[\left(
\begin{array} {ccccccccccccc}
0&0&3q&0&0&0&0&0&0&0&0&2q^2&0\\                   
1&0&0&0&0&8q&3q&2q&0&0&0&0&2q^2\\                 
0&1&0&0&0&0&0&0&21q&9q&6q&0&0\\                   
0&0&0&0&0&0&0&0&0&-q&0&0&0\\                      
0&0&0&0&0&0&0&0&-33q&-12q&-9q&0&0\\               
0&0&1&0&0&0&0&0&0&0&0&2/3q&0\\                    
0&0&0&1&0&0&0&0&0&0&0&0&0\\                       
0&0&0&0&1&0&0&0&0&0&0&-5/3q&0\\                   
0&0&0&0&0&1&0&0&0&0&0&0&q\\                       
0&0&0&0&0&0&1&0&0&0&0&0&0\\                       
0&0&0&0&0&0&0&1&0&0&0&0&-3q\\                     
0&0&0&0&0&0&0&0&57&27&18&0&0\\                    
0&0&0&0&0&0&0&0&0&0&0&1&0                         
   \end{array}\right)\]
with respect to the basis $\{T_{0},T_{1},\ldots,T_{12}\}$.

\medskip

\noindent{\em Remark.} There is a computer program {\texttt {farsta}}
written by Andrew Kresch \cite{Kr} which can be quite
useful for solving WDVV relations on the GW-level if the
rank of the intersection ring is not too large. Suppose we
have a scheme $S$ and a good description of $A^{\ast}(S)$.
The program will derive all WDVV relations among the
GW-invariants of $S$, and, once enough numbers are specified,
start to solve the equations for new numbers. The number
$\lav T_{11},T_{12}\rav_{0,2}=1$ was first obtained by
{\texttt{farsta}}.

\bigskip

\subsection{The Picard-Fuchs operator}
\label{pf}

Suppose
\[F=\sum_{i=0}^{12}F_{i} T_{i}=
\begin{pmatrix}
F_{0}\\ \cdot\\\cdot\\ F_{12}
\end{pmatrix}\]
\\
is a solution to the system of differential equations
\eqref{e:e7-1-0}.
Solving for $F_{12}$ as in Example~\ref{exPP5} of
Section~\ref {QDmodule}
(with $\hbar=1$), we find that $F_{12}$ is a solution to
 $Q\D_{\NN}=0$, where $Q\D_{\NN}$ is the Picard-Fuchs operator
\begin{equation*}
\begin{split}
&D^7(299D^2-1035D+907)(D-1)^3\\
&-qD^3(10166D^6+5474D^5-7135D^4-5855D^3+1148D^2+2109D+513)\\
&+q^2(-83122D^6-377246D^5-675645D^4
-607063D^3
-289727D^2-70962D-7668)\\
&-243q^3(D+1)(299D^2+1357D+1551)\,,
\end{split}
\end{equation*}
with $D=\partial/\partial t$.

\medskip

\noindent{\em Remark~1.}
 The actual computation of $Q\D_{\NN}$ was carried out in {\texttt
{MAPLE}}, using a procedure written by Ciocan-Fontanine and
van~Straten.

\medskip

\noindent{\em Remark~2.}
By the remark in Section~\ref{QDmodule}, it follows that the
homogeneous degree $12$ part of $Q\D_{\NN}$, with $D=p$, is a
relation in $QH^{\ast}(\NN)$. A closer inspection of the
factors in this expression reveals that in fact
\[p(q+p^3)(p^6-35qp^3-243q^2)\]
is a relation.

 \bigskip

\subsection{Instantons}
\label{instantons}

We use the quantum Lefschetz hyperplane principle to
compute the instantons for the
complete intersection Calabi-Yau threefold $X\subseteq\NN$.
It will be convenient to use a basis $\{p^2,k_{1},k_{2}\}$ for
$A_{\CC}^2(\NN)$ such that
$k_{1}\cdot p^3=0=k_{2}\cdot p^3$. In
other words $k_{1}$ and $k_{2}$ span
$\ker(i^{\ast}:A_{\CC}^2(\NN)\to A_{\CC}^2(X))$. Let
$p^{2\vee}\in A_{\CC}^4(\NN)$ be such that
$\lav p^2,p^{2\vee}\rav=1$ and
$\lav k_i,p^{2\vee}\rav=0$ for $i=1,2$. Then
$p^{2\vee}=(1/57)T_{8}$. Let $p'=i^{\ast}p$ denote the
ample generator of $\Pic(X)$.

First, we need to determine $I_0,I_{1},I_{2}$ such that
\begin{equation}
\begin{split}
i^{\ast}I=&i^{\ast}e^{(t_0+pt)/\smallhbar}
\left(1+\sum_{d\ne 0}q^{d}
\prod_{m=1}^{d}(p+m\hbar)^3e_{1\ast}
\left(\frac{1}{\hbar-c}\right)\right)\\
=&e^{t_0/\smallhbar}(I_0+I_{1}p'+I_{2}{p'}^2+I_{3}{p'}^3)\,.
\end{split}
\end{equation}
The crucial point in what follows, is that the expression (with $\hbar=1$)
\[s_{\NN}=e^{(t_0+pt)/\smallhbar}
\left(1+\sum_{d\ne 0}q^{d}
e_{1\ast}\left(\frac{1}{\hbar-c}\right)\right)\]
is a vector solution to $Q\D_{\NN}$. 
The first few components of $s_{\NN}$ with respect to the
new basis for $A_{\CC}^2(\NN)$ are
\begin{multline*}
s_{\NN}=e^{t_0/\smallhbar}\biggl(\psi_0+
\biggl(\frac{t}{\hbar}\psi_0+\psi_{1}\biggr)p+
\biggl(\frac{t^2}{2\hbar^2}\psi_0+\frac{t}{\hbar}\psi_{1}+
\psi_{2}\biggr)p^2+
\lav s_{\NN},k_{1}^{\vee}\rav k_{1}\\
+\lav s_{\NN},k_{2}^{\vee}\rav k_{2}
+\cdots\biggr)\,,
\end{multline*}
where
\begin{equation*}
\psi_{0}=1+\sum_{d>0}q^{d}\frac{a_{d}}{\hbar^{3d}}
\,,\quad
\psi_{1}=\sum_{d>0}q^{d}\frac{b_{d}}{\hbar^{3d+1}}
\,,\quad{\rm and}\quad
\psi_{2}=\sum_{d>0}q^{d}\frac{c_{d}}{\hbar^{3d+2}}\,,
\end{equation*}
with
\begin{equation*}
a_{d}\!=\!\!\int_{\![\smap_{0,2}(\NN,d)]}\!\!\!\!\!\!c^{3d-1}
e_{1}^{\ast}(T_{12})\,,\;
b_{d}\!=\!\!\int_{\![\smap_{0,2}(\NN,d)]}\!\!\!\!\!\!c^{3d}
e_{1}^{\ast}(T_{11})\,,\;
c_{d}\!=\!\!\int_{\![\smap_{0,2}(\NN,d)]}\!\!\!\!\!\!c^{3d+1}
e_{1}^{\ast}
  (p^{2\vee})\,.
\end{equation*}
Write
$Q\D_{\NN}=P_{0}(D)+qP_{1}(D)+q^2P_{2}D)+q^3P_{3}(D)$. Since
$P_{0}(n)\ne 0$ for $n>1$, it follows from standard theory
about Picard-Fuchs equations (see
Appendix~\ref{PicardFuchs}) that the coefficients of the
hypergeometric series are determined recursively from the
relations
\begin{equation}
\label{e:e9-2}
\begin{split}
&\sum_{i=0}^3 a_{n-i}P_{i}(n-i)=0\\
&\sum_{i=0}^3 a_{n-i}P_{i}^{\prime}(n-i)+
\sum_{i=0}^3 b_{n-i}P_{i}(n-i)=0\\
&\frac{1}{2}\sum_{i=1}^3 a_{n-i}P_i^{\prime\prime}(n-i)+
\sum_{i=1}^3 b_{n-i}P_i^{\prime}(n-i)+
\sum_{i=1}^3 c_{n-i}P_i(n-i)=0
\end{split}
\end{equation}
for all $n>1$, and
\begin{equation*}
a_{1}\!=\!\!\int_{\!\smap_{0,2}(\NN,1)}\!\!\!\!\!\!c^2
e_{1}^{\ast}(T_{12})\,,\qquad
b_{1}\!=\!\!\int_{\!\smap_{0,2}(\NN,1)}\!\!\!\!\!\!c^3
e_{1}^{\ast}(T_{11})\,,\qquad
c_{1}\!=\!\!\int_{\!\smap_{0,2}(\NN,1)}\!\!\!\!\!\!c^4
e_{1}^{\ast}(p^{2\vee})\,.
\end{equation*}
These integrals can be evaluated using Bott's formula. We
explain the computation in Section~\ref{representationsmotwo}.

\begin{lem}
\label{lem9-1}
$a_{1}=3$, $b_{1}=-1$, and $c_{1}=-\frac{65}{19}$.
\end{lem}

Let $R_i(d)$ be the coefficient of $p^i$ in
$\prod_{m=1}^{d}(p+m)^3$. Then
\begin{equation}
\label{e:e9-3}
\begin{split}
&I_{0}=1+\sum_{d>0}q^{d}R_{0}(d)a_{d}\\
&I_{1}=\frac{t}{\hbar}I_{0}+\sum_{d>0}q^{d}R_{1}(d)
    \frac{a_{d}}{\hbar}+
   \sum_{d>0}q^{d}R_{0}(d)
    \frac{b_{d}}{\hbar}\\
&I_{2}=\frac{t^2}{2{\hbar}^2}I_{0}+\frac{t}{\hbar}I_{1}+
    \sum_{d>0}q^{d}R_{2}(d)
    \frac{a_{d}}{{\hbar}^2}+
   \sum_{d>0}q^{d}R_{1}(d)
    \frac{b_{d}}{{\hbar}^2}+\sum_{d>0}q^{d}R_{0}(d)
    \frac{c_{d}}{{\hbar}^2}\,.
\end{split}
\end{equation}
Now, since $A^{\ast}(X)$ is generated by $p'$,
$Q\D(X)=Q\D^{\perp}(X)$, and since $X$ is a Calabi-Yau threefold
as in Example~\ref{exCY},
\[
s_{X}^{^{\perp}}=
 e^{(t_{0}^{\prime}+p't')/{\smallhbar}}\left(
 1+\sum_{d>0}{q'}^{d}(\perp\circ
 e_{1})_{\ast}\left(\frac{1}{\hbar-c}\right) \right)
\]
is a vector solution of the differential equation
\[\bigl(\hbar\frac{d}{dt'}\bigr)^2\frac{1}{K(q')}
(\hbar\frac{d}{dt'}\bigr)^2=0\,,\]
where
\[K(q')=1+\frac{1}{57}\sum_{d>0}d^3n_{d}^{X}\frac{
q^{\prime d}}{1-q^{\prime d}}\,.\]
Therefore,
\begin{equation}
\label{e:e9-4}
s_{X}^{^{\perp}}=
e^{t_{0}^{\prime}/{\smallhbar}}\left(1+\frac{t'}{\hbar}p'+
 \left(\frac{{t'}^2}{2{\hbar}^2}+
 \frac{1}{57}\sum_{d>0}{q'}^{d}\frac{n_{d}^{X}d}{\hbar^2}
  \sum_{k}\frac{{q'}^{kd}}{k^2}\right){p'}^2+\cdots\right)\,.
\end{equation}

The quantum Lefschetz hyperplane principle identifies
(with $\hbar=1$)
\[\frac{I_{1}}{I_{0}}=t'\qquad {\rm and}\qquad
\frac{I_{2}}{I_{0}}=\frac{{t'}^2}{2}+\frac{1}{57}
\sum_{d>0}{q'}^{d}d n_{d}^{X}
  \sum_{k}\frac{{q'}^{kd}}{k^2}\,.\]
Solving this for $n_{d}^{X}$ gives the desired numbers.
Table~\ref{tablenumbers} contains a list of the first ten
numbers.
\begin{table}[ht]\centering
\caption{Numbers\label{tablenumbers}}
\begin{tabular}{||l r||}
\hline\hline
$n_1$&$147 $\\$n_2$&$756 $\\$n_3$&$5\,\,283$\\$n_4$&$56\,\,970 $\\
$n_5$&$738\,\,477$\\$n_6$&$10\,\,964\,\,412$\\$n_7$&$177\,\,916\,\,032 $\\
$n_8$&$3\,\,091\,\,158\,\,090$ \\$n_9$&$56\,\,551\,\,583\,\,952 $\\
$n_{10}$&$1\,\,077\,\,954\,\,415\,\,692$\\
\hline\hline
\end{tabular}
\end{table}

\medskip

\noindent{\em Remark~1.}
 Note that, amazingly, these numbers were
determined just from knowledge of the lines on $\NN$ through:
\begin{enumerate}
\item[i)] associativity relations to find $\lav
T_{11},T_{12}\rav_{0,2}$
\item[ii)] the recursion formulas \eqref{e:e9-2}
\item[iii)] the formulas \eqref{e:e9-3} and \eqref{e:e9-4}.
\end{enumerate}

\noindent{\em Remark~2.} The Picard-Fuchs operator
governing the $I_{i}$'s (with $\hbar=1$) is
\begin{equation*}
\begin{split}
Q\D_X=D^4&+q\left(-\frac{700}{19}D^4-\frac{1238}{19}D^3-\frac{999}{19}D^2
   -20D-3\right)\\
&+q^2\!\left(-\frac{64745}{361}D^4-\frac{368006}{361}D^3-\frac{609133}{361}D^2
   -\frac{21724}{19}D-\frac{5382}{19}\right)\\
&+q^3\!\left(\frac{172719}{361}D^4+\frac{17334}{19}D^3-\frac{321921}{361}D^2
   -\frac{38880}{19}D-\frac{16038}{19}\right)\\
&+q^4\!\left(\!\frac{46656}{361}D^4\!+\!\frac{841266}{361}D^3\!+\!
\frac{1767825}{361}D^2
   \!+\!\frac{1347192}{361}D\!+\!\frac{354294}{361}\right)\\
&-\frac{177147}{361}q^5(D+1)^4\,.
   \end{split}
\end{equation*}

\addtocontents{toc}{\vspace{0.2cm}}
\section{Conics on $\NN$}
\label{conics}
\bigskip

In this section we compute the virtual number of conics on the
Calabi-Yau sections $X$, $Y$, and $Z$ using classical methods,
i.e., intersection theory on
a parameter space of conics. The number obtained for $X$
coincides with $n_{2}^{X}$ in Section~\ref{instantons}. Moreover,
we show that the conics on $X$ are rigid, so the number
$n_{2}^{X}$ is the actual number of conics on $X$.

\bigskip

\subsection{The Hilbert scheme of conics on $N$}
\label{Hilbert}

Let $C\subseteq\NN$ be a smooth conic on $\NN$ and let $f\:\PP^1\to C$
be a parametrization. There are two possibilites for the
splitting type of $f^{\ast}\E_{C}$:
\[f^{\ast}\E_{C}\simeq 2\oh_{\PP^1}\oplus \oh_{\PP^1}(-2)
\quad{\rm and}\quad
f^{\ast}\E_{C}\simeq\oh_{\PP^1}\oplus
2\oh_{\PP^1}(-1)\,.\]
We say a smooth conic is of type $(0,2)$ or $(1,1)$ according
to how $f^{\ast}\E_{C}$ decomposes.

Let $H=H_{2}$ denote the components of the Hilbert scheme of $\NN$
containing the set of conics, and let $H^{0,2}$ and $H^{1,1}$
denote the
closure in $H$ of the locus of smooth conics of type
$(0,2)$ and $(1,1)$ respectively.

\begin{thm}
\label{thmIII-1}
\quad
\begin{enumerate}
\item[i)] $H^{0,2}$ and $H^{1,1}$ are schemes
of dimension $9$.
\item[ii)] $H^{0,2}$ is smooth and irreducible.
\item[iii)] $H=H^{0,2}\cup H^{1,1}$ and the
scheme-theoretical intersection $H^{0,2}\cap H^{1,1}$ is
smooth of dimension $8$. In particular, $H$ is of dimension
$9$, as expected.
\item[iv)]  $H^{1,1}-H^{0,2}\cap H^{1,1}$ is smooth.
\end{enumerate}
\end{thm}

The proof follows from a series of lemmas.

\begin{lem}
\label{lemIII-1}
$H^{0,2}\simeq\PP(S_{2}\V)\,.$
\end{lem}

This was proved in Section~\ref{modulspace}.
\begin{lem}
\label{lemIII-2}
$H=H^{0,2}\cup H^{1,1}\,.$
\end{lem}
\begin{proof}
Let $C$ be a singular conic in $\NN$. We must show that
$[C]\in H^{0,2}\cup H^{1,1}$. A singular conic is either the
union of two lines intersecting at a point, or a line
doubled in a plane. Recall the morphism \eqref{e:e6-2},
$\PP(\V^{\vee})\to\NN$. If $C$ lies in the image of
$\PP(\V_{(p,l)}^{\vee})$ for some
$(p,l)$ in $PP^{\vee}$, the result follows by
Lemma~\ref{lemIII-1}, hence we only need to consider
singular conics not contained in the image of any
$\PP(\V_{(p,l)}^{\vee})$. For each such conic we shall
display a $1$-parameter family of smooth conics of
$(1,1)$-type which specialize to it.

First, we consider the reduced singular conics.
Suppose $L_{1}$ and $L_{2}$ are two lines in $\NN$ that meet
in a point, and $\rho([L_{i}])=(p_{i},l_{i})$.
Then $L_{i}$ lies in the image $P_{i}$ of
$\PP(\V_{(p_{i},l_{i})}^{\vee})$. If $(p_{1},l_{1})=(p_{2},l_{2})$,
then $[L_{1}\cup L_{2}]\in H^{0,2}$. Therefore we assume
$(p_{1},l_{1})\ne(p_{2},l_{2})$. Since the lines
meet in a point, the planes $P_{1}$ and $P_{2}$ must also meet.
In the proof of Theorem~\ref{thmfunctor} we saw that when
this occurs, then $P_{1}\cap P_{2}$ is a point of
projective equivalence type 1)--4) in Table~\ref{tablefixpoints},
and
that the union $P_{1}\cup P_{2}$ is determined up to
projective equivalence by this point. Hence for each point
$P_{1}\cap P_{2}$ of type 1)--4), there is a $\PP^1\times\PP^1$-choice
for the line pair $(L_{1},L_{2})$ such that
$P_{1}\cap P_{2}\in L_{i}\subseteq P_{i}$.
Table~\ref{tablelinepairsM-M} contains a list of all such line pairs
up to projective equivalence. There
$[a,b],[c,d]\in\PP^1$, and the line pairs intersect at the point
$t_{1}=0=t_{2}$.

Consider the family of matrices $M_{k}(t)$ in
Table~\ref{tablelinepairsM-M}. For
general $k\in\CC^{\ast}$, $M_{k}(t)$ is stable, and
 $C_{k}=\{(M_{k}(t))\in\NN\mid t\in\PP^1\}$ is
a $(1,1)$-curve. It is straightforward to obtain the
equations for (the Pl\"ucker embedding of) $C_{k}$ in
$\PP(\wedge^3S_{2}V^{\vee})$, and from these one may check,
for each case in Table~\ref{tablelinepairsM-M}, that
$C_{k}$ specializes to $L_{1}\cup L_{2}$ when $k\mapsto 0$.
This shows that all reduced singular conics are
contained in $H^{0,2}\cup H^{1,1}$.

We now look at the double lines in $\NN$.
Let $L$ be a line in $\NN$. The functor of genus $0$ double
lines in $\NN$ with support on $L$ is represented by
$\PP(H^{0}(N_{L/\NN}(-1))^{\vee})$. If $L$ is a line of
type 1) or 2) in Table~\ref{tablefixlines}, then $\dim
H^{0}(N_{L/\NN}(-1))=1$ so there
is only one genus $0$ double line with support on $L$. This
double line must be the one spanning the image of
$\PP(\V_{\rho([L])}^{\vee})$,
hence it is represented by a point in $H^{0,2}$. If $L$ is of
type 3)--6), then  $\dim H^{0}(N_{L/\NN}(-1))=2$, so
there exists a one-parameter family of genus
$0$ double lines with support on $L$. For each type of line, we
exhibit a family of $(1,1)$-curves $C_{m,k}$, with
$(m,k)\in\PP^1\times
\CC^{\ast}$, such that $C_{m,k}$ specializes to a
double line represented by $m$, as $k\mapsto 0$. The
parametrization $m$ of $\PP^1$
is chosen such that $C_{0,k}$ specializes to the double
line spanning $\PP(\V_{\rho([L])}^{\vee})$.

For general
$k\in\CC^{\ast}$, the family of matrices $M_{m,k}(t)$ in
Table~\ref{tabledoublelines}
are stable, and $C_{m,k}=\{(M_{m,k}(t))\in
\NN\mid t\in\PP^1\}$ are $(1,1)$-curves. Again using the
Pl\"ucker embedding, it is straightforward to verify that
$C_{m,k}$ specialize to lines doubled in the
planes $\PP_{m}^2\subseteq\PP(\wedge^3S_{2}V^{\vee})$, generated by the vectors
listed in Table~\ref{tabledoublelines}. This concludes the proof of
Lemma~\ref{lemIII-2}.
\end{proof}
\begin{table}[h!]\centering
\caption{Line pairs in $M^{1,1}-M^{0,2}\cap
M^{1,1}$\label{tablelinepairsM-M}}
\begin{tabular}{||r|c|c|c||}
\hline\hline
&$M_{k}(t)$&Line pairs&Isotropy\\
\hline&&&\\
$\begin{array}{r} l_{1}=l_{2}\\\end{array}$&
$\begin{pmatrix}
                    x_{1}&cx_{0}+dx_{2}\\
                    k(ax_{2}+bx_{0})&x_{1}\\
                    x_{2}&tx_{0}\end{pmatrix}$&
$\begin{array}{c}
\lav x_{1}^2,x_{1}x_{0},x_{1}x_{2}+t_{1}(ax_{0}x_{2}+bx_{0}^2)\rav\\
\lav x_{1}^2,x_{1}x_{2},x_{1}x_{0}+t_{2}(cx_{0}x_{2}+dx_{2}^2)\rav
\end{array}$
&$x_{0}\leftrightarrow x_{2}$
\\&&&\\\hline&&&\\
$\begin{array}{r}l_{1}\ne l_{2}\\{\rm  i)}\end{array}$&
$\begin{pmatrix}
                     x_{1}&cx_{0}\\
                     kax_{0}&x_{2}\\
                     x_{0}+kbtx_{2}&dx_{1}+tx_{0}\end{pmatrix}$&
$\begin{array}{c}
\lav x_{1}x_{2},x_{0}x_{1},x_{0}x_{2}+t_{1}(ax_{0}^2+bx_{2}^2)\rav\\
\lav x_{1}x_{2},x_{0}x_{2},x_{0}x_{1}+t_{2}(cx_{0}^2+dx_{1}^2)\rav
\end{array}$
&$x_{1}\leftrightarrow x_{2}$
\\&&&\\\hline&&&\\
$\begin{array}{r}l_{1}\ne l_{2}\\{\rm  ii)}\end{array}$&
$\begin{pmatrix}
                     x_{1}&cx_{0}\\
                     kax_{0}&x_{2}\\
                     x_{2}+kbtx_{0}&dx_{1}+tx_{0}\end{pmatrix}$&
$\begin{array}{c}
\lav x_{1}x_{2},x_{0}x_{1},x_{2}^2+t_{1}(ax_{0}^2+bx_{0}x_{2})\rav\\
\lav x_{1}x_{2},x_{2}^2,x_{0}x_{1}+t_{2}(cx_{0}x_{2}+dx_{1}^2)\rav
\end{array}$
&$\{id\}$\\
&&&\\\hline&&&\\
$\begin{array}{r}l_{1}\ne l_{2}\\{\rm  iii)}\end{array}$&
$\begin{pmatrix}
                     x_{1}&cx_{0}\\
                     kax_{0}&x_{2}\\
                     x_{2}+kbtx_{0}&dx_{0}+tx_{1}\end{pmatrix}$&
$\begin{array}{c}\lav
x_{1}x_{2},x_{1}^2,x_{2}^2+t_{1}(ax_{0}x_{1}+bx_{0}x_{2})\rav\\
\lav x_{1}x_{2},x_{2}^2,x_{1}^2+t_{2}(cx_{0}x_{2}+dx_{0}x_{1})\rav
\end{array}$
&$x_{1}\leftrightarrow x_{2}$
\\&&&
\\
\hline\hline
\end{tabular}
\end{table}
\vspace{0.5cm}
\begin{table}[h!]\centering
\caption{Families of double lines\label{tabledoublelines}}
\begin{tabular}{||r|c|c||}
\hline\hline
&$M_{m,k}(t)$&$\PP_{m}^2$\\
\hline
&&\\
${\rm \;3)\;}$
&$\begin{pmatrix}
x_{1}&k(k+m)x_{2}\\
x_{0}&x_{1}\\
x_{2}+kx_{0}&tx_{2}\end{pmatrix}$&
$\!\!\!\begin{array}{c}
x_{1}^2\wedge x_{1}x_{2}\wedge x_{2}^2\\
 x_{1}^2\wedge x_{1}x_{2}\wedge x_{0}x_{2}\\
 x_{1}^2\wedge (x_{1}x_{2}\wedge x_{0}x_{1}+
               m x_{2}^2\wedge x_{0} x_{2})
              \end{array}$\\
              &&\\\hline&&\\
${\rm \;4)\;}$
&$\begin{pmatrix}
x_{1}&0\\
x_{0}&x_{1}+kx_{2}\\
x_{2}&k(k+m)x_{0}+tx_{2}\end{pmatrix}$&
$\!\!\!\begin{array}{c}
x_{1}^2\wedge x_{1}x_{2}\wedge x_{0}x_{2}\\
 x_{1}^2\wedge x_{1}x_{2}\wedge x_{0}x_{1}\\
  x_{1}^2\!\wedge\!( x_{1}x_{2}\!\wedge\!x_{2}^2
           \!+\! m(
           x_{0} x_{1}\!\wedge\!x_{0}x_{2}\!+\!x_{1}x_{2}\!\wedge\!
             x_{0}^2))\!\!
              \end{array}\!\!$\\
              &&\\\hline&&\\
${\rm \;5)\;}$
&$\begin{pmatrix}
x_{1}&0\\
x_{2}&x_{1}+kmx_{0}\\
x_{2}&k(k+m)x_{0}+tx_{2}\end{pmatrix}$&
$\!\!\!\begin{array}{c}
x_{1}^2\wedge x_{1}x_{2}\wedge x_{2}^2\\
 x_{1}^2\wedge x_{1}x_{2}\wedge
                               x_{0} x_{1}\\
 x_{1}^2\wedge (x_{1}x_{2}\wedge x_{0} x_{2}
             +m x_{0} x_{1}\wedge x_{2}^2)
              \end{array}$\\
              &&\\\hline&&\\
${\rm \; 6)\;}$
&$\!\!\begin{pmatrix}
x_{1}&k(k+m)x_{2}\\
x_{0}&x_{1}\\
x_{2}\!+\!kx_{0}&k(k\!+\!m)x_{0}\!+\!tx_{2}\end{pmatrix}\!\!$&
$\!\!\!\begin{array}{c}
x_{1}^2\wedge x_{1}x_{2}\wedge
(x_{2}^2+x_{0}x_{1})\\
x_{1}^2\wedge x_{1}x_{2}\wedge x_{0} x_{2}\\
x_{1}^2\!\wedge\!(x_{1}x_{2}\!\wedge\!x_{2}^2
           \!+\! m( x_{1}x_{2}\!\wedge\!x_{0}^2+ x_{2}^2\!\wedge\!
        x_{0}x_{2}))
              \end{array}$ \\
              &&\\\hline\hline
\end{tabular}
\end{table}

\bigskip

Let $I\subseteq\PP(V)\times\PP(V^{\vee})$ be the incidence
correspondence $I=\{(p,l)\mid p\in l\}$ and let
$\PP(S_{2}\V_{I})\subseteq\PP(S_{2}\V)$ be the restriction
of the $\PP^{5}$-bundle to $I$. We identify 
$\PP(S_{2}\V_{I})$ with its image under the embedding
$\PP(S_{2}\V)\to H$.

\bigskip

\begin{lem}
\label{lemIII-3}
\quad
\begin{enumerate}
\item[i)] $\dim T H_{[C]}=9$ if $[C]\in
H-\PP(S_{2}\V_{I})\,.$
\item[ii)] $\dim T H_{[C]}=10$ if $[C]\in \PP(S_{2}\V_{I})\,.$
\end{enumerate}
\end{lem}

Note that this concludes the proof of Theorem~\ref{thmIII-1}.

\begin{proof}
By standard deformation arguments,
$T H_{[C]}=H^{0}(N_{C/\NN})$. The normal bundle
$N_{C/\NN}$ sits inside an exact sequence
\[0\to TC\to T\NN_{C}\to N_{C/\NN}\to l_{C}\to 0\]
on $C$, where $l_{C}$ is a sheaf with support on the singularities
of $C$, and
 $H^1(TC)=0=H^1(l_{C})$. This implies
$H^1(T\NN_{C})=H^1(N_{C/\NN})$. Since $H$ parametrizes flat
families, $\chi(N_{C/\NN})=\dim H^{0}(N_{C/\NN})-
\dim H^1(N_{C/\NN})$ is constant for $[C]\in H$. By Riemann-Roch for
$C\simeq\PP^1$, $\chi(N_{C/\NN})=9$. Hence, to prove the
lemma, it is enough to prove
\begin{equation}
\label{e13-3}
\dim H^1(T\NN_{C})=
\left\{
     \begin{array}{lcl}
     0&&{\rm if}\quad[C]\in H-\PP(S_{2}\V_{I})\\
     1&&{\rm if}\quad[C]\in\PP(S_{2}\V_{I})\,.
     \end{array}
\right.
\end{equation}
First assume $[C]\in H^{0,2}$. Then
$C\subseteq\PP^2\subseteq\NN$, where
$\PP^2=\PP(V_{(p,l)}^{\vee})$ for some $(p,l)$ in
$PP^{\vee}$. Since
$\E_{\PP^2}\simeq 2\oh_{\PP^2}\oplus\oh_{\PP^2}(-1)$ and
$\F_{\PP^2}\simeq\oh_{\PP^2}\oplus\oh_{\PP^2}(-1)$, the
resolution \eqref{e:e5-1-1} implies that the first and
second cohomology groups of $T\NN_{\PP^2}$ and
$T\NN_{\PP^2}(-1)$ vanish. Consequently, if
$L\subseteq\PP^2$ is any line, the exact sequences
\[0\to T\NN_{\PP^2}(-2))\to T\NN_{\PP^2}\to T\NN_{C}\to
0\]
and
\[0\to T\NN_{\PP^2}(-2)\to T\NN_{\PP^2}(-1)\to T\NN_{L}(-1))\to
0\]
yields
$H^1(T\NN_{C})=H^2(T\NN_{\PP^2}(-2))=H^1(T\NN_{L}(-1))$. In
Section~\ref{positivity} we saw that
\begin{equation*}
\dim H^1(T\NN_{L}(-1))=
\left\{
     \begin{array}{lcl}
     0&&{\rm if}\quad p\not\in l\\
     1&&{\rm if}\quad p\in l\,.
     \end{array}
\right.
\end{equation*}
This proves the lemma for $[C]\in H^{0,2}$.

Now, let $[C]\in H^{1,1}-H^{0,2}\cap H^{1,1}$. By
semi-continuity, it is enough to show \eqref{e13-3}
for degenerated conics. If $C$ is the union of two
lines $L_{1}$ and $L_{2}$, then
$\rho([L_{1}])\ne\rho([L_{2}])$.
Consider the component sequence
\begin{equation}
0\to T\NN_{L_{1}\cup L_{2}}\to T\NN_{L_{1}}\oplus T\NN_{L_{2}}
\to T\NN_{L_{1}\cap L_{2}}\to 0\,.
\end{equation}
Let $\rho([L_{i}])=(p_{i},l_{i})$. If $p_{i}\not\in l_{i}$, then
$T\NN_{L_{i}}$ is generated by global sections, hence
\begin{equation}
\label{leif}
H^{0}(T\NN_{L_{1}})\oplus H^{0}(T\NN_{L_{2}})
\to T\NN_{L_{1}\cap L_{2}}
\end{equation}
is surjective, whenever this is the case for one of the lines.
If $p_{i}\in l_{i}$ for both lines, then \eqref{leif} is still
surjective, since the $\oh(-1)$-component of $T\NN_{L_{i}}$ is
canonically isomorphic to
\[H^0(\I_{l_{i}}(1))^{\vee}\otimes\wedge^2H^0(\I_{p_{i}}(1))^{\vee}
\otimes V/H^0(\I_{p_{i}}(1))\otimes\oh_{L_{i}}(-\check{\sigma})\]
(see Appendix~\ref{A-tangent}), and we are assuming $(p_{1},l_{1})\ne
(p_{2},l_{2})$.
This implies $H^1(T\NN_{L_{1}\cup L_{2}})=0$.

If $C$ is a double line, then it is projectively equivalent
to one of the double lines $D_{m}$ which spans $\PP_{m}^{2}$,
$m\ne 0$, in
Table~\ref{tabledoublelines}. Since lines of type 3)--6) in
Table~\ref{tablefixlines} specialize
 to lines of type 5), and a double line $D_{m}$, $m\ne 0$,
 specializes to
 $D_{\infty}$, it is enough to show \eqref{e13-3} for the
 double line $D_{\infty}$ which spans the plane
 \begin{equation}
 \label{e13-3a}
 \PP_{\infty}^2=\lav
x_{0}^2\wedge x_{0}x_{1}\wedge x_{0}x_{2},
x_{0}^2\wedge x_{0}x_{1}\wedge x_{1}^2,
x_{0}^2\wedge x_{0}x_{2}\wedge x_{1}^2\rav/_{\CC^{\ast}}
\end{equation}
in $\PP(\wedge^3S_{2}V^{\vee})$. All points in
$\PP_{\infty}^2$ are decomposable, hence $\PP_{\infty}^2$ is
contained in the image of $\Grass_{3}(S_{2}V)$.
Consider the restriction of the tautological exact sequence
to $\PP_{\infty}^2$,
\begin{equation}
\label{e13-3/4}
0\to \U_{\PP_{\infty}^2}\to S_{2}V_{\PP_{\infty}^2}\to
\Q_{\PP_{\infty}^2}\to 0\,.
\end{equation}
 From \eqref{e13-3a} one sees that
$0\to \U_{\PP_{\infty}^2}\to S_{2}V_{\PP_{\infty}^2}$
factors through the trivial rank four bundle
$\lav x_{0}^2, x_{0}x_{1}, x_{0}x_{2},x_{1}^2\rav_{\PP_{\infty}^2}$.
Since $\deg \Q_{\PP_{\infty}^2}=1$, this implies that
$\Q_{\PP_{\infty}^2}\simeq 2\oh_{\PP_{\infty}^2}\oplus
\oh_{\PP_{\infty}^2}(1)\,.$
This means that the restriction of \eqref{e13-3/4} to
$D_{\infty}$ is exact on global sections, so
$H^1(\E_{D_{\infty}})=0$. To show
$H^1(T\NN_{D_{\infty}})=0$ it suffices, by \eqref{e:e5-1-1},
to show that
$H^1(\F^{\vee}\otimes\E\otimes V_{D_{\infty}})=0$. But this
is now clear, since
\begin{equation}
\label{e13-4}
(S_{2}V\otimes V)^{\vee}\otimes\E\otimes V_{D_{\infty}}
\doublerightarrow\F^{\vee}\otimes\E\otimes  V_{D_{\infty}}
\end{equation}
is surjective, and
$H^1((S_{2}V\otimes V)^{\vee}\otimes\E_{D_{\infty}}\otimes
V)=0$. Here \eqref{e13-4} is the map induced from
the inclusion $\F_{D_{\infty}}\to\E_{D_{\infty}}\otimes
V$. This completes the proof of
Lemma~\ref{lemIII-3}.
\end{proof}

\medskip

\noindent{\em Remark~1.} We come up just short of proving that
$H^{1,1}$
is also smooth. However, the families
 in Table~\ref{tabledoublelines}
show that if $[C]\in H^{1,1}\cap H^{0,2}$, then the tangent
direction $TH_{[C]}-TH_{[C]}^{0,2}$ is unobstructed. Hence,
for every point $[C]\in H^{1,1}$, there exists an immersion
of $\CC^9$ onto a neighborhood of $[C]$. We expect that
$H^{1,1}$ is smooth, and that the components $H^{0,2}$ and
$H^{1,1}$ intersect transversally.

\medskip

\noindent{\em Remark~2.} We have described the Hilbert scheme of conics
on $\NN$ instead of
 $\moo(\NN,2)$. A similar result holds
for $\moo(\NN,2)$. In fact, if we let $M^{0,2}$ and
$M^{1,1}$ be the compactifications in $\moo(\NN,2)$  of the loci of
$(0,2)$-curves and $(1,1)$-curves respectively, then the
theorem holds for $\moo(\NN,2)=M^{0,2}\cup M^{1,1}$. In
particular we make the following observations:
\begin{enumerate}
\item[i)] There is a fiber square
\begin{equation*}
        \begin{CD}
        M^{0,2}\cap M^{1,1} @>{}>>M^{0,2}\\
        @V{} VV             @VV{\rho}V\\
        I @>>{}>   \PP(V)\times\PP(V^{\vee})
        \end{CD}
\end{equation*}
where the vertical maps are $\moo(\PP^2,2)$-fibrations.
\item[ii)] $M^{1,1}- M^{0,2}\cap M^{1,1}$ contains no
double coverings of lines.
\end{enumerate}

\medskip

\noindent{\em Remark~3.} It will be convenient for us to work with $H$ in
the next section. This will enable us to show that all
conics in $X$ are rigid. When we carry out our computations,
however, we will use Bott's formula over $M^{1,1}$, instead
of $H^{1,1}$. This is slightly easier because the fixpoints
on  $M^{1,1}$ are isolated, whereas on $H^{1,1}$, the
family $D_{m}$ of double lines of type 5) in
Table~\ref{tabledoublelines} results in a
$\PP^1$-component of fixpoints on $H^{1,1}$.

\medskip

\noindent{\em Remark~4.} Since $H^{0,2}$ is a projective bundle we can use
Grothendieck's theorem again to determine the Chow ring. If
we let $\xi$ be the first Chern class of the tautological
line bundle on $\PP(S_{2}\V)\xrightarrow{\rho}PP^{\vee}$,
then
\[A^{\ast}(H^{0,2})=\ZZ[\tau,\check{\tau},\xi]
/(\tau^3,\check{\tau}^{3},r)\,,\]
where
\\
\[r=\sum_{i=0}^6(-1)^i\rho^{\ast}c_{i}(S_{2}\V)\xi^{6-i}\,.\]
\\
The total Chern class $c(S_{2}\V)$ can be expressed in terms
of $\tau$ and $\check{\tau}$ by using a standard formula \cite{F} which
relates $c(S_{2}\V)$ to $c(\V)$.

\bigskip

\subsection{Fixpoints on $M^{1,1}$}
\label{fixconics}

We present a complete description of fixpoints of the torus
action on $M^{1,1}$. If $[f\:C\to\NN]\in M^{1,1}$, then
either $f$ is an embedding, or $f$ is a $2:1$ covering of a
line. If $f$ is an embedding, we identify $C$ with its
image $f(C)$.

We begin by describing the fixpoints in
$M^{1,1}-M^{0,2}\cap M^{1,1}$. As noted in
Remark~2 of Section~\ref{Hilbert}, all points in
$M^{1,1}-M^{0,2}\cap M^{1,1}$ represent embeddings
$f\:C\to\NN$.

First the smooth conics of $(1,1)$-type. Let
$\{\E_{t}\mid t\in\PP^1\}$ be a family of nets represented by
a conic of $(1,1)$-type. If the conic is fixed, then
$\{B(\E_{t})\mid t\in\PP^1\}$ must describe a fixed curve
in $\PP(V)$, hence a union of lines (possibly with multiplicity).
Since the conic is of $(1,1)$-type there exists a conic
$Q\subseteq \PP(V)$ (the generator of $\oh_{\PP^1}$ in
$f^{\ast}\E_{C}\simeq\oh_{\PP^1}\oplus 2\oh_{\PP^1}(-1)$)
that contains $\{B(\E_{t})\mid
t\in\PP^1\}$, hence $Q$ must be either i) a union of two lines,
say $x_{0}x_{1}$, or ii) a double line, say $x_{0}^2$. In
each case there are matrix representatives of type
\begin{equation}
\label{e13-5}
{\rm i)}\quad
\begin{pmatrix}
x_{0}&0\\
ax_{2}&x_{1}\\
M_{t}&x_{0}+tx_{2}
\end{pmatrix}
\qquad\qquad{\rm ii)}\quad
\begin{pmatrix}
x_{0}&0\\
L&x_{0}\\
M_{t}&x_{1}+tx_{2}
\end{pmatrix}
\end{equation}
for the nets $\E_{t}$, where
$L=ax_{1}+bx_{2}$, and
$M_{t}=(cx_{1}+dx_{2})+t(ex_{1}+gx_{2})$, with
$a,b,c,d,e,f,g$ in $\CC$.

Each smooth fixconic must contain two fixpoints. Without loss
of generality, we can assume that these occur at $t=0$ and
$t=\infty$. It is easy to verify that the matrices in
\eqref{e13-5} yield fixed conics of $(1,1)$-type, with fixpoints
at $t=0$ and $t=\infty$, only when
they take the form as in Table~\ref{tablesmoothfixconics}.
Besides the matrix representative, Table~\ref{tablesmoothfixconics}
contains a list of the nets
generated by the $2\times 2$ minors of the matrices,
along with the isotropy subgroups of the smooth fixconics.
\begin{table}[h!]\centering
\caption{Smooth fixconics\label{tablesmoothfixconics}}
\begin{tabular}{||c|c|c|c||}
\hline\hline
&&&\\
Type&Matrix&$\E_t$&Isotropy\\&&&\\\hline&&&\\
{\rm i)}&
$\begin{pmatrix}
                    x_{0}&0\\
                    x_{2}&x_{1}\\
                    0&x_{0}+tx_{2}\end{pmatrix}$&
$\lav
x_{0}x_{1},x_{0}(x_{0}+tx_{2}),x_{2}(x_{0}+tx_{2})\rav$
&$\{id\}$\\&&&\\\hline&&&\\
{\rm ii)}&
$\begin{pmatrix}
                    x_{0}&0\\
                    x_{1}&x_{0}\\
                    0&x_{1}+tx_{2}\end{pmatrix}$&
$\lav
x_{0}^2,x_{0}(x_{1}+tx_{2}),x_{1}(x_{1}+tx_{2})\rav$
&$\{id\}$
\\&&&
\\
\hline\hline
\end{tabular}
\end{table}

If $[C]\in M^{1,1}-M^{0,2}\cap M^{1,1}$ is a union
of two lines, then $C$ is projectively equivalent to one of
the line pairs in Table~\ref{tablelinepairsM-M}. A line pair is
fixed if and only if both lines and the intersection point
are fixed. Hence the line pairs in Table~\ref{tablelinepairsM-M} are fixed
only for the four values
$\{(0,0),(0,\infty),(\infty,0),(\infty,\infty)\}$ of
$([a,b],[c,d])$. Table~\ref{tablelinepairsM-M} contains a
list of the line pairs,
and the isotropy subgroup of $S_{3}$ acting on the set
$\{(0,0),(0,\infty),(\infty,0),(\infty,\infty)\}$. This
takes care of all fixpoints in $M^{1,1}\setminus M^{0,2}\cap M^{1,1}$.

We now look at fixpoints in
$M^{0,2}\cap
M^{1,1}$. Suppose
$[f\:C\to\NN]\in M^{0,2}\cap M^{1,1}$ is fixed. Since
$M^{0,2}\cap M^{1,1}$ is a $\moo(\PP^2,2)$ fibration over
$I$, the map $f$ factors through
$\PP(\V_{(p,l)}^{\vee})\to\NN$ for some $(p,l)\in I$, which must
be fixed. Without loss of generality, we may assume that
$(p,l)$ is the point-line pair $(x_{1}^2,x_{1}x_{2})$.
The plane $\PP(\V_{(p,l)}^{\vee})$ contains no smooth fixconics,
hence $f(C)$ must be supported on the triangle of
fixlines 3)--5) in Table~\ref{tablefixlines}. If $f$ is a
double covering, then the branch points must be fixed.
For each of the fixlines 3)--5) there
are three possibilities for a 2:1 covering which is fixed:
If $C\simeq\PP^1$, then $f$ has two branch points, which
must be the fixpoints on the line. If
$C\simeq\PP^1\cup\PP^1$, then $f$ has one branch point,
which can be either of the two fixpoints on the line.
The isotropy subgroup of a fixpoint in $M^{0,2}\cap
M^{1,1}$ is trivial.

\bigskip

\subsection{Conics on the Calabi-Yau threefolds}\label{conicsCY}

Suppose $\B$ is one of the bundles $3\oh(p)$,
$\F^{\vee}\oplus\oh(2p)$, $S_{2}\F^{\vee}$, and $W$ is a
general Calabi-Yau sections of $\B$. Let
\begin{equation}
        \begin{CD}
                {}@.   \B\\
                @. @V{}VV \\
                \C @>e'>>\NN \\
                 @V{\pi}'VV \\
                H
        \end{CD}
\end{equation}
be the universal curve over $H$, let
$\B_{2}^{\prime}=\pi_{\ast}^{\prime}e^{\prime\ast}\B$, and let
$B_{2}^{\prime}=c_{9}(\B_{2}^{\prime})$. The virtual number of conics on
$W$ is by definition
\begin{equation}
n_{2}^W=\int_{H}B_{2}^{\prime}=
\int_{H^{1,1}}B_{2}^{\prime}+\int_{H^{0,2}}B_{2}^{\prime}\,,
\end{equation}
where
$\int_{H^{1,1}}B_{2}^{\prime}$ and $\int_{H^{0,2}}B_{2}^{\prime}$
are virtual numbers of $(1,1)$-conics and $(0,2)$-conics,
respectively.

Let $j^d\:\moo(W,d)\to\moo(\NN,d)$. Recall from the
Kleiman-Bertini argument in Section~\ref{linesCY} that if
$L\subseteq W$, then $j_{\ast}^1[L]\in\PP(\V)-\PP(\V_{I})$.
Suppose $[f\:C\to L]\in\moo(W,2)$ is a double covering of $L$.
Remark~2 of Section~\ref{Hilbert} implies that
$j_{\ast}^2[f\:C\to L]\in\moo(\NN,2)-M^{1,1}$. In other
words, there is no contribution from the double coverings to
$\int_{M^{1,1}}B_{2}$, where
$B_{2}=c_{9}(\pi_{\ast}e^{\ast}\B)$, hence
\[
\int_{M^{1,1}}B_{2}=\int_{H^{1,1}}B_{2}^{\prime}\,.\]

\begin{thm}
\label{thm14-1}
\quad
\begin{enumerate}
\item[i)] The number of conics on $X$ is 756. These are
all smooth of $(1,1)$-type.
\item[ii)] The virtual number of conics on $Y$ is 1674. All of
these are of $(1,1)$-type.
\item[iii)] The virtual number of conics on $Z$ is 504. Of
these are 468 of $(1,1)$-type and 36 of $(0,2)$-type.
\end{enumerate}
\end{thm}

\begin{proof}
We must evaluate
\\
\begin{equation}
\int_{M^{1,1}}B_{2}\qquad {\rm and}\qquad
\int_{H^{0,2}}B_{2}^{\prime}
\end{equation}
for each of the three vector bundles $3\oh(p)$,
$\F^{\vee}\oplus\oh(2p)$, and $S_{2}\F^{\vee}$. We shall use
Bott's formula to determine
$\int_{M^{1,1}}B_{2}$.
This computation is explained in the next section. The
integral $\int_{H^{0,2}}B_{2}^{\prime}$ will be evaluated in
the Chow ring $A^{\ast}(H^{0,2})$. This goes
exactly as our computation in Section~\ref{linesCY}. We
need to express the top Chern classes of (the
restrictions of) the vector bundles
\[\oh(p)_{2}^{\prime}\,,\quad \oh(2p)_{2}^{\prime}\,,\quad
\F_{2}^{\vee\prime}\,,\quad {\rm and}\quad (S_{2}\F^{\vee})_{2}
^{\prime}\]
in terms of the generators $\tau$, $\check{\tau}$, and $\xi$.
Let $\C_2$ denote the restriction of the universal curve
$\C$ to $H^{0,2}$. As seen in Section~\ref{modulspace}, the map
$\C_2\to\NN$ factors through the map $\PP(V^{\vee})\to\NN$, hence
\begin{equation}
        \begin{CD}
        \C @>{}>>\PP(\V^{\vee})@>{}>>\NN\\
        @V\pi' VV             @VV{}V@.\\
       H^{0,2} @>>{}>    PP^{\vee}
        \end{CD}
\end{equation}
is commutative. This identifies
\begin{equation}
\begin{split}
&\E_{\C_{2}}^{\vee}=pr_{\ast}\I_{\Tilde{p}}
\I_{\Tilde{l}}(2)_{\C_{2}}^{\vee}
\oplus\oh_{\C_{2}}(\check{\sigma})\\
&\F_{\C_{2}}^{\vee}\simeq pr_{\ast}\I_{\Tilde{l}}(1)_{\C_{2}}^{\vee}
\otimes\bigl(\wedge^2 pr_{\ast}\I_{\Tilde{p}}(1)_{\C_{2}}^{\vee}
\oplus\oh_{\C_{2}}(\check{\sigma})\bigr)\,,
\end{split}
\end{equation}
hence, on $H^{0,2}$, we have
\begin{equation}
\label{d3a}
\begin{split}
&\oh(p)_{2}^{\prime}=\pi_{\ast}^{\prime}\wedge^3\E_{\C_{2}}^{\vee}=
\wedge^2pr_{\ast}\I_{\Tilde{p}}\I_{\Tilde{l}}(2)^{\vee}
\otimes\pi_{\ast}^{\prime}\oh_{\C_{2}}(\check{\sigma})\\
&\oh(2p)_{2}^{\prime}=\pi_{\ast}^{\prime}\wedge^3\E_{\C_{2}}^{\vee\otimes 2}=
\wedge^2pr_{\ast}\I_{\Tilde{p}}\I_{\Tilde{l}}(2)^{\vee\otimes 2}
\otimes\pi_{\ast}^{\prime}\oh_{\C_{2}}(2\check{\sigma})\\
&\F_{2}^{\vee\prime} =pr_{\ast}\I_{\Tilde{l}}(1)^{\vee}
\otimes\bigl(\wedge^2 pr_{\ast}\I_{\Tilde{p}}(1)^{\vee}
\oplus\pi_{\ast}^{\prime}\oh_{\C_{2}}(\check{\sigma})\bigr)\\
&(S_{2}\F^{\vee})_{2}^{\prime}
=pr_{\ast}\I_{\Tilde{l}}(1)^{\vee\otimes 2}
\otimes\bigl(\wedge^2
pr_{\ast}\I_{\Tilde{p}}(1)^{\vee\otimes 2}
\oplus\wedge^2pr_{\ast}\I_{\Tilde{p}}(1)^{\vee}\otimes
\pi_{\ast}^{\prime}\oh_{\C_{2}}(\check{\sigma})
\oplus\pi_{\ast}^{\prime}\oh_{\C_{2}}(2\check{\sigma})\bigr)\,.
\end{split}
\end{equation}
There remains to determine
$\pi_{\ast}^{\prime}\oh_{\C_{2}}(\check{\sigma})$ and
$\pi_{\ast}^{\prime}\oh_{\C_{2}}(2\check{\sigma})$.
Consider the twisted ideal sequence
\begin{equation}
\label{d4}
0\to\I_{\C_{2}}(m\check{\sigma})\to\oh_{P}(m\check{\sigma})\to
\oh_{\C_{2}}(m\check{\sigma})\to 0\,,
\end{equation}
where $P=\PP(S_{2}\V^{\vee})\times_{PP^{\vee}}\PP(\V^{\vee})$.
Again, by a diagonal argument as in \eqref{diagonal},
$\I_{\C_{2}}\simeq\oh(-\xi)\otimes\oh(-2\check{\sigma})$, hence
pushing down \eqref{d4} with $m=1$ and $m=2$, we find
\begin{equation}
\V^{\vee}=\pi_{\ast}^{\prime}\oh_{\C_{2}}(\check{\sigma})
\end{equation}
and an exact sequence
\begin{equation}
0\to\oh(-\xi)\to S_{2}\V^{\vee}\to
\pi_{\ast}^{\prime}\oh_{\C_{2}}(2\check{\sigma})\to 0\,.
\end{equation}
Along with the expressions in Section~\ref{chowen}, this enables us to
express the Chern
classes of \eqref{d3a} in terms of the generators $\sigma$,
$\check{\sigma}$, and $\xi$.

Finally, the last assertions in i) follow (as in
Section~\ref{linesCY}) from Kleiman's Bertini theorem
and from the surjectivity of
\[H^0(\PP^{19},\oh_{\PP^{19}}(1))\to H^0(C,\oh_{C}(p))\]
for all conics $C\subseteq\PP^{19}$.
\end{proof}

\addtocontents{toc}{\vspace{0.2cm}}
\section{Computations. Bott's formula on the space of maps}
\label{calculations}
\bigskip

The integrands in the integrals of Theorem~\ref{thm7-1},
Lemma~\ref{lem9-1}, and
Theorem~\ref{thm14-1} are all polynomials in the
Chern classes of $T$-equivariant bundles. In Section~\ref{fixlines} and
Section~\ref{fixconics}  we saw that fixpoints on $\moo(\NN,1)$ and
$M^{1,1}$ were isolated. In
Section~\ref{representationsmotwo} we shall see that
fixpoints on $\smap_{0,2}(\NN,1)$ are also isolated.

Let $x_{i}$ denote the basis of $V$ with characters
$\lambda_{i}$ as in Section~\ref{torusaction}. If $U$ is a
finite dimensional representation of $T$ we shall let $[U]$
denote the class of $U$ in the (orbifold) representation ring
generated by $\lambda_{i}^s$ for non-zero
rational numbers $s$. Hence we write $[U]=\sum
a_{n_{0}n_{1}n_{2}}\lambda_{0}^{n_{0}}\lambda_{1}^{n_{1}}\lambda_{2}^{n_{2}}$
with $a_{n_{0}n_{1}n_{2}}\in\ZZ$ and $n_{i}\in\QQ$. We also
use the notation $[U]^{(k)}=\sum a_{n_{0}n_{1}n_{2}}
(\lambda_{0}^{n_{0}}\lambda_{1}^{n_{1}}\lambda_{2}^{n_{2}})^k$
 for $k\in\QQ$. Let $\CC^{\ast}\subseteq T$ be a
 one-parameter subgroup of $T$ with weights $\omega_{i}$.
 The induced weights of the $\CC^{\ast}$-action on $U$ are
 $n_{0}\omega_{0}+n_{1}\omega_{1}+n_{2}\omega_{2}$. Choose
 weights $\omega_{i}$ such that the induced weights on the
 tangent spaces of $\smap_{0,2}(\NN,1)$ (resp. $M_{1,1}$)
 at fixpoints are all non-zero. Then the fixpoints of the
 induced $\CC^{\ast}$-action on $\smap_{0,2}(\NN,1)$ (resp.
 $M_{1,1}$) are the same as the fixpoints of the $T$-action.
 In particular, they are isolated, hence we can apply the
 version of Bott's formula in the Introduction to evaluate
 the integrals.

 Suppose $\int_{M}P$ is any of the relevant integrals,
 where $P$ is a polynomial in the Chern classes of vector
 bundles $A$, $B$, etc.. To evaluate the integral using
 Bott's formula we must determine the weights of the vector
 spaces $T_{f}M$, $A_{f}$, $B_{f}$, etc. at each fixpoint
 $f\in M$. Our method will be as follows: For each fixpoint
 $f\in M$, we shall obtain formulas (in terms of the
 $\lambda_{i}$) for the tangent representation $[T_{f}M]$
 and the integrand representations $[A_{f}]$, $[B_{f}]$, etc..
 Once these are in store, they can readily be translated
 into the relevant lists of weights. In
 Section~\ref{representationsmotwo} we explain how to obtain
 formulas for the tangent and integrand representations of
 the integrals over $\smap_{0,2}(\NN,1)$ (Theorem~\ref{thm7-1}
 and Lemma~\ref{lem9-1}). In Section~\ref{representationsmoneone}
 we do the same for the integrals over $M^{1,1}$
 (Theorem~\ref{thm14-1}).
 The actual computation is carried out in
 {\texttt{MAPLE}} (see Appendix~\ref{maplecode}), where
 lists of fixpoints are constructed from our descriptions
 in Sections~\ref{fixlines} and \ref{fixconics}, along with
 the corresponding representations of
 Sections~\ref{representationsmotwo}
 and \ref{representationsmoneone}.

 Determining the tangent space representations is the most
 laborious. Let $(C,f,s_{i})$ be
 a fixed stable $n$-marked map of degree $d$ to $\NN$.
 We shall abuse notation and use $f$ to
 also denote the corresponding point in $\nk(\NN,d)$,
 i.e., $f=[C,f,s_{i}]$. Suppose
 $C=\cup C^{\alpha}$ and all $C^{\alpha}\simeq\PP^1$. By standard
 deformation arguments \cite{K,FP},
\begin{multline}
\label{c:0}
[T_{f}\nk(\NN,d)]=
[H^0(T\NN_{C})]-\sum_{\alpha}[H^0(TC^{\alpha})]+
\sum_{i=1}^{n}[T_{s_{i}}C]\\
+\sum_{\underset{\alpha\ne\beta}{y\in C^{\alpha}\cap
C^{\beta}}}
[T_{y}C^{\alpha}]\cdot [T_{y}C^{\beta}]+
[T_{y}C^{\alpha}]+ [T_{y}C^{\beta}]
\,.
\end{multline}
We shall only need to consider cases with $C\simeq\PP^1$, and
$C=\PP^1\cup\PP^1$.
To determine the terms of \eqref{c:0}, we shall use the
resolution \eqref{e:e5-1-1}, along with the component sequence
\begin{equation}
\label{c:1}
0\to T\NN_{C}\to T\NN_{C^1}\oplus T\NN_{C^2}\to T_{y}\NN\to
0\,,
\end{equation}
when $C\simeq C^1\cup C^2$ and $y=C^1\cap C^2$.  Also, if
$\L$ is a line bundle of degree one on $\PP^1$, then there
is an exact sequence
\begin{equation}
\label{c:2}
0\to\oh_{\PP^1}\to H^0(\L)^{\vee}\otimes\L\to T\PP^1\to
0\,.
\end{equation}
On $M^{1,1}$, matters are complicated by the fact that
$T_{f}M^{1,1}$ is only a proper subspace of
$T_{f}\moo(\NN,2)$ when $f\in M^{0,2}\cap M^{1,1}$.

\medskip

\noindent {\em Remark~1.} The number
$\lav T_{11},T_{12}\rav_{0,2}$, and the number of lines and
$(0,2)$-conics on the Calabi-Yau sections can also be
computed using Bott's formula. However, fixpoints
on $H^{0,2}$ and $M^{0,2}$ are not isolated, hence a
more general version of Bott´s formula than the one stated
in the introduction is required.

\medskip

\noindent{\em Remark~2.} Recall that we are using Bott's formula over the space
$M^{1,1}$ under the assumption that it is smooth. Even if this is not
true, our computation still holds. There is a theorem due
to Graber and Pandharipande \cite{GrP} which gives a
localization formula for integrals over virtual fundamental
classes. Applied to the space $\moo(\NN,2)$, the formula
says that the numerator in the contributions to Bott's
formula from singular points $f$ must be multiplied by the weight
of the $\CC^{\ast}$-action on the obstruction space
$H^1(f^{\ast}T\NN)$. A computation using this formula
reduces to the computation obtained by assuming
$M^{1,1}$ is smooth.

\medskip

\noindent{\em Remark~3.} We should not forget about the
automorphism group of a stable map. If
$[C,f,s_{1},s_{2}]\in\smap_{0,2}(\NN,1)$ then $\Aut(f)$ is
trivial. If $[C,f]\in M^{1,1}$, then $\Aut(f)$ is trivial
 unless $f$ is a $2:1$-covering, and if so $|\Aut(f)|=2$.

\bigskip

\subsection{Representations on $\smap_{0,2}(\NN,1)$}
\label{representationsmotwo}

First we describe the fixpoints on $\smap_{0,2}(\NN,1)$.
Suppose $f=[C,f,s_{1},s_{2}]\in\smap_{0,2}(\NN,1)$ is a
fixpoint. Then the line $f(C)$ and the points $f(s_{1})$
and $f(s_{2})$ must be fixed. There are essentially two
possibilities:

\begin{enumerate}
\item[i)] $f(s_{1})\ne f(s_{2})$. Then $C\simeq\PP^1$ and
$f(s_{1})$, $f(s_{2})$ are two fixpoints on $f(C)$.
\item[ii)] $f(s_{1})= f(s_{2})$. Then $C$ is the union of
two $\PP^1$-components $C^1$ and $C^2$, such that $f(C^1)$
is a fixline, and $C^2$ contains the two marked points.
Hence $f(C^2)$ is a fixpoint on $f(C^1)$.
\end{enumerate}
In other words, for each fixline $L$ in $\NN$ with
fixpoints $f_{1},f_{2}\in L$, there are four fixpoints in
$\smap_{0,2}(\NN,1)$ corresponding to the configurations
$(L,f_{1},f_{2})$, $(L,f_{2},f_{1})$, $(L,f_{1})$, and
$(L,f_{2})$. In
particular, the fixpoints on $\smap_{0,2}(\NN,1)$ are
isolated.

The integrals in Theorem~\ref{thm7-1} and Lemma~\ref{lem9-1}
are of type
\begin{equation*}
\int_{\smap_{0,2}(\NN,1)}\!\!e_{1}^{\ast}(T_{r})e_{2}^{\ast}(T_{s})\qquad
{\rm and}\qquad
\int_{\smap_{0,2}(\NN,1)}\!
\frac{e_{1}^{\ast}(T_{r})}{1-c}\,,
\end{equation*}
where $0\le r,s\le 12$. The integrands are polynomials in
the Chern classes of the bundles $\L_{1}$,
$e_{k}^{\ast}(\E)$, and $e_{k}^{\ast}(\F)$, for $k=1,2$,
whose fibers over $f$ are $T_{s_{1}}C^{\vee}$,
$\E_{f(s_{k})}$, and $\F_{f(s_{k})}$, respectively. Since
$f(s_{k})$ is a fixpoint in $\NN$, the representations
$[\E_{f(s_{k})}]$ and $[\F_{f(s_{k})}]$ are immediately
obtained from Table~\ref{tablefixpoints}. To determine the
representations $[T_{s_{1}}C^{\vee}]$ and
$[T_{f}\smap_{0,2}(\NN,1)]$, we must consider the fixpoint
types i) and ii) separately. Suppose $f(C)=L$, a fixline
with fixpoints $f_{1}$ and $f_{2}$.
\\\\
{\em Case~i).} Then
$[T_{s_{1}}C^{\vee}]=[T_{f(s_{1})}L]^{-1}$, and \eqref{c:0}
yields
\begin{equation}
\label{c:-1}
[T_{f}\smap_{0,2}(\NN,2)]=[H^0(T\NN_{L})]-[H^0(TL)]+[T_{f_{1}}(L)]+
[T_{f_{2}}(L)]\,.
\end{equation}\\
{\em Case~ii).} Since $C^2$ is contracted by $f$, the
induced action on $C^2$ is trivial and
$H^0(T\NN_{C^2})=H^0(T\NN_{y})$, where $y=C^1\cap C^2$.
Hence $[T_{s_{1}}C^{\vee}]=1=[T_{y}C^2]$, and \eqref{c:0},
\eqref{c:1} yields
\begin{equation}
\label{c:-2}
[T_{f}\smap_{0,2}(\NN,1)]=[H^0(T\NN_{L})]-[H^0(TL)]+2[T_{y}L]\,.
\end{equation}
Note that since $[T_{s_{1}}C^{\vee}]=1$, the
contributions from these points to the integral
 $\int_{\smap_{0,2}(\NN,1)}
\frac{e_{1}^{\ast}(T_{r})}{1-c}$ is zero.

We shall express the different parts of
\eqref{c:-1} and \eqref{c:-2} in terms of the {\em elementary}
representations
\begin{equation}
\label{c:3}
[H^0(\I_{p}(1))]\,,\qquad[H^0(\I_{l}(1))]\,,\,\footnotemark\qquad
[\oh_{f_{i}}(\check{\sigma})]\,,\qquad{\rm and }\qquad
[H^0(\oh_{L}(\check{\sigma}))]\,,
\end{equation}
\footnotetext{Here
$H^0(\I_{l}(1))=H^0(\PP(V),\I_{l}(1))$, where
$l\in\PP(V^{\vee})$ is identified with the line
$l\subseteq\PP(V)$.}
where $\rho([L])=(p,l)$. From Remark 2) in
Section~\ref{modulspace}, there are $T$-equivariant
isomorphisms:
\begin{equation}
\label{c:4}
\begin{split}
\E_{L}=&H^0(\I_{p}\I_{l}(2))\oplus\oh(-\check{\sigma})\\
\F_{L}\simeq& H^0(\I_{l}(1))\otimes\bigl(\wedge^2
H^0(\I_{p}(1))\oplus
\oh_{L}(-\check{\sigma})\bigr)\,.
\end{split}
\end{equation}
For each fixline, the representations \eqref{c:3} can
be read off Table~\ref{tablefixlines}.

As $\oh_{L}(\check{\sigma})$ is of degree one on $L$, there
is a $T$-equivariant exact sequence
\begin{equation}
0\to\oh_{L}\to H^0(\oh_{L}(\check{\sigma}))^{\vee}\otimes
\oh_{L}(\check{\sigma})\to TL\to 0\,,
\end{equation}
hence
\begin{equation}
[T_{f_{i}}L]=[\oh_{f_{j}}(\check{\sigma})]^{-1}\cdot
[\oh_{f_{i}}(\check{\sigma})]\,,
\end{equation}
with $i\ne j$, and
\begin{equation}
[H^0(TL)]=[H^0(\oh_{L}(\check{\sigma}))]^{(-1)}\cdot
[H^0(\oh_{L}(\check{\sigma}))]-1\,.
\end{equation}
Since
$H^1(\E_{L}^{\vee}\otimes\E_{L})=0=H^1(\F_{L}^{\vee}\otimes\F_{L})$,
the exact sequence \eqref{e:e5-1-1} implies
\begin{equation}
\label{c:8}
[H^0(T\NN_{L})]=[H^0(\F_{L}^{\vee}\otimes\E_{L})]\cdot
(\lambda_{0}+\lambda_{1}+\lambda_{2})
-[H^0(\E_{L}^{\vee}\otimes\E_{L})]-[H^0(\F_{L}^{\vee}\otimes\F_{L})]+1\,.
\end{equation}
The terms in this expression can be determined by using the
identifications \eqref{c:4} and
$[H^0(\I_{p}\I_{l}(2))]=[H^0(\I_{p}(1))]\cdot[H^0(\I_{l}(1))]$.
This yields
\begin{equation}
\label{c:9}
\begin{split}
[H^0(\F_{L}^{\vee}\otimes\E_{L})]=&
[H^0(\I_{p}(1))]^{(-1)}+[H^0(\I_{l}(1))]^{(-1)}+
[H^0(\oh_{L}(\check{\sigma}))]\cdot[H^0(\I_{p}(1))]\\
[H^0(\E_{L}^{\vee}\otimes\E_{L})]=&
[H^0(\I_{p}(1))]^{(-1)}\cdot[H^0(\I_{p}(1))]+
[H^0(\oh_{L}(\check{\sigma}))][H^0(\I_{p}\I_{l}(2))]+1\\
[H^0(\F_{L}^{\vee}\otimes\F_{L})]=&
[H^0(\oh_{L}(\check{\sigma}))]\cdot[\wedge^2H^0(\I_{p}(1))]\,.
\end{split}
\end{equation}
In the next section we shall need the representations
$[H^0(\wedge^3\E_{L}^{\vee})]$, $[H^0(\F_{L}^{\vee})]$,
$[H^0(\wedge^3\E_{L}^{\vee\otimes 2})]$, and
$[H^0(S_{2}\F_{L}^{\vee})]$. Again, these are obtained from
\eqref{c:4}:
\begin{equation}
\label{c:-3}
\begin{split}
[H^0(\wedge^3\E_{L}^{\vee})]=&
[\wedge^2H^0(\I_{p}(1))]^{-1}\cdot[H^0(\I_{l}(1))]^{-2}
\cdot[H^0(\oh_{L}(\check{\sigma}))]\\
[H^0(\F_{L}^{\vee})]=&
[H^0(\I_{l}(1))]^{-1}\cdot\bigl([\wedge^2H^0(\I_{p}(1))]^{-1}+
[H^0(\oh_{L}(\check{\sigma}))]\bigr)\\
[H^0(\wedge^3\E_{L}^{\vee\otimes 2})]=&
[S_{2}H^0(\wedge^3\E_{L}^{\vee})]\\
[H^0(S_{2}\F_{L}^{\vee})]=&
[S_{2}H^0(\F_{L}^{\vee})]\,.
\end{split}
\end{equation}

\bigskip

\subsection{Representations on $M^{1,1}$}
\label{representationsmoneone}

The integrals in Theorem~\ref{thm14-1} can be expressed as
\[\int_{M^{1,1}}c_{3}(\oh(p)_{2})^3\,,\qquad
\int_{M^{1,1}}c_{4}(\F_{2}^{\vee})c_{5}(\oh(2p)_{2})\,,\qquad
{\rm and}\qquad
\int_{M^{1,1}}c_{9}(S_{2}\F_{2}^{\vee})_{2}\,,\]
hence for each fixpoint $f=[f\:C\to\NN]\in M^{1,1}$, we must
determine $[T_{f}M^{1,1}]$, and the integrand representations
\begin{equation}
\label{c:10}
[H^0\oh_{C}(p)],\quad
[H^0(\F_{C}^{\vee})],\quad[H^0(\oh_{C}(2p))],\quad{\rm
and}\quad[H^0(S_{2}\F_{C}^{\vee})]\,.
\end{equation}
The fixpoints of $M^{1,1}$ are described in Section~\ref{fixconics}.
We determine the representations \eqref{c:10} for
each type separately. First, we consider the fixpoints in
$M^{1,1}-M^{0,2}\cap M^{1,1}$. As seen, if $f\in M^{1,1}-
M^{0,2}\cap M^{1,1}$, then $f$ is an embedding and
$C\simeq\PP^1$, or $C\simeq\PP^1\cup\PP^1$, and
$T_{f}M^{1,1}=T_{f}\moo(\NN,2)$.
\\\\
\underline{$M^{1,1}-M^{0,2}\cap M^{1,1}$. $\quad C\simeq\PP^1$}. If $f\in
M^{1,1}-M^{0,2}\cap M^{1,1}$ and $C\simeq\PP^1$, then we may
assume $f(C)$ is one of the smooth conics in
Table~\ref{tablesmoothfixconics}.
For each type, consider its matrix representation.
Let $a_{ij}$ denote the $(i,j)$-entry and let $\L$ be the
line bundle on $C$ with fibers $\L_{t}=a_{32}(t)$. Then the
bundles $\E_{C}$ and $\F_{C}$ decomposes as:
\begin{equation}
\label{c:20a}
\begin{split}
&\E_{C}=a_{11}a_{22}\oplus a_{11}\L\oplus
a_{21}\L\\
&\F_{C}\simeq a_{11}a_{21}\L\oplus a_{11}a_{22}\L\,,
\end{split}
\end{equation}
where we abuse notation and write
$a_{11}a_{22}$, $a_{11}\L_{t}$, etc., for the bundles with
fibers $a_{11}a_{22}$, $a_{11}\L_{t}$, etc..
The isomorphism follows by observing that the map
\begin{equation*}
a_{11}a_{21}\L\oplus a_{11}a_{22}\L\longrightarrow \E_{C}\otimes V\,,
\end{equation*}
given on the fibers by
\begin{align*}
&a_{11}a_{21}\L_{t}\longmapsto a_{11}\L_{t}\otimes a_{21}-
a_{21}\L_{t}\otimes a_{11}
\\
&a_{11}a_{22}\L_{t}\longmapsto a_{11}a_{22}\otimes\L_{t}-
a_{11}\L_{t}\otimes a_{22}
\end{align*}
is injective, and lies in the kernel of the
multiplication map
$S_{2}V\otimes V\to S_{3}V$.

 From \eqref{c:0} we have
\begin{equation}
[T_{f}M^{1,1}]=[H^0(T\NN_{C})]-[H^0(TC)]\,.
\end{equation}
The line bundle $\L^{\vee}$ is of degree one on
$C\simeq\PP^1$, hence the exact sequence of \eqref{c:2}
yields
\begin{equation}
[H^0(TC)]=
[H^0(\L^{\vee})]^{(-1)}\cdot[H^0(\L^{\vee})]-1\,.
\end{equation}
Since $H^1(\E_{C}^{\vee}\otimes\E_{C})=0=
H^1(\F_{C}^{\vee}\otimes\F_{C})$,
\eqref{c:8} still holds with $L=C$. Using \eqref{c:20a}, we
find the following
representations for the terms in \eqref{c:8}:
\begin{equation}
\begin{split}
[H^0(\F_{C}^{\vee}\otimes\E_{C})]=&[H^0(\L^{\vee})]
(1+[a_{21}]^{-1}[a_{22}])+ [a_{21}]^{-1}+[a_{11}]^{-1}+
[a_{22}]^{-1}\\
&+[a_{11}]^{-1}[a_{22}]^{-1}[a_{21}]\\
[H^0(\E_{C}^{\vee}\otimes\E_{C})]=&[H^0(\L^{\vee})]
([a_{22}]+[a_{11}][a_{22}][a_{21}]^{-1})+[a_{21}]^{-1}[a_{11}]\\
&+[a_{11}]^{-1}[a_{21}]+3\\
[H^0(\F_{C}^{\vee}\otimes\F_{C})]=&[a_{21}]^{-1}[a_{22}]+
[a_{22}]^{-1}[a_{21}]+2\,.
\end{split}
\end{equation}

Since $\oh_{C}(p)=\wedge^3\E_{C}^{\vee}$,
$\oh_{C}(2p)=\wedge^3\E_{C}^{\vee\otimes 2}$, and
$S_{2}\F_{C}^{\vee}\simeq S_{2}(a_{21}\oplus
a_{22})^{\vee}\oplus a_{11}^{2\vee}\oplus\L^{\vee\otimes 2}$,
the representations \eqref{c:10} are:
\begin{equation*}
\begin{split}
[H^0(\oh_{C}(p))]=&[a_{11}]^{-2}[a_{22}]^{-1}[a_{21}]^{-1}
[S_{2}H^0(\L^{\vee})]\\
[H^0(\oh_{C}(2p))]=&[a_{11}]^{-4}[a_{22}]^{-2}[a_{21}]^{-2}
[S_{4}H^0(\L^{\vee})]\\
[H^0(\F_{C}^{\vee})]=&[a_{11}]^{-1}[a_{21}
\oplus a_{22}]^{(-1)}[H^0(\L^{\vee})]\\
[H^0(S_{2}\F_{C}^{\vee})]=&[S_{2}(a_{21}\oplus a_{22})]^{(-1)}
[a_{11}]^{-2}\cdot[S_{2}H^0(\L^{\vee})]\,.
\end{split}
\end{equation*}
In these formulas the representations $[a_{ij}]$ and
$[H^0(\L^{\vee})]$ can be read off from
Table~\ref{tablesmoothfixconics}.
\\\\
\underline{ $M^{11}-M^{0,2}\cap M^{11}$.
$\quad C\simeq \PP^1\cup\PP^1$.}
 If $f\in M^{1,1}-M^{0,2}\cap
M^{1,1}$ and $C\simeq C^1\cup C^2$ with
$C^1\simeq\PP^1\simeq C^2$, then we may assume that $f(C)$
is one of the fixed line pairs in
Table~\ref{tablelinepairsM-M}.
Let $f^i$ be the restriction of $f$ to $C^i$. Applying \eqref{c:1},
we find
\begin{multline}
\label{c:14}
[T_{f}M^{11}]=[T_{f^1}\moo(\NN,1)]+[T_{f^2}\moo(\NN,1)]+
[T_{y}C^1]\cdot[T_{y}C^2]\\+[T_{y}C^1]
+[T_{y}C^2]-[T_{y}\NN]\,,
\end{multline}
Formulas for
all terms, except $[T_{y}\NN]$, in \eqref{c:14} were
essentially found in Section~\ref{representationsmotwo}.
The exact sequence \eqref{e:e5-1-1} yields
\begin{equation}
[T_{y}\NN]=[\F_{y}]^{(-1)}[\E_{y}]
(\lambda_{0}+\lambda_{1}+\lambda_{2})
-[\E_{y}]^{(-1)}[\E_{y}]-[\F_{y}]^{(-1)}[\F_{y}]+1\,.
\end{equation}
\\
The representations in these expressions obtained from
Table~\ref{tablefixpoints}.

The integrand representations \eqref{c:10} are found using
component sequences similar to \eqref{c:1} on the integrand
bundles:
\begin{equation}
\label{c:-4}
\begin{split}
[H^0(\oh_{C}(p))]=&[H^0(\oh_{C^1}(p))]+[H^0(\oh_{C^2}(p))]-
[\oh_{y}(p)]\\
[H^0(\oh_{C}(2p))]=&[H^0(\oh_{C^1}(2p))]+[H^0(\oh_{C^2}(2p))]-
[\oh_{y}(p)]^2\\
[H^0(\F_{C}^{\vee})]=&[H^0(\F_{C^1}^{\vee})]+
[H^0(\F_{C^2}^{\vee})]-[\F_{y}]^{(-1)}\\
[H^0(S_{2}\F_{C}^{\vee})]=&[H^0(S_{2}\F_{C^1}^{\vee})]+
[H^0(S_{2}\F_{C^2}^{\vee})]-[S_{2}\F_{y}]^{(-1)}\,,
\end{split}
\end{equation}
where all representations on the right-hand sides can be
obtained from \eqref{c:-3} and Table~\ref{tablefixlines}.\\

We now turn our attention to fixpoints $f$ in $M^{0,2}\cap
M^{1,1}$. By our reasoning in Section~\ref{modulspace}, we
may assume without loss of generality that $f\:C\to\NN$ factors
through $\PP(\V_{(p,l)}^{\vee})\to\NN$, where $(p,l)$ is
the point-line pair $(x_{1}^2,x_{1}x_{2})$.
Since we are assuming that $M^{0,2}$ and $M^{1,1}$
intersect transversally,
\begin{equation}
[T_{f}M^{1,1}]=[T_{f}\moo(\NN,2)]-[N_{f}]\,,
\end{equation}
where $N$ is the normal bundle
of $M^{0,2}\cap M^{1,1}$ in $M^{0,2}$. From the fiber square
in Remark~2, Section~\ref{Hilbert}, we obtain a $T$-equivariant
diagram of exact sequences on $M^{0,2}\cap M^{1,1}$
\begin{equation*}
        \begin{CD}
          {} @.{} @. 0 @. 0 @.{} \\
         @. @. @V{}VV  @V{}VV         @. \\
          0 @>{}>>TR @>{}>> T(M^{0,2}\cap M^{1,1})
          @>{}>>\rho^{\ast} TI @>{}>> 0 \\
          @. @|  @V{}VV  @V{}VV         @.\\
            0 @>{}>>TR @>{}>> TM^{0,2}_{M^{0,2}\cap M^{1,1}}
            @>{}>>\rho^{\ast}
            TPP_{M^{0,2}\cap M^{1,1}}^{\vee}
            @>{}>>0\\
             @. @. @V{}VV  @V{}VV         @.  \\
               {} @.{} @. N_{f} @=\rho^{\ast}N_{I/PP^{\vee}}@.{}\\
         @. @. @V{}VV  @V{}VV         @.\\
          {} @.{} @. 0 @. 0 @.{}
        \end{CD}
\end{equation*}
\\
where $TR$ is the
relative tangent bundle of $\rho$.
The fibers of $\rho^{\ast}N_{I/PP^{\vee}}$ are
$\Hom(\oh_{p}(-\tau),\oh_{l}(\check{\tau}))$, so with our
choice of $(p,l)$,
\begin{equation}
[N_{f}]=\lambda_{0}\lambda_{1}^{-1}\,.
\end{equation}

We shall discuss the representations separately for the three
types of fixpoints in $M^{0,2}\cap M^{1,1}$. However, before
we start, note that
if $L$ is a line of type 3)--5) in
Table~\ref{tablefixlines}, then
$T\NN_{L}\simeq\oh_{\PP^1}(2)\oplus 2\oh_{\PP^1}(1)\oplus
2\oh_{\PP^1}\oplus\oh_{\PP^1}(-1)$. In Appendix~\ref{A-tangent},
we show that the line bundle of degree $-1$ in the
split is $T$-equivariant to the line bundle
\begin{equation}
\label{c:18}
H^0(\I_{l}(1))^{\vee}\otimes \wedge^2H^0(\I_{p}(1))^{\vee}
\otimes V/H^0(\I_{p}(1))\otimes
\oh_{L}(-\check{\sigma})\,.
\end{equation}
\\\\
\underline{$M^{0,2}\cap M^{1,1}$.
$\quad C\simeq C^1\cup C^2$ and $f$ an embedding.}
Assume $f$ maps $C\simeq C^1\cup C^2$ onto two of the three
lines 3)--5) in Table~\ref{tablefixlines}. Then
\eqref{c:1} and \eqref{c:0} yields
\begin{equation}
\label{c:18a}
\begin{split}
[H^0(T\NN_{C})]=&[H^0(T\NN_{C^1})]+[H^0(T\NN_{C^2})]-
[T_{y}\NN]+[H^1(T\NN_{C})]\,,\\
[T_{f}\moo(\NN,2)]=[T_{f^1}\moo(\NN,1)]&+[T_{f^2}\moo(\NN,1)]+
[T_{y}C^1][T_{y}C^2]\\
&+[T_{y}C^1]+[T_{y}C^2]-[T_{y}\NN]+[H^1(T\NN_{C})]\,,
\end{split}
\end{equation}
where $f^i\:C^i\to\NN$ are restriction maps. Only
$[H^1(T\NN_{C})]$ remains to be determined.
Since $H^1(T\NN_{C})$ is the cokernel of
\[H^0(T\NN_{C^1})\oplus H^0(T\NN_{C^2})\to T_{y}N\,,\]
it follows from \eqref{c:18} that, with our choice of $(p,l)$,
\begin{equation}
[H^1(T\NN_{C})]=\lambda_{0}\lambda_{1}^{-2}\lambda_{2}^{-1}
[\oh_{y}(-\check{\sigma})]\,.
\end{equation}

The integrand representations for these
fixpoints receive the same formulas as in \eqref{c:-4}.
\\\\
\noindent\underline{$M^{0,2}\cap M^{11}$.$\quad$ $2:1$ covering with
$C\simeq C^1\cup C^2$.}
This is similar to the previous case, except that the
restriction maps $f^1$ and $f^2$ are the same. Hence we
obtain the formulas \eqref{c:18a} and \eqref{c:-4} with
$C^1=C^2$ and $f^1=f^2$.
\\\\
\noindent\underline{$M^{0,2}\cap M^{11}$.$\quad$ $2:1$ covering with
$C\simeq \PP^1$.}
Let $\oh(1)$ be the tautological line bundle on $C$. If
$f\:C\to L$ is a $2:1$ covering, then
$f^{\ast}\oh_{L}(\check{\sigma})=\oh(2)$, and $f$ induces a
$T$-action on $\oh(1)$ with
$[H^0(\oh(1))]=[H^0(\oh_{L}(\check{\sigma}))]^{(1/2)}$.
By \eqref{c:0}, we have
\[[T_{f}\moo(\NN,2)]=[H^0(T\NN_{C})]-
[H^0(TC)]\,.\]
Since $\oh(1)$ is of degree $1$ on $C$,
\[ [H^0(TC)]=[H^0(\oh_{L}(\check{\sigma}))]^{(1/2)}
[H^0(\oh_{L}(\check{\sigma}))]^{(-1/2)}-1\]
The representation $H^0(T\NN_{C})]$ and the integrand
representations are found by using the pullbacks of \eqref{c:4}
and \eqref{c:18} to $C$.
One may check that
\[[H^1(\F_{C}^{\vee}\otimes\E_{C})](\lambda_{0}+\lambda_{1}+
\lambda_{2})=
[H^1(T\NN_{C})]+[H^1(\E_{C}^{\vee}\otimes\E_{C})]\,,\]
so the exact sequence \eqref{e:e5-1-1} yields
\begin{multline}
\label{c:20}
[H^0(T\NN_{C})]=[H^0(\F_{C}^{\vee}\otimes\E_{C})](\lambda_{0}+
\lambda_{1}+\lambda_{2})-[H^0(\E_{C}^{\vee}\otimes\E_{C})]-
[H^0(\F_{C}^{\vee}\otimes\F_{C})]\\
+1+
[H^1(\F_{C}^{\vee}\otimes\F_{C})]\,.
\end{multline}
\\
Here $[H^1(\F_{C}^{\vee}\otimes\F_{C})]=
\lambda_{1}^{-1}\lambda_{2}^{-1}
[\wedge^2H^0(\oh_{L}(\check{\sigma}))]^{-1/2}$. All the
other terms in \eqref{c:20} are obtained by replacing
$[H^0(\oh_{L}(\check{\sigma}))]$ with
$[H^0(\oh(2))]=[S_{2}H^0(\oh_{L}(\check{\sigma}))]^{(1/2)}$
in \eqref{c:9}. For the integrand representations, we obtain:
\begin{equation*}
\begin{split}
[H^0(\oh_{C}(p))]=&[\wedge^2H^0(\I_p(1))]^{-1}\cdot
[H^0(\I_l(1))]^{-2}\cdot[S_{2}H^0(\oh_{L}(\check{\sigma}))]^{(1/2)}\\
[H^0(\F_{C}^{\vee})]=&[H^0(\I_l(1))]^{-1}\cdot
([\wedge^2H^0(\I_p(1))]^{-1}+[S_{2}H^0(\oh_{L}(\check{\sigma}))]^{(1/2)})\\
[H^0(\oh_{C}(2p))]=&[\wedge^2H^0(\I_p(1))]^{-2}\cdot
[H^0(\I_l(1))]^{-4}\cdot[S_{4}H^0(\oh_{L}(\check{\sigma}))]^{(1/2)}\\
[H^0(S_2\F_{C}^{\vee})]=&[H^0(\I_l(1))]^{-2}\cdot
\bigl([\wedge^2H^0(\I_p(1))]^{-2}+[\wedge^2H^0(\I_{p}(1))]^{-1}
\cdot[S_{2}H^0(\oh_{L}(\check{\sigma}))]^{(1/2)}\\
&\qquad\qquad\qquad\qquad\qquad\qquad\qquad\qquad+
[S_{4}H^0(\oh_{L}(\check{\sigma}))]^{(1/2)}
\bigr)\,.
\end{split}
\end{equation*}

\addtocontents{toc}{\vspace{0.2cm}}
\appendix
\section{Tangent bundles over lines}
\label{A-tangent}
\bigskip

In this appendix, we compute the restriction of the tangent
bundle $T\NN$ to lines $L$ on $\NN$.

If $L$ is a line in $\NN$ with homogeneous coordinates
$[u,t]$, then the differential map
\[\HomC(\F_{L},\E_{L}\otimes V)\to T\NN_{L}\] can be
represented by a $6\times 6$-matrix with entries in
$V^{\vee}=\Hom(V,\CC)$ and $\CC[u,t]$. For lines 1)--6) in
Section~\ref{fixlines}, we can decompose the vector bundles
\begin{equation}
\label{A0}
\E_{L}\simeq\oh_{L}\oplus\oh_{L}\oplus\oh_{L}(-1)\qquad
{\rm and}\qquad
\F_{L}\simeq\oh_{L}\oplus\oh_{L}(-1)\,,
\end{equation}
such that the universal maps $\A_{L}\:\F_{L}\to\E_{L}\otimes
V$ are represented by the matrices\footnotemark
\footnotetext{Note that, up to a reparametrization of $L$,
the $2\times 2$-minors of the matrices in
Table~\ref{tableuniversalmap} generate the nets in
Section~\ref{fixlines}.}
in Table~\ref{tableuniversalmap}.
\begin{table}[h!]
\centering
\caption{Universal map\label{tableuniversalmap}}
$\!\!\!\!\!\!\!\!$\begin{tabular}
{||c|cccccc||}
\hline\hline
Type&1)&2)&3)&4)&5)&6)\\\hline
$\A_{L}$&$\begin{pmatrix}
x_{2}&ux_{1}\\
x_{1}&tx_{2}\\
0&-x_{0}
\end{pmatrix}$
&$\!\!\!\!\!\!\begin{pmatrix}
x_{2}&ux_{2}+tx_{1}\\
x_{1}&0\\
0&-x_{0}
\end{pmatrix}$
&$\!\!\!\!\!\!\begin{pmatrix}
x_{2}&ux_{0}\\
x_{1}&tx_{2}\\
0&-x_{1}
\end{pmatrix}$
&$\!\!\!\!\!\!\!\!\begin{pmatrix}
x_{2}&0\\
x_{1}&ux_{2}+tx_{0}\\
0&-x_{1}
\end{pmatrix}$
&$\!\!\!\!\!\!\begin{pmatrix}
x_{2}&tx_{0}\\
x_{1}&ux_{0}\\
0&-x_{1}
\end{pmatrix}$
&$\!\!\!\!\!\!\begin{pmatrix}
x_{2}&ux_{0}\\
x_{1}&ux_{2}+tx_{0}\\
0&-x_{1}
\end{pmatrix}$
\\\hline\hline
\end{tabular}
\end{table}

 Let $e_{1},e_{2},e_{3}$,
and $f_{1},f_{2}$, respectively, be bases for the graded
$\CC[u,t]$-modules in \eqref{A0}, with $\deg e_{1}=\deg e_{2}=
\deg f_{1}=0$ and $\deg e_{3}=
\deg f_{2}=1$. If we represent elements in the
$\CC[u,t]$-module
$\HomC(\F_{L},\E_{L})\simeq\F_{L}^{\vee}\otimes\E_{L}$ as
\begin{equation}
\label{A1}
\begin{pmatrix}
a_{11}\\a_{12}\\a_{13}\\a_{21}\\a_{22}\\a_{23}
\end{pmatrix}=\sum a_{ij}f_{i}^{\vee}\otimes e_{j}\,,
\end{equation}
then the differential maps over lines 1)--6) can be
represented by the matrices in
Table~\ref{tabledifferentialmap}.
The procedure for finding the matrices in
Table~\ref{tabledifferentialmap} is as follows:
\begin{table}[b]
\centering
\caption{Differential map\label{tabledifferentialmap}}
\begin{tabular}{||c|c||}
\hline\hline
Type&$D$\\\hline
1)&$\begin{pmatrix}
x_{0}^{\vee}&0&0&0&0&0\\
0&x_{0}^{\vee}&0&0&0&0\\
0&0&ux_{2}^{\vee}&0&0&-x_{1}^{\vee}\\
-tx_{1}^{\vee}&ux_{2}^{\vee}&0&x_{2}^{\vee}&-x_{1}^{\vee}&0\\
utx_{2}^{\vee}&-utx_{1}^{\vee}&0&-tx_{1}^{\vee}&ux_{2}^{\vee}&0\\
0&0&-tx_{1}^{\vee}&0&0&x_{2}^{\vee}
\end{pmatrix}$
\\\hline
2)&$\begin{pmatrix}
x_{0}^{\vee}&0&0&0&0&0\\
0&x_{0}^{\vee}&0&0&0&0\\
0&0&tx_{2}^{\vee}&0&0&-x_{1}^{\vee}\\
0&t^2x_{2}^{\vee}&0&tx_{2}^{\vee}-ux_{1}^{\vee}&-tx_{1}^{\vee}&0\\
0&x_{1}^{\vee}&0&tx_{2}^{\vee}-ux_{1}^{\vee}&-tx_{1}^{\vee}&0\\
0&0&0&0&0&tx_{2}^{\vee}-ux_{1}^{\vee}
\end{pmatrix}$
\\\hline
3)&$\begin{pmatrix}
utx_{2}^{\vee}&-utx_{1}^{\vee}&0&-tx_{0}^{\vee}&ux_{2}^{\vee}&0\\
0&x_{0}^{\vee}&0&0&0&0\\
0&0&x_{0}^{\vee}&0&0&0\\
tx_{0}^{\vee}&-ux_{2}^{\vee}&0&0&x_{0}^{\vee}&0\\
0&-ux_{2}^{\vee}&-utx_{1}^{\vee}&0&x_{0}^{\vee}&ux_{2}^{\vee}\\
0&0&ux_{2}^{\vee}&0&0&x_{0}^{\vee}
\end{pmatrix}$
\\\hline
4)&$\begin{pmatrix}
x_{0}^{\vee}&0&0&0&0&0\\
0&x_{0}^{\vee}&tx_{1}^{\vee}&0&0&-x_{0}^{\vee}\\
0&0&x_{0}^{\vee}&0&0&0\\
-tx_{1}^{\vee}&0&0&x_{0}^{\vee}&0&0\\
0&0&0&0&tx_{2}^{\vee}-ux_{0}^{\vee}&0\\
0&0&0&0&0&tx_{2}^{\vee}-ux_{0}^{\vee}
\end{pmatrix}$
\\\hline
5)&$\begin{pmatrix}
utx_{0}^{\vee}&-tx_{0}^{\vee}&0&0&0&0\\
0&x_{0}^{\vee}&ux_{1}^{\vee}+tx_{2}^{\vee}&0&0&-x_{0}^{\vee}\\
0&0&x_{0}^{\vee}&0&0&0\\
-u^2x_{1}^{\vee}&0&0&ux_{0}^{\vee}&0&-tx_{0}^{\vee}\\
0&0&0&0&x_{2}^{\vee}&0\\
0&0&0&0&0&x_{2}^{\vee}
\end{pmatrix}$
\\\hline
6)&$\begin{pmatrix}
tx_{0}^{\vee}&-ux_{0}^{\vee}&0&0&0&0\\
x_{0}^{\vee}&0&ux_{1}^{\vee}&0&0&-x_{2}^{\vee}\\
0&0&x_{0}^{\vee}&0&0&0\\
-ux_{2}^{\vee}-tx_{1}^{\vee}&ux_{1}^{\vee}&0&x_{0}^{\vee}&-x_{2}^{\vee}&0\\
-u^2x_{0}^{\vee}&u^2x_{2}^{\vee}&0&0&tx_{2}^{\vee}-ux_{0}^{\vee}&0\\
0&0&u^2x_{2}^{\vee}&0&0&tx_{2}^{\vee}-ux_{0}^{\vee}
\end{pmatrix}$
\\\hline\hline
\end{tabular}
\end{table}

\medskip

1) For each line $L$, use the matrix representation in
Table~\ref{tableuniversalmap} to write down a $6\times
13$-matrix $M$, with entries in $V$ and $\CC[u,t]$,
representing the map
\begin{gather*}
\CEnd(\E_{L})\oplus\CEnd(\F_{L})\to\F_{L}^{\vee}\otimes\E_{L}\otimes
V\\
(g,h)\mapsto(g\otimes id_{V})\circ\A_{L}-\A_{L}\circ h
\end{gather*}
with respect to our choice of basis \eqref{A1} for
$\F_{L}^{\vee}\otimes\E_{L}$.

\medskip

2) For each matrix $M$, determine a $6\times 6$-matrix $D$,
with entries in $V^{\vee}$ and $\CC[u,t]$, such that
\begin{enumerate}
\item[i)] $D\cdot M=0$
\item[ii)] $\rank D=6$ for all $[u,t]\in L$.
\end{enumerate}

\medskip

\noindent{\em Remark.} From
Table~\ref{tabledifferentialmap}, we can read off the
decomposition of $T\NN_{L}$. Since $\deg
f_{1}^{\vee}\otimes  e_{1}=\deg f_{1}^{\vee}\otimes e_{2}=
\deg f_{2}^{\vee}\otimes e_{3}=0$, $\deg f_{1}^{\vee}\otimes
e_{3}=1$, and $\deg f_{2}^{\vee}\otimes e_{1}=
\deg f_{2}^{\vee}\otimes e_{2}=-1$, we find that for lines
1)--2)
\[T\NN_{L}\simeq\oh_{L}(2)\oplus\oh_{L}(1)\oplus 4\oh_{L}\,,\]
and for lines 3)--6)
\[T\NN_{L}\simeq\oh_{L}(2)\oplus 2\oh_{L}(1)\oplus 2\oh_{L}
\oplus\oh_{L}(-1)\,,\]
where the $\oh_{L}(-1)$-component is the image of
$f_{1}^{\vee}\otimes e_{3}\otimes x_{0}$.

More intrinsically, if $\rho([L])=(p,l)$, then there
are canonical isomorphisms (see remark in Section~\ref{modulspace})
\begin{equation*}
\begin{split}
\E_{L}=&H^0(\I_{p}\I_{l}(2))\oplus\oh(-\check{\sigma})\\
\F_{L}\simeq& H^0(\I_{l}(1))\otimes\biggl(\wedge^2
H^0(\I_{p}(1))\oplus
\oh_{L}(-\check{\sigma})\biggr)\,.
\end{split}
\end{equation*}
This implies that the $\oh_{L}(-1)$-component of $T\NN_{L}$
is canonically isomorphic to
\begin{equation*}
H^0(\I_{l}(1))^{\vee}\otimes \wedge^2H^0(\I_{p}(1))^{\vee}
\otimes V/H^0(\I_{p}(1))\otimes
\oh_{L}(-\check{\sigma})\,.
\end{equation*}

\section{Picard-Fuchs operators}
\label{PicardFuchs}

\bigskip

This is an elaboration of sections 4.1--4.3 in \cite{BvS}.

Let $\DD=\CC[D,q]$ and $\DD[t]=\CC[D,q,t]$ be algebras of linear
differential operators acting on $\CC[[q]]$ and $\CC[[q]][t]$, where
$t$ is a formal variable, $q=e^t$, and $D=\frac{\partial}{\partial t}=
q\frac{\partial}{\partial q}$. Hence
$
[D,q]=q$ and $[D,t]=1$.
A differential operator $\P\in\D$ of order $d$ and degree $m$ can be
written in the form
\begin{equation*}
\P=P_{0}(D)+qP_{1}(D)+\cdots+q^mP_{m}(D)
\end{equation*}
where $P_{0},\ldots, P_{m}$ are polynomials of degree $\le d$ with at
least one equality. It is Picard-Fuchs if $\deg P_{0}=d$.

Let
$\psi_{0}=\sum a_{i}q^{i}$ in $\CC[[q]]$.
 If $P$ is a
polynomial then $P(D)q^n=q^nP(n)$, hence the coefficient of $q^n$ in
$\P\psi_{0}$ is
\begin{equation*}
a_{n}P_{0}(n)+a_{n-1}P_{1}(n-1)+\cdots+a_{n-m}P_{m}(n-m)\,.
\end{equation*}
Thus
$\psi_{0}=\sum a_{i}q^{i}$
 is a solution to $\P=0$ if and only if
\\
\begin{equation}
\label{e:e8-1}
a_{n}P_{0}(n)+a_{n-1}P_{1}(n-1)+\cdots+a_{n-m}P_m(n-m)=0
\end{equation}
\\
for all $n$.
Note that the $a_{i}$'s are determined recursively by
\eqref{e:e8-1}
once we have specified all the $a_{n}$ where $P_{0}(n)=0$.

\begin{lem}
$[\P,t]=\P'$ where $\P'$ is the formal derivative of $\P$ with respect
to $D$.
\end{lem}
\begin{proof} See section~4 in \cite{BvS}.
\end{proof}

Next, suppose we are looking for
$\psi_{1}=\sum b_{i}q^{i}$
such that
$t\cdot\psi_{0}+\psi_{1}$ and $\psi_{0}$
are solutions to $\P$.
If $\P\psi_{0}=0$ we have
\[\P(t\cdot\psi_{0}+\psi_{1})=(\P\circ
t)\psi_{0}+\P\psi_{1}=\P'\psi_{0}+\P\psi_{1}\,.\]
Hence $\psi_{1}$ is determined recursively by the relations
\begin{equation}
\label{e:e8-2}
a_{n}P_{0}^{\prime}(n)+\cdots+a_{n-m}P_{m}^{\prime}(n-m)+
b_{n}P_{0}(n)+\cdots+b_{n-m}P_{m}(n-m)=0
\end{equation}
for all $n$ where the $a_{i}\,$'s satisfy \eqref{e:e8-1}.
To get started, we must specify $b_{n}$ for all $n$ such that
$P_{0}(n)=0$.

Finally, suppose we are looking for $\psi_{2}=\sum c_{i}q^i$ such that
$\frac{t^2}{2}\psi_{0}+t\psi_{1}+\psi_{2}$,
$\, \, t\psi_{0}+\psi_{1}\,$ and, $\psi_0$
are solutions.
If $\P\psi_{0}=0$ and $\P(t\psi_{0}+\psi_{1})=0$, then
\\
\begin{align*}
\P&(\frac{t^2}{2}\psi_{0}+t\psi_{1}+\psi_{2})=
\frac{1}{2}(\P\circ t)(t\psi_{0}+\psi_{1})+\frac{1}{2}(\P\circ t)\psi_{1}
+\P\psi_{2}\\
&=\frac{1}{2}\P'(t\psi_{0}+\psi_{1})+\frac{1}{2}(\P\circ t)\psi_{1}
 +\P\psi_{2}\\
&=\frac{1}{2}(\P'\circ t)\psi_{0}+\frac{1}{2}\P'\psi_{1}
   +\frac{1}{2}(\P\circ t)\psi_{1}+\P\psi_{2}\\
&=\frac{1}{2}\P''\psi_{0}+\frac{1}{2}(t\circ \P)\psi_{0}+
\P'\psi_{1}+
  \frac{1}{2}(t\circ \P)\psi_{1}+\P\psi_{2}\\
&=\frac{1}{2}\P''\psi_{0}+\P'\psi_{1}+\P\psi_{2}\,.
\end{align*}
\\
Hence $\psi_{2}$ is determined recursively by the relations
\begin{equation}
\frac{1}{2}\sum_{i=1}^m a_{n-i}P_i^{\prime\prime}(n-i)+
\sum_{i=1}^m b_{n-i}P_i^{\prime}(n-i)+
\sum_{i=1}^m c_{n-i}P_i(n-i)=0
\end{equation}
for all $n$ where the $a_i$'s and $b_i$'s satisfy \eqref{e:e8-1}
and \eqref{e:e8-2}. Again, $c_n$ when $P_0(n)=0$ must be specified.

\section{The {\texttt {MAPLE}} code}

Available from the author upon request.
\addtocontents{toc}{\vspace{0.2cm}}


\begin{thebibliography}{[AAAAAAA]}

\bibitem[AB]{AB}  M.~Atiyah, R.~Bott, {\em The moment map and
equivariant cohomology}, Topology {\bf 23}, 1--28 (1984)

\bibitem[AM]{AM} P.~Aspinwall, D.~Morrison, {\em Topological
field theory and rational curves}, Comm.~Math.~Phys. {\bf
151}, 245--262 (1993)

\bibitem[B]{B} K.~Behrend, {\em Gromov-Witten invariants in
algebraic geometry}, Invent.~Math. {\bf 127}, 601--617 (1996),
{\texttt alg-geom/9601011}

\bibitem[BB]{BB} V.~Batyrev, L.~Borisov, {\em Dual cones and mirror
symmetry for generalized Calabi-Yau manifolds}, in {\em Essays on
Mirror Manifolds, II} (B.~Greene and S.-T.~Yau, eds.) to appear,
{\texttt alg-geom/9402002}

\bibitem[BCKvS]{BCKvS} V.~V.~Batyrev, I.~Ciocan-Fontanine, B.~Kim,
D.~van Straten, {\em Conifold transitions and mirror symmetry for
Calabi-Yau complete intersections in Grassmanians}, preprint
(1997)\\
      {\em Mirror symmetry and toric degenerations of
      partial flag manifolds}, preprint (1997)

\bibitem[BF]{BF} K.~Behrend, B.~Fantechi, {\em The intrinsic
normal cone}, Invent.~Math. (to appear), {\texttt alg-geom/9601010}

\bibitem[Bi]{Bi} A.~Bialynicki-Birula, {\em Some theorems on actions
of algebraic groups}, Annals~of~Math. {\bf 98}, 480--497 (1973)\\
              A.~Bialynicki-Birula, {\em Some properties of the
decomposition of algebraic varieties by action on a torus},
Bull.~Acad.~Polon.~Sci., S\'er. Sci.~Math.~Astronom.~Phys. no.~9,
667--674 (1976)

\bibitem[BM]{BM} K.~Behrend, Y.~Manin, {\em Stacks of stable
maps and Gromov-Witten invariants}, Duke~Math.~J. {\bf 85}, 1--60
(1996)

\bibitem[Bo]{Bo} C.~Borcea, {\em $K3$ surfaces with involution and
mirror pairs of Calabi-Yau manifolds}, in {\em Essays on Mirror
Manifolds, II} (B.~Greene and S.-T.~Yau, eds.) to appear

\bibitem[BvS]{BvS} V.~V.~Batyrev, D.~van Straten, {\em
Generalized hypergeometric functions and rational curves on
Calabi-Yau complete intersections in toric varieties},
Comm.~Math.~Phys. {\bf 168}, 493--533 (1995), {\texttt
alg-geom/9307010}

\bibitem[CdOGP]{CdOGP} P.~Candelas, X.~C.~de la Ossa,
P.~S.~Green, L.~Parkes, {\em A pair of Calabi-Yau manifolds
as an exactly soluble superconformal field theory}, Nucl.~Phys. {\bf
B~359}, 21--74 (1991), and in {\em Essays on Mirror
Manifolds} (S.-T.~Yau, ed., International Press, Hong Kong, 1992) 31--95

\bibitem[CK]{CK} D.~Cox, S.~Katz, {\em Mirror symmetry and algebraic
geometry}, in preparation

\bibitem[D]{D} B.~Dubrovin, {\em The geometry of $2D$
topological field theories}, in {\em Integrable Systems and
Quantum Groups}, Lecture Notes in Mathematics, volume~1620
(Springer-Verlag, New York, Berlin, Heidelberg, 1996) 120--348

\bibitem[DM]{DM} P.~Deligne, D.~Mumford, {\em The
irreducibility of the space of curves of given genus}, Publ. Math.
Inst. Hautes Etudes Sci. {\bf 36}, 75--110 (1969)

\bibitem[Dr]{Dr} J.~M.~Drezet, {\em Fibr\'es exceptionnels et
vari\'et\'es de modules de faisceaux semi-stables sur $\PP_2(C)$},
J.~reine~ang.~Math. {\bf 380}, 14--58 (1987)\\
             J.~M.~Drezet, {\em Cohomologie des vari\'et\'es de
modules de hauteur nulle}, preprint (1987)

\bibitem[E]{E} D.~Eisenbud, {\em Commutative Algebra (with a
view towards algebraic geometry}, Graduate
Texts in Mathematics, volume~15 (Springer-Verlag, New York,
Berlin, Heidelberg, 1997)

\bibitem[EPS]{EPS} G.~Ellingsrud, R.~Piene, S.~A.~Str\o mme,
{\em On the variety of nets of quadrics defining twisted
cubic curves}, in {\em Space Curves}, Lecture Notes in
Mathematics, volume~1266
(F.~Ghione, C.~Peskine, E.~Sernesi, eds., Springer-Verlag, New York,
Berlin, Heidelberg, 1987)

\bibitem[ES1]{ES1}  G.~Ellingsrud, S.~A.~Str\o mme, {\em On
the Chow ring of a geometric quotient}, Annals~of~Math. {\bf
130} 159--187 (1989)

\bibitem[ES2]{ES2} G.~Ellingsrud, S.~A.~Str\o mme, {\em The
number of twisted cubic curves on the general quintic
threefold}, Math.~Scand. {\bf 76} 5--34 (1995)

\bibitem[ES3]{ES3} G.~Ellingsrud, S.~A.~Str\o mme, {\em Bott's
formula and enumerative geometry}, J.~Am.~Math.~Soc. {\bf 9}, 175--193 (1996)

\bibitem[F]{F}  W.~Fulton, {\em Intersection Theory},
(Springer-Verlag, New York, Berlin, Heidelberg, 1984)

\bibitem[FP]{FP} W.~Fulton, R.~Pandharipande, {\em Notes on
stable maps and quantum cohomology}, Proc.~Am.~Math.~Soc. (to appear),
{\texttt alg-geom/9608011}

\bibitem[G]{G} A.~Givental, {\em Homological geometry I: Projective
hypersurfaces}, Selecta~Math., New Series {\bf 1}, No~2, 325--345
(1995)\\
               A.~Givental, {\em Equivariant Gromov-Witten
invariants}, Int.~Math.~Res.~Not. {\bf 13},  613--663 (1996),
{\texttt alg-geom/9603021}\\
               A.~Givental, {\em A mirror theorem for toric complete
intersections},
preprint (1996), {\texttt alg-geom/9701016}

\bibitem[GMP]{GMP} B.~R.~Greene, D.~R.~Morrison, M.~R.~Plessner, {\em
Mirror manifolds in higher dimension}, Comm.~Math.~Phys. {\bf 173},
559--598 (1995)

\bibitem[GP]{GP} L.~G\"ottsche, R.~Pandharipande, {\em The quantum 
cohomology of blow-ups of $\PP^2$ and enumerative geometry}, 
{\texttt alg-geom/9611012} (1996)

\bibitem[GrP]{GrP} T.~Graber, R.~Pandharipande, {\em
Localization of virtual classes}, {\texttt alg-geom/9708001} (1997)

\bibitem[Gr]{Gr} A.~Grothendieck, {\em Sur quelques
propri\'et\'es fondamentales en th\'eorie des
intersections}, S\'eminaire Chevalley 2\`eme ann\'ee.
Secr.~Math.~Paris, 1958

\bibitem[H]{H}  R.~Hartshorne, {\em Algebraic Geometry},
(Springer-Verlag, New York, Berlin, Heidelberg, 1977)


\bibitem[Katz]{Katz} S.~Katz, {\em Rational curves on Calabi-Yau
threefolds}, in {\em Essays on mirror manifolds} (S.-T.~Yau, ed.,
International Press, Hong Kong, 1992) 168--180

\bibitem[Ke]{Ke}  S.~Keel, {\em Intersection theory on moduli spaces of
stable n-pointed curves of genus zero}, Trans.~Am.~Math.~Soc. {\bf 330}
545--574 (1992)

\bibitem[Kim]{Kim} B.~Kim, {\em Quantum hyperplane section
theorem for homogeneous spaces}, preprint (1997)

\bibitem[Kl]{Kl} S.~Kleiman, {\em The transversality of a
general translate}, Compositio~Math. {\bf 38}, 287--297 (1974)

\bibitem[K]{K}  M.~Kontsevich, {\em Enumeration of rational curves via
torus actions}, in {\em The moduli Space of Curves}
(R.~Dijkgraaf, C.~Faber, and G.~van~der~Geer, eds., Birkhauser,
Boston, Basel, Berlin, 1995) 335--368, , {\texttt hep-th/9405035}


\bibitem[KM]{KM} M.~Kontsevitch, Y.~Manin, {\em
Gromov-Witten classes, quantum cohomology and enumerative
geometry}, Comm.~Math.~Phys. {\bf 164}, 525--562 (1994)

\bibitem[Kn]{Kn} F.~Knudsen, {\em Projectivity of the moduli space of
stable curves. II}, Math.~Scand. {\bf 52}, 1225--1265 (1983)

\bibitem[Kr]{Kr} A.~Kresch, {\em The computer program}
{\texttt farsta} in {\em Quantum cohomology at the Mittag-Leffler
institute, 1996-1997, first semester} (P.Aluffi, preprint 1996)
141--147.


\bibitem[KS]{KS} S.~Katz, S.~A.~Str\o mme, {\em
``Schubert'', a Maple package for intersection theory and
enumerative geometry}. Software and documentation available from
the authors, or by anonymous ftp from
{\texttt ftp.math.okstate.edu} (1992)

\bibitem[LiT]{LiT} J.~Li, G.~Tian, {\em Virtual moduli cycles
and Gromov-Witten invariants of algebraic varieties},
{\texttt
alg-geom/9602007} (1996)\\
            J.~Li, G.~Tian, {\em Virtual moduli cycles
and Gromov-Witten invariants of general symplectic manifolds},
{\texttt alg-geom/9608032} (1996)

\bibitem[LT]{LT} A.~Libgober, J.~Teitelbaum, {\em Lines on Calabi-Yau
complete intersections, mirror symmetry, and Picard-Fuchs equations},
Int.~Math.~Res.~Not. {\bf 1}, 29--39 (1993), {\texttt alg-geom/9301001}

\bibitem[M]{M} B.~V.~Char, K.~O.~Geddes, G.~H.~Gonnet,
B.~L.~Leon, M.~B.~Monagan, S.~M.~Watt, {\em Maple~V
reference manual} (Springer-Verlag, New York, Berlin,
Heidelberg, 1991)

\bibitem[M1]{M1} D.~R.~Morrison, {\em Mirror symmetry and rational
curves on quintic threefolds: A guide for mathematicians},
J.~Am.~Math.~Soc. {\bf 6}, 223--247 (1993)\\
      D.~R.~Morrison, {\em Compactifications of moduli spaces inspired
by mirror symmetry}, in {\em Journ\'ees de G\'eom\'etrie Alg\'ebrique
d'Orsay} (Juillet 1992), Ast\'erisque {\bf 218}, Soc.~Math.~France,
Paris, 243--271 (1993)

\bibitem[M2]{M2} D.~R.~Morrison, {\em Mathematical aspects of mirror
symmetry}, {\texttt alg-geom/9609021} (1996)

\bibitem[Me]{Me} P.~Meurer, {\em Gromov-Witten numbers of
rational curves on Calabi-Yau complete intersections in
weighted projective space $\PP(2,1^n)$}, Doctoral
dissertation, Department of Mathematics, University of Bergen,
Norway (1996)

\bibitem[MF]{MF} D.~Mumford, J.~Fogarty, {\em Geometric
invariant theory}, in {\ em Ergebnisse der Matematik und ihrer
Grenz-Gebiete}, volume 34 (Springer-Verlag, New York, Berlin,
Heidelberg, 1982)

\bibitem[P]{P} R.~Pandharipande, {\em Notes on gravitational
 descendents}, personal notes (1996)

\bibitem[R]{R} E.~R\o dland, {\em A non-complete
intersection mirror manifold pair}, in preparation

\bibitem[RT] {RT} Y.~Ruan, G.~Tian, {\em A mathematical theory
of quantum cohomology}, J.~Diff.~Geom. {\bf 42}, 259--367 (1995)

\bibitem[V]{V} A.~Vistoli, {\em Intersection theory on algebraic stacks
and their moduli}, Invent.~Math. {\bf 97}, 613--670 (1989)

\bibitem[Vo1]{Vo1} C.~Voisin, {\em Miroirs et involutions sur les
surfaces $K3$}, in {\em Journ\'ees de G\'eom\'etrie Alg\'ebrique
d'Orsay} (Juillet 1992), Ast\'erisque {\bf 218}, Soc.~Math.~France,
Paris, 273--323 (1993)

\bibitem[Vo2]{Vo2} C.~Voisin, {\em Sym\'etrie miroir},
Panramas et synth\'eses 2, Soc.~Math.~France,
Paris, (1996)\\
 C.~Voisin, {\em Variations of Hodge structure of Calabi-Yau
 threefolds}, lecture notes (1996)



\bibitem[W]{W} E.~Witten, {\em Two-dimensional gravity and
intersection theory on moduli space}, Surveys~in~Diff.~Geom.
{\bf 1}, 243--310
(1991)

\end{thebibliography}
\end{document}